\documentclass[a4paper,10pt]{article}

\def\zibreport{0}
\def\guidelines{0}
\def\expensiveFigures{1}
\def\longtitle{The~SCIP~Optimization~Suite~8.0}
\def\shortfunding{The work for this article has been partly conducted within the
  \emph{Research Campus MODAL} funded by the German Federal Ministry of
  Education and Research (BMBF grant number 05M14ZAM) and has received funding from
  the European Union's Horizon 2020 research and innovation programme under
  grant agreement No 773897. It has also been
  partly supported by
  the German Research Foundation (DFG) within the
  Collaborative Research Center 805, Project A4,
  and the EXPRESS project of the priority program CoSIP (DFG-SPP 1798),
  the German Research Foundation (DFG) within the project HPO-NAVI (project number 391087700).
  }

\makeatletter
\DeclareRobustCommand*{\escapeus}[1]{%
  \begingroup\@activeus\scantokens{#1 }\endgroup}
\begingroup\lccode`\~=`\_\relax
   \lowercase{\endgroup\def\@activeus{\catcode`\_=\active \let~\_}}
\makeatother

\usepackage[utf8]{inputenc}
\usepackage[font=small, labelfont=bf, margin=1cm]{caption}
\usepackage{amsmath,amsfonts,amssymb}
\usepackage{tabularx,booktabs}
\usepackage{url}
\usepackage{rotating}
\usepackage{microtype}

\usepackage[textsize=small,textwidth=3.3cm]{todonotes}

\usepackage{xspace}
\usepackage{etex}
\usepackage{enumitem}
\usepackage{dsfont}
\usepackage{mathtools}
\usepackage[linesnumbered,ruled,vlined]{algorithm2e}

\SetCommentSty{AlgoCommentStyle}
\usepackage{changes}
\usepackage{longtable}
\usepackage{tabularx}
\usepackage{varwidth}
\usepackage{listings}
\usepackage{diagbox}
\usepackage{makecell}

\usepackage{geometry}
\usepackage{amsthm}
\usepackage{multirow}

\usepackage[breaklinks,colorlinks=true,citecolor=black,linkcolor=black,urlcolor=black,hyperfootnotes=false,
  pdftitle={\longtitle},
  pdfsubject={\longtitle}]{hyperref}

\usepackage[capitalize]{cleveref}
\usepackage[numbers]{natbib}
\usepackage{doi}

\usepackage{tikz}
\usepackage{tikz-3dplot}
\usepackage{tikzscale}
\usetikzlibrary{arrows}
\usetikzlibrary{shapes}
\usetikzlibrary{snakes}
\usetikzlibrary{patterns}
\usetikzlibrary{backgrounds,topaths}
\usetikzlibrary{matrix,chains,scopes,positioning,arrows,fit}

\usepackage{pgfplots}
\pgfplotsset{compat=1.15}
\usepackage{subcaption}
\usepackage{textcomp}
\usepgfplotslibrary{groupplots}
\usepgfplotslibrary{units}

\usepackage{graphicx}

\makeatletter
\newcommand*\rel@kern[1]{\kern#1\dimexpr\macc@kerna}
\newcommand*\widebar[1]{%
  \begingroup
  \def\mathaccent##1##2{%
    \rel@kern{0.8}%
    \overline{\rel@kern{-0.8}\macc@nucleus\rel@kern{0.2}}%
    \rel@kern{-0.2}%
  }%
  \macc@depth\@ne
  \let\math@bgroup\@empty \let\math@egroup\macc@set@skewchar
  \mathsurround\z@ \frozen@everymath{\mathgroup\macc@group\relax}%
  \macc@set@skewchar\relax
  \let\mathaccentV\macc@nested@a
  \macc@nested@a\relax111{#1}%
  \endgroup
}
\makeatother

\newcommand{\LP}{{LP}\xspace}
\newcommand{\LPs}{{LPs}\xspace}
\newcommand{\CIP}{{CIP}\xspace}
\newcommand{\CIPs}{{CIPs}\xspace}
\newcommand{\MIP}{{MIP}\xspace}
\newcommand{\MIPs}{{MIPs}\xspace}
\newcommand{\MILP}{{MILP}\xspace}
\newcommand{\MILPs}{{MILPs}\xspace}

\newcommand{\MINLP}{{MINLP}\xspace}
\newcommand{\MINLPs}{{MINLPs}\xspace}

\newcommand{\T}{\top}

\newcommand{\abs}[1]{\lvert{#1}\rvert}

\newcommand{\defi}{\coloneqq}

\newcommand{\linobj}{c}
\newcommand{\nonlinobj}{f}
\newcommand{\linmatrix}{A}
\newcommand{\nonlincons}{g}
\newcommand{\rhs}{b}
\newcommand{\lb}{\ell}
\newcommand{\ub}{u}
\newcommand{\low}[1]{\underline{#1}}
\newcommand{\upp}[1]{\overline{#1}}
\newcommand{\consindex}{\mathcal{M}}
\newcommand{\varindex}{\mathcal{N}}
\newcommand{\intvarindex}{\mathcal{I}}

\newcommand{\N}{\mathds{N}}
\newcommand{\R}{\mathds{R}}

\newcommand{\Z}{\mathds{Z}}
\newcommand{\Rinf}{\ensuremath{\widebar{\mathds{R}}}\xspace}

\usepackage{siunitx}
\newcommand{\cleaninst}{all}
\newcommand{\affected}{affected}
\newcommand{\alloptimal}{{both-solved}\xspace}
\newcommand{\difftimeouts}{{diff-timeouts}\xspace}
\newcommand{\bracket}[2]{[#1,#2]}

\usepackage{xparse}
\ExplSyntaxOn
\def\myround#1{\num{\fp_eval:n {round(#1, 2)}}}
\ExplSyntaxOff

\definecolor{c1}{HTML}{000060}
\definecolor{c2}{HTML}{0000FF}
\definecolor{c3}{HTML}{36648B}
\definecolor{c4}{HTML}{4682B4}
\definecolor{c5}{HTML}{5CACEE}
\definecolor{c6}{HTML}{00FFFF}
\definecolor{c7}{HTML}{008888}
\definecolor{c8}{HTML}{00DD99}
\definecolor{c9}{HTML}{527B10}
\definecolor{c10}{HTML}{7BC618}
\definecolor{c11}{HTML}{33AA00}

\definecolor{scipoldcol}{HTML}{36648B}
\definecolor{scipnewcol}{HTML}{7BC618}

\newcommand{\solver}[1]{\textsc{#1}\xspace}
\newcommand{\scipopt}{\scip Optimization Suite\xspace}
\newcommand{\scipversion}{8.0}

\newcommand{\scipoptv}{\scipopt~\scipversion\xspace}
\newcommand{\scip}{\solver{SCIP}}
\newcommand{\scipv}{\solver{SCIP}~\scipversion\xspace}
\newcommand{\soplex}{\solver{SoPlex}}
\newcommand{\soplexversion}{6.0}
\newcommand{\soplexv}{\solver{SoPlex}~\soplexversion\xspace}
\newcommand{\papilo}{\solver{PaPILO}}
\newcommand{\papiloversion}{2.0}
\newcommand{\papilov}{\solver{PaPILO}~\papiloversion\xspace}
\newcommand{\zimpl}{\solver{Zimpl}} 

\newcommand{\ug}{\solver{UG}}
\newcommand{\presollib}{\solver{PaPILO}}

\newcommand{\gcg}{\solver{GCG}}
\newcommand{\gcgversion}{3.5}

\newcommand{\scipsdp}{\solver{SCIP-SDP}}
\newcommand{\scipsdpversion}{4.0}
\newcommand{\scipsdpv}{\scipsdp~\scipsdpversion\xspace}
\newcommand{\scipjack}{\solver{SCIP-Jack}}

\newcommand{\scipjackversion}{2.0}

\newcommand{\param}[1]{\texttt{#1}\xspace}
\newcommand{\method}[1]{\texttt{#1}\xspace}

\newcommand{\cplex}{\solver{CPLEX}}
\newcommand{\ipopt}{\solver{Ipopt}}
\newcommand{\cppad}{\solver{CppAD}}

\newcommand{\gurobi}{\solver{Gurobi}}

\newcommand{\worhp}{\solver{WORHP}}

\newcommand{\filtersqp}{\solver{FilterSQP}}

\newcommand{\nbsc}[1]{\mbox{#1}\xspace}
\newcommand{\miplib}{\nbsc{MIPLIB}}

\newcommand{\coral}{\nbsc{COR@L}}

\newcommand{\minlplibtwo}{\nbsc{MINLPLib}}

\definecolor{darkgreen}{HTML}{008800}

\newcommand{\fa}{\text{ for all }}

\newcommand{\perm}{\gamma}
\newcommand{\inv}[1]{{#1}^{-1}}

\newcommand{\group}{\Gamma}

\newcommand{\bliss}{\solver{bliss}}

\theoremstyle{plain}

\setlist[itemize]{leftmargin=3.45ex}
\setlist[itemize,1]{label=$-$,itemsep=0ex,topsep=0.9ex}
\setlist[itemize,2]{label=$\cdot$,topsep=0.5ex,leftmargin=2.75ex}
\setlist[enumerate]{leftmargin=3ex,itemsep=0.1ex,parsep=1ex,topsep=0.9ex}

\definecolor{tabcolor}{HTML}{6666AA}
\definecolor{f1}{HTML}{000060}
\definecolor{f2}{HTML}{0000FF}
\definecolor{f3}{HTML}{36648B}
\definecolor{f4}{HTML}{4682B4}
\definecolor{f5}{HTML}{5CACEE}
\definecolor{f6}{HTML}{00FFFF}
\definecolor{f7}{HTML}{00DD99}
\definecolor{f8}{HTML}{008888}
\definecolor{f9}{HTML}{000000}

\usepackage{float}
\newfloat{program}{thp}{lop}
\floatname{program}{Program}
\crefname{program}{program}{programs}

\newcommand{\inputExpensiveFigure}[1]{
\ifthenelse{\expensiveFigures = 1}{\input{#1}}{}
}

\usepackage{titling}
\pretitle{\begin{center}\linespread{1.05}\Large}
\posttitle{\par\end{center}\vskip 0.5em}
\preauthor{\begin{center}\linespread{1.15}\large \lineskip 0.6em\begin{tabular}[t]{c}}
\postauthor{\end{tabular}\par\end{center}\vskip 0.6em}
\predate{\begin{center}}
\postdate{\par\end{center}}
\usepackage[]{titlesec}
\usepackage{etoolbox}
\makeatletter
\patchcmd{\ttlh@hang}{\parindent\z@}{\parindent\z@\leavevmode}{}{}
\patchcmd{\ttlh@hang}{\noindent}{}{}{}
\makeatother
\titleformat{\section}
{\normalfont\large\bfseries}{\thesection}{0.9ex}{}
\titleformat{\subsection}
{\normalfont\normalsize\bfseries}{\thesubsection}{0.9ex}{}
\titleformat{\subsubsection}
{\normalfont\normalsize\upshape}{\thesubsubsection}{0.9ex}{}
\titleformat{\paragraph}[runin]
{\normalfont\normalsize\itshape}{\theparagraph}{1em}{}
\titleformat{\subparagraph}[runin]
{\normalfont\normalsize\itshape}{\theparagraph}{1em}{}
\titlespacing*{\section}     {0pt}{21dd plus 8pt minus 4pt}{10.5dd}
\titlespacing*{\subsection}   {0pt}{21dd plus 8pt minus 4pt}{10.5dd}
\titlespacing*{\subsubsection}{0pt}{19dd plus 8pt minus 4pt}{10.5dd}
\titlespacing*{\paragraph}   {0pt}{13pt plus 8pt minus 4pt}{1em}
\titlespacing*{\subparagraph}   {0pt}{13pt plus 8pt minus 4pt}{1em}

\usepackage[page]{appendix} 

\newcommand{\myand}{$\cdot$\xspace}
\newcommand{\myorcidlink}[1]{\,\href{https://orcid.org/#1}{\raisebox{-0.45ex}{\includegraphics[width=1.8ex]{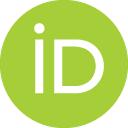}}}}

\newcommand{\scipoptauthors}{%
Ksenia~Bestuzheva\protect\myorcidlink{0000-0002-7018-7099},
Mathieu~Besançon\protect\myorcidlink{0000-0002-6284-3033},
Wei-Kun~Chen\protect\myorcidlink{0000-0003-4147-1346},
Antonia~Chmiela\protect\myorcidlink{0000-0002-4809-2958},
Tim~Donkiewicz\protect\myorcidlink{0000-0002-5721-3563},
Jasper~van~Doornmalen\protect\myorcidlink{0000-0002-2494-0705},
Leon~Eifler\protect\myorcidlink{0000-0003-0245-9344},
Oliver~Gaul\protect\myorcidlink{0000-0002-2131-1911},
Gerald~Gamrath\protect\myorcidlink{0000-0001-6141-5937},
Ambros~Gleixner\protect\myorcidlink{0000-0003-0391-5903},
Leona~Gottwald\protect\myorcidlink{0000-0002-8894-5011},
Christoph~Graczyk,
Katrin~Halbig\protect\myorcidlink{0000-0002-8730-3447},
Alexander Hoen\protect\myorcidlink{0000-0003-1065-1651},
Christopher~Hojny\protect\myorcidlink{0000-0002-5324-8996},
Rolf~van~der~Hulst,
Thorsten Koch\protect\myorcidlink{0000-0002-1967-0077},
Marco~L\"ubbecke\protect\myorcidlink{0000-0002-2635-0522},
Stephen~J.~Maher\protect\myorcidlink{0000-0003-3773-6882},
Frederic Matter\protect\myorcidlink{0000-0002-0499-1820},
Erik~M\"uhmer\protect\myorcidlink{0000-0003-1114-3800}, 
Benjamin~M\"uller\protect\myorcidlink{0000-0002-4463-2873},
Marc~E.~Pfetsch\protect\myorcidlink{0000-0002-0947-7193},
Daniel~Rehfeldt\protect\myorcidlink{0000-0002-2877-074X},
Steffan Schlein,
Franziska~Schl\"osser,
Felipe~Serrano\protect\myorcidlink{0000-0002-7892-3951},
Yuji~Shinano\protect\myorcidlink{0000-0002-2902-882X},
Boro~Sofranac\protect\myorcidlink{0000-0003-2252-9469},
Mark Turner\protect\myorcidlink{0000-0001-7270-1496},
Stefan~Vigerske,
Fabian~Wegscheider,
Philipp~Wellner,
Dieter~Weninger\protect\myorcidlink{0000-0002-1333-8591},
Jakob~Witzig\protect\myorcidlink{0000-0003-2698-0767}
}

\ifthenelse{\zibreport = 1}{
  \let\pdfoutorg\pdfoutput
  \let\pdfoutput\undefined
  \usepackage{zibtitlepage}
  \let\pdfoutput\pdfoutorg

  \ZTPTitle{\longtitle}
  \ZTPAuthor{\scipoptauthors}
  \ZTPPreprint
  \ZTPNumber{21-?}
  \ZTPMonth{?}
  \ZTPYear{2021}
  \ZTPInfo{}
}{}

\begin{document}

\ifthenelse{\guidelines = 1}{\input{guidelines}}{}

\title{\longtitle}

\author{%
  Ksenia Bestuzheva\protect\myorcidlink{0000-0002-7018-7099} \myand
  Mathieu Besançon\protect\myorcidlink{0000-0002-6284-3033} \myand
  Wei-Kun Chen \protect\myorcidlink{0000-0003-4147-1346} \and
  Antonia Chmiela\protect\myorcidlink{0000-0002-4809-2958} \myand
  Tim Donkiewicz\protect\myorcidlink{0000-0002-5721-3563} \myand
  Jasper van Doornmalen\protect\myorcidlink{0000-0002-2494-0705} \and
  Leon Eifler\protect\myorcidlink{0000-0003-0245-9344} \myand
  Oliver Gaul\protect\myorcidlink{0000-0002-2131-1911} \myand
  Gerald Gamrath\protect\myorcidlink{0000-0001-6141-5937} \and
  Ambros Gleixner\protect\myorcidlink{0000-0003-0391-5903} \myand
  Leona Gottwald\protect\myorcidlink{0000-0002-8894-5011} \myand
  Christoph Graczyk \and
  Katrin Halbig\protect\myorcidlink{0000-0002-8730-3447} \myand
  Alexander Hoen\protect\myorcidlink{0000-0003-1065-1651} \myand
  Christopher Hojny\protect\myorcidlink{0000-0002-5324-8996} \and
  Rolf van der Hulst \myand
  Thorsten Koch\protect\myorcidlink{0000-0002-1967-0077} \myand
  Marco L\"ubbecke\protect\myorcidlink{0000-0002-2635-0522} \and
  Stephen J.~Maher\protect\myorcidlink{0000-0003-3773-6882} \myand
  Frederic Matter\protect\myorcidlink{0000-0002-0499-1820} \myand
  Erik M\"uhmer\protect\myorcidlink{0000-0003-1114-3800} \and
  Benjamin M\"uller\protect\myorcidlink{0000-0002-4463-2873} \myand
  Marc E. Pfetsch\protect\myorcidlink{0000-0002-0947-7193} \myand
  Daniel Rehfeldt\protect\myorcidlink{0000-0002-2877-074X} \and
  Steffan Schlein \myand
  Franziska Schlösser \myand
  Felipe Serrano\protect\myorcidlink{0000-0002-7892-3951} \and
  Yuji Shinano\protect\myorcidlink{0000-0002-2902-882X} \myand
  Boro Sofranac\protect\myorcidlink{0000-0003-2252-9469} \myand
  Mark Turner\protect\myorcidlink{0000-0001-7270-1496} \and
  Stefan Vigerske \myand
  Fabian Wegscheider \myand
  Philipp Wellner \and
  Dieter Weninger\protect\myorcidlink{0000-0002-1333-8591} \myand
  Jakob Witzig\protect\myorcidlink{0000-0003-2698-0767}%
  \thanks{Extended author information is available at the end of the paper.
    \shortfunding}}


\ifthenelse{\zibreport = 1}{\zibtitlepage}{}

\newgeometry{left=38mm,right=38mm,top=35mm}

\maketitle

\paragraph{\bf Abstract}


The \scipopt provides a collection of software packages for mathematical
optimization centered around the constraint integer programming framework \scip.
This paper discusses enhancements and extensions contained in version~8.0 of the
\scipopt.
Major updates in SCIP include improvements in symmetry handling and decomposition algorithms,
new cutting planes, a new plugin type for cut selection, and a complete rework of the way
nonlinear constraints are handled.
Additionally, \scipv now supports interfaces for Julia as well as Matlab.
Further, UG now includes a unified framework to parallelize all solvers, a utility to analyze
computational experiments has been added to GCG, dual solutions can be postsolved by PaPILO,
new heuristics and presolving methods were added to SCIP-SDP, and additional problem classes and
major performance improvements are available in SCIP-Jack.



\newpage

\section{Introduction}
\label{sect:introduction}
The \scipopt comprises a set of complementary software packages designed to model and
solve a large variety of mathematical optimization problems:
%
\begin{itemize}
\item the modeling language \zimpl~\cite{Koch2004},
\item the presolving library \presollib for linear and mixed-integer linear programs, a new addition in version~7.0 of the \scipopt~\cite{papilogithub},
\item the simplex-based linear programming solver
  \soplex~\cite{Wunderling1996},
\item the constraint integer programming solver
  \scip~\cite{Achterberg2009}, which can be used as a fast standalone solver for
  mixed-integer linear and nonlinear programs and a flexible
  branch-cut-and-price framework,
\item the automatic decomposition solver \gcg~\cite{GamrathLuebbecke2010}, and
\item the \ug framework for parallelization of branch-and-bound
  solvers~\cite{Shinano2018}.
\end{itemize}
All six tools can be downloaded in source code and are freely available
for members of noncommercial and academic institutions. 
They are accompanied by several extensions for solving specific problem-classes
such as the award-winning Steiner tree solver
\scipjack~\cite{Gamrath2017scipjack} and the mixed-integer semidefinite
programming solver \scipsdp~\cite{GallyPfetschUlbrich2018}.
This paper describes the new features and enhanced algorithmic components
contained in version~8.0 of the \scipopt.


\paragraph{Background}

\scip has been designed as a branch-cut-and-price framework to
solve different types of optimization problems, most importantly,
\emph{mixed-integer linear programs} (\MILPs) and
\emph{mixed-integer nonlinear programs} (\MINLPs).
\MILPs are optimization problems of the form
\begin{equation}
  \begin{aligned}
    \min \quad& \linobj^\T x \\
    \text{s.t.} \quad& \linmatrix x \geq \rhs, \\
    &\lb_{i} \leq x_{i} \leq \ub_{i} && \fa i \in \varindex, \\
    &x_{i} \in \Z && \fa i \in \intvarindex,
  \end{aligned}
  \label{eq:generalmip}
\end{equation}
defined by $c \in \R^n$, $A \in\R^{m\times n}$, $ \rhs\in \R^{m}$, $\lb$, $\ub \in
\Rinf^{n}$, and the index set of integer variables $\mathcal{I} \subseteq \mathcal{N} \defi \{1, \ldots, n\}$.  The usage of $\Rinf \defi \R \cup
\{-\infty,\infty\}$ allows for variables that are free or bounded only in
one direction (we assume that variables are not fixed to~$\pm \infty$).

Another focus of \scip's research and development are \emph{mixed-integer nonlinear
  programs} (\MINLPs).
\MINLPs can be written in the form
\begin{equation}
  \begin{aligned}
    \min \quad& \nonlinobj(x) \\
    \text{s.t.} \quad& \nonlincons_{k}(x) \leq 0 && \fa k \in \consindex, \\
    &\lb_{i} \leq x_{i} \leq \ub_{i} && \fa i \in \varindex, \\
    &x_{i} \in \Z && \fa i \in \intvarindex,
  \end{aligned}
  \label{eq:generalminlp}
\end{equation}
where the functions $\nonlinobj : \R^n \rightarrow \R$ and $\nonlincons_{k} : \R^{n}
\rightarrow \R$, $k \in \consindex \defi \{1,\ldots,m\}$, are possibly nonconvex.
Within \scip, we assume that $\nonlinobj$ and $\nonlincons_k$ are specified explicitly in
algebraic form using base expressions that are known to \scip.

\scip is not restricted to solving \MILPs and \MINLPs, but is a framework for solving \emph{constraint integer programs} (\CIPs),
a generalization of the former two problem classes.
The introduction of \CIPs was motivated by the modeling flexibility of
constraint programming and the algorithmic requirements of integrating it with
efficient solution techniques available for \MILPs. Later on, this framework
allowed for an integration of \MINLPs as well.
Roughly speaking, \CIPs are finite-dimensional optimization problems with arbitrary constraints and a linear objective function
that satisfy the following property: if all integer variables are fixed, the remaining subproblem must form a linear or nonlinear program.

In order to solve \CIPs, \scip constructs relaxations---typically \LP relaxations.
If the relaxation solution is not feasible for the current subproblem, the enforcement callbacks of the constraint handlers
need to take measures to eventually render the relaxation solution infeasible for the updated relaxation,
for example by branching or separation.
%
Being a framework for solving \CIPs, \scip can be extended by plugins to be able to solve any \CIP.
The default plugins included in the \scipopt provide tools to solve \MILPs and many \MINLPs as
well as some classes of instances from constraint programming, satisfiability testing, and
pseudo-Boolean optimization. Additionally, SCIP-SDP allows to
solve mixed-integer semidefinite programs.


The core of \scip coordinates a central branch-cut-and-price algorithm.
The methods for processing constraints of a given type are implemented in a corresponding constraint handler, and
advanced methods like primal heuristics, branching rules, and cutting plane separators can be integrated as plugins with a pre-defined interface.
\scip comes with many such plugins needed to achieve a good \MILP and \MINLP performance.
In addition to plugins supplied as part of the SCIP distribution, new plugins can be created by users.
This basic design and solving process is described in more detail by
Achterberg~\cite{Achterberg2007a}.

By design, \scip interacts closely with the other components of the \scipopt.
Optimization models formulated in \zimpl can be read by \scip.
\presollib provides an additional fast and effective presolving procedure that is called from a \scip presolver plugin.
The linear programs (\LPs) solved repeatedly during the branch-cut-and-price algorithm
are by default optimized with \soplex.
Interfaces to several external \LP solvers exist, and new interfaces can be added by users.
\gcg extends \scip to automatically detect problem structure
and generically apply decomposition algorithms based on the Dantzig-Wolfe or the
Benders' decomposition scheme.
And finally, the default instantiations of the \ug framework use \scip as a base
solver in order to perform branch-and-bound in parallel computing
environments with shared or distributed memory architectures.

\paragraph{New Developments and Structure of the Paper}

This paper focuses on two main aspects.
The first one is to explain the changes and progress made in the solving process of \scip
and analyze the resulting improvements on \MILP and \MINLP instances, both in terms of performance and robustness.
A performance comparison of \scipv against \scip 7.0 is carried out in Section \ref{sect:performance}.
Improvements to the core of \scip are presented in Section~\ref{sect:scip} and include
\begin{itemize}
   \item a new framework for handling nonlinear constraints,
   \item symmetry handling on general variables and improved orbitope detection,
   \item a new separator for mixing cuts,
   \item improvements to decomposition-based heuristics and the Benders decomposition framework, and
   \item a new type of plugins for cut selection, and several technical improvements.
\end{itemize}
A more detailed explanation of the changes to the \MINLP solving process and the new expression framework is given in
Section \ref{sect:minlp}.
Improvements to the default \LP solver \soplex and presolver \papilo are explained in Section \ref{sect:soplex}
and \ref{sect:papilo} respectively.
This aspect will be of interest to the optimization community working on methods and algorithms related to these building
blocks and to practitioners willing to understand the performance they observe on their particular instances.

The second aspect of this paper is to present the evolving possibilities for working with the \scipoptv for optimization practitioners.
This includes improvements and changes to the interfaces in Section \ref{sect:interfaces} and the modeling language \zimpl in Section~\ref{sect:zimpl};
to \scip extensions specialized for other computational settings
such as distributed computing with UG in Section \ref{sect:ug} and Dantzig-Wolfe decompositions with GCG in Section \ref{sect:gcg};
and finally to \scip extensions for particular problem classes such as the mixed-integer semidefinite solver SCIP-SDP in Section \ref{sect:SCIP-SDP} and the Steiner tree solver SCIP-Jack in Section \ref{sect:SCIP-Jack}.

\section{Overall Performance Improvements for MILP and MINLP}
\label{sect:performance}


\scip is used extensively to solve mixed-integer linear and nonlinear programs out of the box.
In this section, we present computational experiments conducted by running \scip without parameter tuning or
algorithmic variations to assess the performance changes since the 7.0.0 release.
We detail below the methodology and results of these experiments.

The indicators of interest to compare the two versions of \scip
on a given subset of instances are
the number of solved instances, the shifted geometric mean
of the number of branch-and-bound nodes, and the shifted geometric mean
of the solving time.
The \emph{shifted geometric mean} of values $t_1, \dots, t_n$ is
\[
\big(\prod_{i=1}^n(t_i + s)\big)^{1/n} - s.
\]
The shift~$s$ is set to 100~nodes and 1~second, respectively.

\subsection{Experimental Setup}

We use the \scipopt 7.0.0 as the baseline, including \soplex 5.0.0 and \papilo 1.0.0,
and compare it with \scipoptv including \soplexv and \papilov.
Both were compiled using GCC 7.5, use \solver{Ipopt} 3.12.13 as NLP subsolver built with the \solver{MUMPS} 4.10.0
numerical linear algebra solver, \solver{CppAD} 20180000.0 as algorithmic differentiation library,
and \solver{bliss} 0.73 for graph automorphisms to detect symmetry in \MIPs.
The time limit was set to 7200 seconds in all cases.

The \MILP instances are selected from the \miplib 2003, 2010, and 2017~\cite{miplib2017} and
\coral instance sets, including all instances
previously solved by \scip 7.0.0 with at least one of five random seeds or newly solved by \scipv with at least one of five random seeds;
this amounts to 347~instances.
The \MINLP instances are similarly selected from the \minlplibtwo\footnote{\url{https://www.minlplib.org}}
with newly solvable instances added to the ones previously solved by \scip 7 for a total of 113 instances.

All performance runs are carried out on identical machines with Intel Xeon CPUs E5-2690 v4 @ 2.60GHz
and 128GB in RAM. A single run is carried out on each machine in a single-threaded mode.
Each optimization problem is solved with SCIP using five different seeds for random number generators.
This results in a testset of 565 \MINLPs and 1735 \MILPs.
Instances for which the solver reported numerically inconsistent results are excluded from the presented results.

\subsection{\MILP Performance}

The results of the performance runs on \MILP instances are presented in \cref{tbl:rubberband_table_mip}.
The changes introduced with
\scip 8.0 improved the performance on \MILPs both in terms of number of solved instances
and shifted geometric mean of the time.
Furthermore, the difference in terms of geometric mean time is starker on harder instances, with an improvement of up to $52\%$ on instances taking more than 1000 seconds to solve.
The improvement is more limited on \alloptimal instances that were solved by both solvers, for which the relative improvement is only of $11\%$.
This indicates that the overall speedup is due to newly solved instances more than to improvement on instances that were already solved by \scip 7.0.

\begin{table}
\caption{Performance comparison for MILP instances}
\label{tbl:rubberband_table_mip}
\scriptsize

\begin{tabular*}{\textwidth}{@{}l@{\;\;\extracolsep{\fill}}rrrrrrrrr@{}}
\toprule
&           & \multicolumn{3}{c}{\scip~8.0.0+\soplex~6.0.0} & \multicolumn{3}{c}{\scip~7.0.0+\soplex~5.0.0} & \multicolumn{2}{c}{relative} \\
\cmidrule{3-5} \cmidrule{6-8} \cmidrule{9-10}
Subset                & instances &                                   solved &       time &        nodes &                                   solved &       time &        nodes &       time &        nodes \\
\midrule
\cleaninst            &      1708 &                                     1478 &      231.3 &         3311 &                                     1445 &      271.3 &         4107 &        1.17 &          1.24 \\
\affected             &      1475 &                                     1424 &      173.8 &         2843 &                                     1391 &      209.7 &         3611 &        1.21 &          1.27 \\
\cmidrule{1-10}
\bracket{0}{tilim}    &      1529 &                                     1478 &      154.4 &         2512 &                                     1445 &      184.6 &         3167 &        1.20 &          1.26 \\
\bracket{1}{tilim}    &      1470 &                                     1419 &      185.9 &         2870 &                                     1386 &      223.8 &         3647 &        1.20 &          1.27 \\
\bracket{10}{tilim}   &      1361 &                                     1310 &      248.1 &         3612 &                                     1277 &      303.1 &         4661 &        1.22 &          1.29 \\
\bracket{100}{tilim}  &      1000 &                                      949 &      537.1 &         7270 &                                      916 &      702.6 &        10262 &        1.31 &          1.41 \\
\bracket{1000}{tilim} &       437 &                                      386 &     1566.2 &        17973 &                                      353 &     2383.1 &        31707 &        1.52 &          1.76 \\
\difftimeouts         &       135 &                                       84 &     2072.7 &        19597 &                                       51 &     5062.1 &        69354 &        2.44 &          3.54 \\
\alloptimal           &      1394 &                                     1394 &      119.9 &         2048 &                                     1394 &      133.8 &         2330 &        1.12 &          1.14 \\
\bottomrule
\end{tabular*}
\end{table}


\subsection{\MINLP Performance}
\label{sect:perfminlp}

With the major revision of the handling of nonlinear constraints,
the performance of \scip on \MINLPs has changed a lot on the instance set compared to \scip~7.0.
The \MINLP performance results are summarized in \cref{tbl:rubberband_table_minlp}.
On all subsets of the instances selected by runtime, more instances are solved by \scipv than by \scip~7.0.
Furthermore, \scipv solves the instances for each of these subsets with a shorter shifted geometric mean time
even though it produces more nodes in the branch-and-bound tree.

On the 382 instances solved by both versions, \scipv requires fewer nodes and less time.
The number of instances solved by only one of the two versions (\difftimeouts) is much higher than
reported in previous release reports with similar experiments, with 66 instances newly solved by \scipv
and 46 instances previously solved that \scipv did not succeed on.


A finer comparison of the two \scip versions on additional subsets of instances is provided in
\cref{tab:detailedminlp}.
Instances are split according to mixed-integer and continuous, nonconvex and convex problems.
They are classified as mixed-integer if at least one integer or binary constraint is present
in the original problem. The convexity of instances is identical to the information provided on the \minlplibtwo website.

\cref{tab:detailedminlp} shows that \scipv brings most significant improvements for nonconvex problems with 41 more instances solved and a drastic speedup factor of 3.54 on the purely continuous nonconvex problems.
Performance has, however, degraded on convex problems with 21 instances that are not solved anymore
and the shifted geometric mean runtime more than tripled.

\begin{table}
\caption{Performance comparison for MINLP}
\label{tbl:rubberband_table_minlp}
\scriptsize

\begin{tabular*}{\textwidth}{@{}l@{\;\;\extracolsep{\fill}}rrrrrrrrr@{}}
\toprule
&           & \multicolumn{3}{c}{\scip~8.0.0+\soplex~6.0.0} & \multicolumn{3}{c}{\scip~7.0.0+\soplex~5.0.0} & \multicolumn{2}{c}{relative} \\
\cmidrule{3-5} \cmidrule{6-8} \cmidrule{9-10}
Subset                & instances &                                   solved &       time &        nodes &                                   solved &       time &        nodes &       time &        nodes \\
\midrule
\cleaninst            &       558 &                                      454 &       39.1 &         2427 &                                      435 &       45.7 &         1845 &        1.17 &          0.76 \\
\affected             &       487 &                                      438 &       23.5 &         1748 &                                      419 &       28.4 &         1456 &        1.21 &          0.83 \\
\cmidrule{1-10}
\bracket{0}{tilim}    &       503 &                                      454 &       21.7 &         1585 &                                      435 &       25.9 &         1326 &        1.19 &          0.84 \\
\bracket{1}{tilim}    &       375 &                                      326 &       56.1 &         3994 &                                      307 &       71.0 &         3113 &        1.27 &          0.78 \\
\bracket{10}{tilim}   &       293 &                                      244 &      121.6 &         7450 &                                      225 &      169.3 &         5393 &        1.39 &          0.72 \\
\bracket{100}{tilim}  &       195 &                                      146 &      307.6 &        14204 &                                      127 &      433.9 &         6696 &        1.41 &          0.47 \\
\bracket{1000}{tilim} &       153 &                                      104 &      466.9 &        23425 &                                       85 &      565.3 &         8382 &        1.21 &          0.36 \\
\difftimeouts         &       117 &                                       68 &      451.4 &        29142 &                                       49 &      461.8 &         6275 &        1.02 &          0.22 \\
\alloptimal           &       386 &                                      386 &        8.2 &          609 &                                      386 &       10.4 &          806 &        1.27 &          1.32 \\
\bottomrule
\end{tabular*}
\end{table}


As can be seen in the table, this is mostly due to worse performance on a specific group of instances, the \texttt{syn} group of instances.
The \texttt{syn} group includes specifically instances \texttt{syn40m04h}, \texttt{rsyn0840m03h}, \texttt{rsyn0820m02m}, \texttt{syn20h}, \texttt{syn30m03h}, and \texttt{rsyn0840m04h}.
The solving time for \texttt{syn} instances has degraded significantly on \scipv while
the degradation is moderate on other convex instances.
The much higher time on the \texttt{syn} instances explains alone the degradation on the total convex subset.
We presume that the new expression simplification at the moment obfuscates some structure on some instances of the \texttt{syn} group that was exploited with \scip~7.0.

An MINLP performance evaluation that focuses only on the changes in handling nonlinear constraints is given in Section~\ref{sect:perfconsexpr}.

\begin{table}
\caption{Detailed MINLP performance comparison}
\label{tab:detailedminlp}
\centering
\scriptsize
\begin{tabular*}{\textwidth}{@{}ll@{\;\;\extracolsep{\fill}}rrrrrr@{}}
\toprule
\multicolumn{2}{c}{Subset} &   & \multicolumn{2}{c}{\scip~8+\soplex~6} & \multicolumn{2}{c}{\scip~7+\soplex~5} & relative \\
\cmidrule{1-2} \cmidrule{4-5} \cmidrule{6-7} \cmidrule{8-8}
convexity & integrality & instances & solved & time & solved & time & time\\
\midrule
nonconvex & continuous & 95 & 88 & 66.60 & 71 & 235.86 & 3.54 \\
nonconvex & mixed-integer & 320 & 258 & 236.43 & 244 & 315.74 & 1.34 \\
nonconvex & total & 415 & 346 & 186.43 & 315 & 295.92 & 1.59 \\
\midrule
convex & continuous & 5 & 5 & 0.14 & 5 & 0.06 & 0.45 \\
convex & mixed-integer & 160 & 118 & 329.37 & 134 & 167.62 & 0.51 \\
convex-no-syn & mixed-integer & 130 & 108 & 209.04 & 108 & 190.70 & 0.91 \\
convex-syn-only & mixed-integer & 30 & 10 & 1685.14 & 26 & 87.01 & 0.05 \\
convex & total & 165 & 123 & 310.84 & 139 & 159.76 & 0.51 \\
\bottomrule
\end{tabular*}
\end{table}

\section{SCIP}
\label{sect:scip}

\subsection{A New MINLP Framework}\label{sect:minlp_summary}

The \scip~8.0 release comes with a major change in the way that nonlinear constraints are handled.
The main motivation for this change is twofold: First, it aims at increasing the reliability of the
solver and alleviating a numerical issue that arose from problem reformulations and led to SCIP
returning solutions that are feasible in the reformulated problem, but infeasible in the original
problem.
Second, the new design of the nonlinear framework reduces the ambiguity of expression and structure
types by implementing different plugin types for low-level expression types that define expressions,
and high-level structure types that add functionality for particular, often overlapping structures.
Finally, a number of new features for improving the solver's performance on MINLPs were introduced.
A detailed description of the changes can be found in Sections~\ref{sect:minlp} and~\ref{sect:symmetrynl}.

\subsection{Improvements in Symmetry Handling}
\label{sec:symmetry}

Symmetries are well-known to have an adverse effect on the performance of
\MILP and \MINLP solvers, because symmetric subproblems are treated
repeatedly without providing new information to the solver.
For this reason, there exist different methods to handle symmetries in
\scip.
Until version~7.0, \scip was only able to handle symmetries in \MILP{}s.
With the release of~\scip~8.0 also symmetries in \MINLP{}s can be handled.
Furthermore, the release of
\scip~8.0 features several algorithmic enhancements of existing
as well as the implementation of further symmetry handling methods.
In the following, we describe the kind of symmetries \scip can handle and
list the techniques used in \scip~7.0.
Afterwards, we describe the novel symmetry handling methods and highlight
algorithmic enhancements.

Let us start with some preliminary remarks.
For a permutation~$\perm$ of the variable index set~$\{1,\dots,n\}$ and a
vector~$x \in \R^n$, we define~$\perm(x) = (x_{\inv\perm(1)}, \dots,
x_{\inv\perm(n)})$.
We say that~$\perm$ is a \emph{symmetry} of~\eqref{eq:minlp} if the
following holds: $x \in \R^n$
is feasible for~\eqref{eq:minlp} if and only if~$\perm(x)$ is feasible, and~$c^\T x =
c^\T\perm(x)$.
The set of all symmetries forms a group~$\bar{\group}$, the \emph{symmetry group}
of~\eqref{eq:minlp}.
Since computing~$\bar{\group}$ is NP-hard, see Margot~\cite{Margot2010}, one typically
refrains from handling all symmetries.
Instead, one only computes a subgroup~$\group$ of~$\bar{\group}$ that keeps
the constraint system of~\eqref{eq:minlp} invariant.
Computing this \emph{formulation group}~$\group$ for \MILP{}s can be
accomplished by computing symmetries of an auxiliary graph, see Salvagnin~\cite{Salvagnin2005}.
In \scip~8.0, the already existing routine for computing symmetries of
\MILP{}s has been extended to handle also nonlinear constraints.
To detect symmetries of the auxiliary graphs, \scip uses the graph
isomorphism package \solver{bliss}~\cite{bliss}.

\subsubsection{Previously Existing Symmetry Handling Methods in \scip}

\scip~7.0 used two paradigms to handle symmetries of binary variables:
a constraint-based approach or the pure propagation-based approach \emph{orbital
  fixing}~\cite{Margot2002,Margot2003,OstrowskiEtAl2011}.
The constraint-based approach is implemented via three different constraint
handler plugins to deal with different kinds of matrix symmetries.
The \emph{symresack} constraint handler~\cite{HojnyPfetsch2019} provides separation and
propagation routines for general permutations~$\perm$, whereas the \emph{orbisack}
constraint handler~\cite{KaibelLoos2011} uses specialized separation and propagation
methods if~$\perm$ is a composition of 2-cycles.
The \emph{orbitope} constraint handler~\cite{BendottiEtAl2021,HojnyPfetsch2019} handles symmetries of special
subgroups of~$\group$.
These subgroups are required to act on binary matrices and to be able to
reorder their columns arbitrarily.
Moreover, if the variables affected by the corresponding permutations or
groups interact with set packing or partitioning constraints in a certain
way, all constraint handlers provide specialized separation and propagation
mechanisms to find stronger cutting planes and reductions~\cite{Hojny2020,KaibelPfetsch2008,KaibelEtAl2011}.
The common ground of these constraint handlers is that they enforce
solutions that are lexicographically maximal in their orbit of
symmetric solutions.

The integer parameter \param{misc/usesymmetry} can be used to
enable/disable these two methods.
In \scip~7.0, the parameter ranged between~0 and~3, where the~Bit~1
enables/disables the usage of the constraint-based approach and Bit~2 enables/disables
orbital fixing.
If the group~$\group$ is a product group~$\group = \group_1 \otimes \dots
\otimes \group_k$, the variables affected by one factor of~$\group$
are not affected by any other factor.
In this case, one can apply different symmetry
handling methods for the different factors.
The sets of variables affected by the different factors are called the
\emph{components} of~$\group$.
Thus, if both methods are enabled, \scip searches for independent components of
the symmetry group~$\group$ and, depending on structural properties of the
component, either uses cutting planes or orbital fixing:
if a component can be completely handled by orbitopes, \scip uses orbitopes
and orbital fixing otherwise; see the SCIP Optimization Suit 7.0 release report~\cite{SCIP7} for further details.

\subsubsection{Symmetry Detection Extended to Nonlinear Constraints}
\label{sect:symmetrynl}

In \scip~8.0, symmetry detection has been extended to handle two types of symmetries in nonlinear constraints.

The detection of permutation-based symmetries is performed by analyzing expression graphs as first
proposed by Liberti~\cite{Liberti2012a}.
The detected automorphisms are then projected onto problem variables, which yields a
permutation group.

Symmetries of a different type, referred to as complementary symmetries, are detected in quadratic
problems by considering affine transformations $\gamma: \mathbb{R}^b \rightarrow \mathbb{R}^n$,
$\gamma(x) = Rx + s$, where $R \in \mathbb{R}^{n\times n}$, $s \in \mathbb{R}^n$, and $R_{ij}$ may
have the values $1$, $-1$ or $0$ and $s_i$ may be equal to either some constant $d_i$, or $0$.
Such a transformation defines a complementary symmetry if
it preserves the objective and constraint functions.
The detection is performed by solving an auxiliary problem that compares coefficients before and
after substituting variables for their complements.

For more details, see Wegscheider~\cite{Wegscheider2019}.

\subsubsection{New Symmetry Handling Methods}

One drawback of the mentioned approaches is that they can only handle
symmetries of binary variables, but not of general integer or continuous
variables.
Moreover, \scip~7.0 can only detect orbitopes when a component of~$\group$
can be completely handled by orbitopes, but not if some part of the
component allows applying orbitopes.
In \scip~8.0, both issues are resolved by the implementation of further
symmetry handling methods and a refined detection and handling mechanism
for orbitopes.

\paragraph{Symmetry Handling of General Variables}

To handle symmetries of general variables, we have implemented symmetry
handling inequalities derived from the Schreier-Sims table (SST cuts) as
described by Salvagnin~\cite{Salvagnin2018}; see also Liberti and Ostrowski~\cite{LibertiOstrowski2014} for a more
general version.
These inequalities are defined using the following procedure.
Let~$\group$ be the symmetry group of~\eqref{eq:minlp} and let~$A(\group) = \{ i \in
\{1,\dots, n\} : \exists\perm\in\group \text{ with } \perm(i) \neq i\}$ be
its set of affected variables.
Select a variable index~$\ell \in A(\group)$ and compute its orbit~$O =
\{\perm(\ell) : \perm \in \group\}$.
We call~$\ell$ the \emph{leader} of its orbit.
Afterwards, replace the initial group~$\group$ by the stabilizer group~$\{
\perm \in \group : \perm(\ell) = \ell\}$, compute the set of affected
variables, and iterate the procedure of selecting a leader and replacing
groups by their stabilizer until the set of affected variables becomes
empty.
At termination, we are given a list of leaders~$\ell_1,\dots,\ell_k$ with
associated orbits~$O_1,\dots,O_k$.
Salvagnin~\cite{Salvagnin2018} shows that the inequalities
\begin{align*}
  x_{\ell_i} &\geq x_j, && j \in O_i,\; i \in \{1,\dots, k\},
\end{align*}
can be used to handle symmetries of general variables.

The above procedure allows to select the orbit leaders arbitrarily.
In \scip~8.0, several rules exist to select the
leader using the parameters \param{propagating/symmetry/<Xyz>}, where
\param{Xyz} is one of the following: \param{sstleadervartype}, \param{sstleaderule},
\param{ssttiebreakrule}, and \param{sstmixedcomponents}.
The bitset \param{sstleadervartype} controls which (combinations of)
variable types can be used as leaders; if several variable types are
allowed, \scip selects the one with the most affected variables.
The rule \param{sstleaderule} determines whether the first or last variable
(according to SCIP's variable ordering) in an orbit shall be used as
leader.
If a binary variable shall be used as leader, the parameter also allows to
use the number of other binary variables it is in conflict with as a
selection criterion.
We say that two binary variables are in conflict if not both can be~1
simultaneously.
The rule \param{ssttiebreakrule} selects a leader whose orbit is as small or as
large as possible, or that contains the most conflicting binary variables.
Finally, \param{sstmixedcomponents} controls whether SST cuts are allowed
to be added to components of~$\group$ that contain variables of different
types.
By default, we allow adding SST cuts for non-binary variables whose orbit
is as small as possible.
We select the first variable per orbit as leader and allow different
variable types in a component.

Since SST cuts extend the class of previous symmetry handling methods, the
range of parameter \param{misc/usesymmetry} has been extended to~$\{0,\dots,7\}$,
where Bit~4 controls whether SST cuts are enabled.
This in particular means that SST cuts can be used in combination with
other symmetry handling methods.
Below we will describe how SST cuts are used when also other symmetry
handling methods are active.

\paragraph{Improved Orbitope Detection}

As mentioned previously, \scip~7.0 uses orbitopes for a component
of~$\group$ only if all permutations within the component form an
orbitope structure.
This, however, can be rather restrictive as illustrated next.
Consider the problem of coloring an undirected graph~$G = (V,E)$ with~$k$
colors.
Every feasible coloring can be encoded by a matrix~$X \in \{0,1\}^{V \times
  k}$, where~$X_{vi} = 1$ if and only if node~$v$ is colored by color~$i$.
We can transfer the coloring~$X$ into another equivalent coloring~$Y$ by taking
an arbitrary permutation~$\pi$ of~$\{1,\dots,k\}$ and defining~$Y_{vi} =
X_{v\pi(i)}$.
That is, the symmetry group~$\group$ of the coloring problem can reorder the columns
of binary matrices arbitrarily, and thus, allows the application of
orbitopes as indicated above.
If the graph~$G$ is symmetric, however, $\group$ will also contain
permutations that reorder the rows of~$X$ according to automorphisms
of~$G$.
Since these row permutations interact with the variables affected by column
permutations, they form a common component.
Hence, not all permutations within this component are permutations
necessary for an orbitope and the detection routine of \scip~7.0 will not
recognize the applicability of orbitopes.

In \scip~8.0 we have refined the orbitope detection routine to be able to
heuristically find such hidden orbitopes.
In the following, we call such orbitopes \emph{suborbitopes}, because they are
defined by a subgroup of a component.
To explain the procedure, note that a subgroup of a component defines an
orbitope for matrix~$X$ if the component contains permutations that swap
adjacent columns of~$X$, see Hojny and Pfetsch~\cite{HojnyPfetsch2019}.
Such a swap of two columns is a permutation that decomposes into~2-cycles.
Therefore, our routine sieves all permutations~$P$ from a component that has
such a decomposition.
Then, we iteratively build a set of permutations~$Q \subseteq P$ that
define an orbitope or several independent orbitopes.
Initially, $Q = \emptyset$ and we check, one after another, whether
adding~$\perm \in P$ to~$Q$ allows to define independent orbitopes.
If this is possible, $Q$ is updated; otherwise, we continue with the next permutation
in~$P$ and discard~$\perm$.
To check whether we can add~$\perm$ to~$Q$, we maintain a list of the
orbitopes defined by the permutations in~$Q$ so far.
Then, $\perm$ is added to~$Q$ either if the variables affected by~$\perm$
are not contained in any of the already known orbitopes, or~$\perm$ adds a new
column to an already existing orbitope, or it merges two existing orbitopes.

If a component can not completely be handled by a single orbitope, there
might exist variables that are not contained in one of the detected
orbitopes, or several independent suborbitopes are found that are linked
via permutations not contained in~$Q$.
To partially add the missing link, and thus to handle more symmetries,
\scip selects one of the found orbitopes with variable matrix~$X \in
\{0,1\}^{s \times t}$ and computes one round of SST cuts with~$X_{11}$ as leader.
We refer to these inequalities as \emph{weak} inequalities, because they
weakly connect the found orbitopes without exploiting any further group
structure.
Besides weak inequalities, we can also add \emph{strong}
inequalities~$X_{11} \geq X_{12} \geq \dots \geq X_{1t}$ for every found
orbitope.
These cuts are called strong because they also exploit the group structure
allowing to arbitrarily reorder the columns of the orbitope.
Note that the strong inequalities are implicitly added by orbitope
constraints.
In some situations, however, \scip adds strong inequalities instead of
orbitopes as we explain next.

The detection of suborbitopes and application of strong and weak
inequalities can be controlled via the Boolean parameter
\param{propagating/symmetry/detectsubgroups}.
If the parameter value is \param{TRUE} (default), SCIP searches for
suborbitopes using the above mechanism.
A found orbitopes is called \emph{useful} if it has at least three columns.
The reason for this classification is that an orbitope with just two
columns can also be handled by orbisack constraints, which can more easily
be combined with other symmetry handling constraints.
Moreover, the Boolean parameters \param{propagating/symmetry/addstrongsbcs} and
\param{propagating/symmetry/addweaksbcs} enable/disable whether strong
inequalities are used if suborbitopes are not handled and whether weak
inequalities are used to handle more group structure, respectively.

\paragraph{\scip's Symmetry Handling Strategy}

As explained above, \scip allows to handle symmetries using different
strategies depending on the parameter \param{misc/usesymmetry}.
If a mixed strategy is used, \scip analyzes the structure of the symmetry
group's components and decides which strategy is used for which component.
In one case, however, these strategies can also be combined and applied
to the same component:
If both symresacks and SST cuts for binary variables are enabled, \scip
computes SST cuts first.
The leaders of SST cuts then play a special role, because they need to
attain the largest values in their orbits.
To make these cuts compatible with symresacks, one thus needs to adapt the
lexicographic order used by symresacks giving the leaders the highest rank.

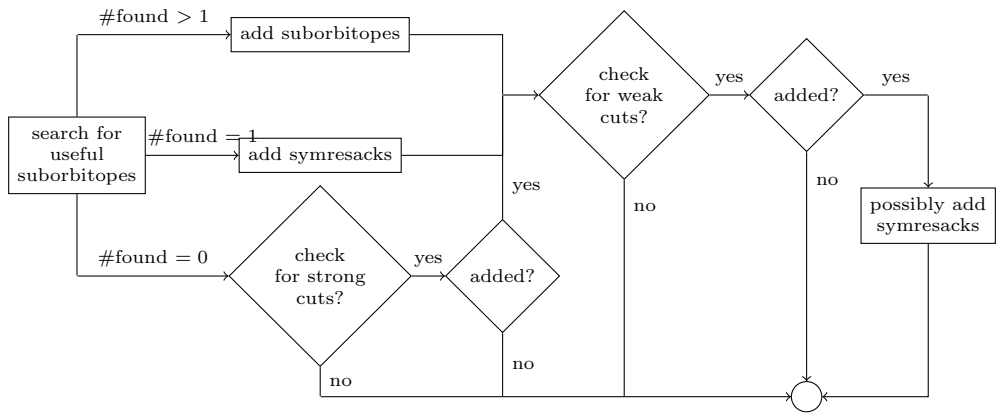
\begin{figure}[t]
  \centering
  \begin{tikzpicture}
    [scale=0.8,font=\scriptsize,
    terminal/.style={circle, minimum size=4mm, draw=black, align=center},
    decision/.style={diamond, minimum size=4mm, draw=black, align=center},
    nonterminal/.style={rectangle, minimum size=4mm, draw=black, align=center},
    pseudo/.style={inner sep=0mm, minimum size=0mm}]

    \node (orbitopes) at (2,0) [nonterminal] {search for\\ useful \\ suborbitopes};

    \node (checkstrong) at (6,-2) [decision] {check\\ for strong\\ cuts?};
    \node (symresacks) at (6,0) [nonterminal] {add symresacks};
    \node (suborbitopes) at (6,2) [nonterminal] {add suborbitopes};

    \node (foundstrong) at (9,-2) [decision] {added?};
    \node (checkweak) at (11,1) [decision] {check\\ for weak\\ cuts?};
    
    \node (foundweak) at (14,1) [decision] {added?};
    \node (subsymresacks) at (16,-1) [nonterminal] {possibly add\\ symresacks};

    \node (end) at (14,-4) [terminal] {};

    \node (p1) at (2,2) [pseudo] {};
    \node (p2) at (2,-2) [pseudo] {};
    \node (p3) at (6,-4) [pseudo] {};
    \node (p4) at (9,-4) [pseudo] {};
    \node (p5) at (9,1) [pseudo] {};
    \node (p6) at (11,-4) [pseudo] {};
    \node (p7) at (16,1) [pseudo] {};

    \draw[-] (orbitopes) -- (p1);
    \draw[->] (p1) to node [above] {\#{}found $> 1$} (suborbitopes);
    \draw[->] (orbitopes) to node [above] {\phantom{aa}\#{}found $= 1$} (symresacks);
    \draw[-] (orbitopes) -- (p2);
    \draw[->] (p2) to node [above] {\#{}found $= 0$} (checkstrong);

    \draw[->] (checkstrong) to node [above] {yes} (foundstrong);
    \draw[-] (checkstrong) to node [right] {no} (p3);
    \draw[->] (p3) -- (end);
    \draw[-] (foundstrong) to node [right] {no} (p4);

    \draw[-] (suborbitopes) -| (p5);
    \draw[-] (symresacks) -| (p5);
    \draw[->] (p5) -- (checkweak);
    \draw[-] (foundstrong) to node [right,near start] {yes} (p5);

    \draw[-] (checkweak) to node [right,very near start] {no} (p6);
    \draw[->] (checkweak) to node [above] {yes} (foundweak);

    \draw[->] (foundweak) to node [right, very near start] {no} (end);

    \draw[-] (foundweak) to node [above] {yes} (p7);
    \draw[->] (p7) -- (subsymresacks);
    \draw[->] (subsymresacks) |- (end);
  \end{tikzpicture}
  \caption{Rules to decide whether strong/weak inequalities are added.}
  \label{fig:strongweak}
\end{figure}

Similarly, if suborbitopes are detected, the orbitopes can be made
compatible with symresacks for the permutations not used by the orbitopes
by adapting the variable order in a specific way: the variables of the
first orbitope get the highest rank in the lexicographic order, afterwards
the variables of the succeeding orbitopes are listed, and finally the
variables not contained in any orbitope are added to the lexicographic
order.
The exact mechanism how suborbitopes are combined with weak and strong
inequalities as well as symresacks is explained in Figure~\ref{fig:strongweak}.

\scip's strategy to decide on which symmetry handling methods are used is
carried out in the following order, depending on the enabled strategies;
by default, \scip is allowed to use all implemented symmetry handling
methods (\param{misc/usesymmetry = 7}).
First, \scip checks whether a component can be fully handled by orbitopes
or whether suborbitopes can be detected.
If the component is handled by (sub)orbitopes, it gets blocked and no other
symmetry handling method can be applied to this component.
Second, \scip adds SST cuts to all applicable non-blocked components and
blocks these components.
If the selected leaders are binary, also symresacks can applied to this
component.
Third, if a component has not been blocked yet, either symresacks or
orbital fixing is used to handle symmetries, depending on whether orbital
fixing is active.

\subsubsection{Algorithmic Enhancements}

Besides new symmetry handling methods, \scip~8.0 also contains more
efficient implementations of previously available methods that we describe
in turn.

First, as mentioned above, orbisack constraints allow to apply stronger cutting
planes or reductions if they interact with set packing or partitioning
constraints in a certain way.
\scip automatically checks whether such an upgrade is possible.
The implementation of this upgrade has been revised and is more efficient
in \scip~8.0.

Second, the symresack constraint handler separates so-called minimal cover
inequalities for symresacks.
In \scip~7.0, we used a quadratic time separation routine for these
inequalities.
With the release of \scip~8.0, these inequalities can be separated in
linear time, which also improves on the almost linear running time
procedure by Hojny and Pfetsch~\cite{HojnyPfetsch2019}.
The linear time procedure makes use of the observation~\cite{HojnyPfetsch2019} that
minimal cover inequalities for symresacks can be separated by merging
connected components of an auxiliary graph.
Using a disjoint-set data structure, an almost linear running time could be
achieved.
In our new implementation, we exploit that the graph's connected components
are either paths or cycles.
Merging such connected components can be realized using more efficient data
structures based on a few arrays.

Finally, both the symresack and orbisack constraint handler provide
routines to propagate their constraints.
While the previous implementation could miss some variable fixings, the
implementation in \scip~8.0 allows to find all variable fixings that can be
derived from local variable bound information.

\subsubsection{Further Features}

Although symresack and orbitope constraints have been available in \scip since
version~5.0, these constraints could not be parsed in any file format.
With the release of \scip~8.0, these constraints can be parsed when reading
a \param{cip} file.
Thus, users can easily tell \scip about the symmetries that are present in
their problems and how to handle them.

For a permutation~$\perm$ of~$\{1,\dots,n\}$ and a vector~$x \in
\{0,1\}^n$, a symresack constraint enforces that~$x$ is lexicographically
not smaller than its permutation~$\perm(x)$.
This can be encoded in a \param{cip} file using the line
\[
  \texttt{symresack([varName1,\dots,varNameN],[$\perm(1)$,\dots,$\perm(n)$])}.
\]
Since orbisacks are symresacks for permutations that decompose into
2-cycles, this structure can directly be encoded using an~$\frac{n}{2}
\times 2$ matrix, where each row encodes the variables that can be
interchanged.
The \param{cip} encoding is then given by
\[
  \texttt{fullOrbisack(varName1-1,varName1-2.varName2-1,varName2-2. \dots)}.
\]
If users know that in each row of the orbitope matrix at most or
exactly one variable can attain value~1, they can provide this
information to \scip by replacing \texttt{fullOrbisack} by \texttt{packOrbisack} or
\texttt{partOrbisack}, respectively.

Finally, an orbitope constraint for a variable matrix~$X \in \{0,1\}^{m
  \times n}$ can be encoded similarly to an orbisack by the line
\[
  \texttt{fullOrbitope(varName1-1,\dots,varName1-N. \dots .varNameM-1,\dots,varNameM-N)}.
\]
If in each row of the orbitope at most or exactly one variable can attain
value~1, \texttt{fullOrbitope} can be replaced by \texttt{packOrbitope}
or~\texttt{partOrbitope}, respectively, to provide this information to \scip.


\subsection{Mixing Cuts}
\label{subsect:mixcuts}

%
Mixing cuts \cite{AtamturkEtAl2000,GunlukiPochet2001} can effectively reduce the computational time
to solve \MIP formulations of chance constrained programs (CCPs), especially for those in which
the uncertainty appears only in the right-hand side~\cite{LuedtkeEtAl2010,Kucukyavuz2012,AbdiFukasawa2016,ZhaoEtAl2017}.
In order to enhance the capability of employing \scip as a black box to solve 
such CCPs, \scipv includes a new separator called \texttt{mixing}, which 
leverages the \emph{variable bound relations} 
\cite{Achterberg2007a,MaherEtAl2017} to construct mixing cuts.
It is worthwhile remarking that though the development of this feature is motivated by CCPs,
the mixing separator can, however, be applied for other \MIPs
as long as the related variable bound relations can be detected by \scip.

Let us first review the variable bound relations in \scip; for more details, see Achterberg~\cite{Achterberg2007a}
and the SCIP Optimization Suite 4.0 release report~\cite{MaherEtAl2017}.
A variable bound relation in \scip is a linear constraint on two variables.
As such, it is of the form $y \star ax + b$ with $a$, $b \in
\R$ and $\star \in \{\leq, \geq\}$. 
During the presolving process, \scip derives these relations either from 
two-variable linear constraints or general constraints by probing 
\cite{Savelsbergh1994} and stores them in a data structure called \emph{variable bound graph}.
Such relations can be used to, for example, tighten the bounds of variables through propagation \cite{MaherEtAl2017}
or enhance the MIR cuts separation \cite{MarchandWolsey2001} in the subsequent main solution process.
The mixing cut separator uses a subclass of these relations, that is, those in
which $x$ is a binary variable and $y$ is a non-binary variable. Thereby, three 
families of cuts are constructed which is discussed in detail in the following.
For simplicity, we only consider the case that $y$ is a continuous variable but 
the result can also be applied to the case that $y$ is an integer variable.

\paragraph{$\geq$-Mixing Cuts}

Consider the variable lower bounds of  variable $y \in [\lb, \ub]$:
\begin{equation}\label{eq:variable_lower_bounds}
		y \geq a_i x_i + b_i,~x_i \in \{0,1\},~i \in \mathcal{N}.
\end{equation}
Without loss of generality, we impose the following assumption:
\begin{equation}\label{ass:Mixing}
\tag{A} 0 < a_i \leq \ub - \lb \text{ and } b_i = \lb \text{ for all }i \in \mathcal{N}.
\end{equation}
Indeed, assumption~\eqref{ass:Mixing} can be guaranteed by applying the 
following preprocessing steps in order:
\begin{enumerate}[label=(\roman*),leftmargin=5ex]
	\item If $a_i<0$, variable $x_i$ can be complemented by $1 - x_i$.
	If $a_i = 0$, $y \geq a_i x_i + b_i$ can be removed from
	\eqref{eq:variable_lower_bounds} and $\lb' : = \max\{\lb,  b_i\}$ is the new
	lower bound for $y$.
	\item If $a_i + b_i \leq \lb$, by $a_i > 0$ (from (i)), constraint $y \geq a_ix_i + b_i$ is implied by $y \geq
	\lb$ and hence can be removed from \eqref{eq:variable_lower_bounds}.
	\item If $b_{i} > \lb$, by $a_i > 0$ (from (i)), $\lb' := b_{i}$ is the new
	lower bound for $y$; if $b_i < \lb$, by $a_i + b_i > \lb$ (from (ii)),
	relation $y \geq a_ix_i + b_i$ can be changed into $y \geq (a_i + b_i - \lb)x_i +
	\lb$.
	\item If $a_i > \ub - \lb$, by $b_i = \lb$ (from (iii)), $x_i = 0$ must hold and
      constraint $y \geq a_i x_i + \lb$ can be removed from \eqref{eq:variable_lower_bounds}.
\end{enumerate}
By assumption~\eqref{ass:Mixing}, \eqref{eq:variable_lower_bounds} can be
presented in normalized form:
\begin{equation}
\label{eq:standard_variable_lower_bounds}
y \geq a_i x_i + \lb ,~x_i \in \{0,1\},~i \in \mathcal{N}.
\end{equation}
Let $\{  i_1, \ldots, i_s\} \subseteq \mathcal{N}$ with $s \in \N$
such that $a_{i_1} \leq \cdots \leq a_{i_s}$, and define $a_{i_0} := 0$. 
Then the $\geq$-mixing inequality \cite{AtamturkEtAl2000,GunlukiPochet2001} is given by
\begin{equation}\label{eq:mixingineq1}
	y - \lb \geq \sum_{\tau = 1}^s (a_{i_\tau} - a_{i_{\tau - 1}})\,x_{i_\tau}.
\end{equation}

\paragraph{$\leq$-Mixing Cuts}

Using a similar analysis as that in variable lower bounds, the variable upper bounds of
variable $y$ can be presented in normalized form:
\begin{equation}
\label{eq:standard_variable_upper_bounds}
y \leq \ub - a_j x_j,~x_j \in \{0,1\},~j \in \mathcal{M},
\end{equation}
where $0 < a_j \leq \ub - \lb$, $j\in \mathcal{M}$.
Let $\{j_1, \ldots, j_t\}\subseteq \mathcal{M}~(t \in \N)$ such that $a_{j_1}\leq \cdots \leq a_{j_t}$,
and define $a_{j_0} := 0$.
Then the $\leq$-mixing inequality \cite{AtamturkEtAl2000,GunlukiPochet2001} is given by
\begin{equation}\label{eq:mixingineq2}
y  \leq \ub - \sum_{\tau = 1}^t (a_{j_\tau} - a_{j_{\tau - 1}})\,x_{j_\tau}.
\end{equation}

\paragraph{Conflict Cuts}

Besides the $\geq$- and $\leq$-mixing cuts, the mixing separator also constructs conflict cuts,
which are derived by jointly considering
\eqref{eq:standard_variable_lower_bounds} and \eqref{eq:standard_variable_upper_bounds}.
To be more specific, let $i' \in \mathcal{N}$ and $j' \in \mathcal{M}$ such that $a_{i'} + \lb > \ub - a_{j'}$.
By $y\geq a_{i'} x_{i'} + \lb$ and $y \leq u - a_{j'} x_{j'}$, variables $x_{i'}$ and $x_{j'}$ cannot
simultaneously take values at one, and hence the conflict inequality
\begin{equation}\label{eq:conflictineq}
	x_{i'} + x_{j'} \leq 1
\end{equation}
can be derived.

\paragraph{Separation}

Given a fractional point $(x^*,y^*)$, the separation problem of \eqref{eq:mixingineq1},
\eqref{eq:mixingineq2} or \eqref{eq:conflictineq} asks to find an inequality violated by $(x^*,y^*)$ or
prove that no such one exists.
To separate the $\geq$-mixing inequalities \eqref{eq:mixingineq1}, G\"{u}nl\"{u}k and Pochet \cite{GunlukiPochet2001}
provided the following algorithm, which selects the subset $\mathcal{S} = \{i_1, \ldots, i_\tau\} \subseteq \mathcal{N}$ such that
$\sum_{\tau = 1}^s (a_{i_\tau} - a_{i_{\tau - 1}})\,x^*_{i_\tau}$ is maximized.
\begin{enumerate}
	\item Reorder variables $x_{i}$, $i \in \mathcal{N}$, such that
      $x^*_{1} \geq x^*_{2} \geq \cdots \geq x^*_{|\mathcal{N}|}$.
	\item Add $1$ to set $\mathcal{S}$.
	\item For each $i \in \mathcal{N} \backslash \{1\}$, set $\mathcal{S} := \mathcal{S} \cup \{i\}$ if $a_i > a_{k}$,
      where $k$ is last index added into $\mathcal{S}$.
	\item If the $\geq$-mixing inequality corresponding to $\mathcal{S}$ is violated by $(x^*,y^*)$, output it.
\end{enumerate}
Obviously, the above algorithm can be implemented to run in $\mathcal{O}(|\mathcal{N}|\log(|\mathcal{N}|))$.
Similarly, the $\leq$-mixing inequalities \eqref{eq:mixingineq2} can also be separated in
$\mathcal{O}(|\mathcal{M}|\log(|\mathcal{M}|))$.
Finally, for the conflict inequalities  \eqref{eq:conflictineq},
since number of them is bounded by $\mathcal{O}(|\mathcal{M}||\mathcal{N}|)$,
by enumeration, they can be separated in $\mathcal{O}(|\mathcal{M}||\mathcal{N}|)$.

The performance impact of the mixing separator is neutral on the internal \MIP benchmark testset.
However, when applied to the chance constrained lot sizing instances used by Zhao et al.~\cite{ZhaoEtAl2017} ($90$ in total),
a speedup of $20\%$ can be observed and $15$ more instances can be solved.

\subsection{Primal Decomposition Heuristics}
\label{subsect:heuristics}


\scipv comes with an improvement of the heuristic Penalty Alternating Direction Method (\method{PADM})
and introduces the new heuristic Dynamic Partition Search (\method{DPS}).
Both heuristics explicitly require a decomposition provided by the user
and therefore belong to the class of so-called \emph{decomposition heuristics}.
A decomposition consisting of $k\geq0$ blocks is a partition
\begin{equation*}
\mathcal{D}\coloneqq (D^{\text{row}},D^{\text{col}})
\text{ with }
D^{\text{row}} \coloneqq (D^{\text{row}}_{1},\dots,D^{\text{row}}_{k},L^{\text{row}}), \ D^{\text{col}} \coloneqq (D^{\text{col}}_{1},\dots,D^{\text{col}}_{k},L^{\text{col}})
\end{equation*}
of the rows/columns of the constraint matrix $\linmatrix$ into $k+1$ pieces each.
The distinguished rows $L^{\text{row}}$ and columns $ L^{\text{col}}$ are called \emph{linking rows} and \emph{linking columns}, respectively.
If $\linmatrix$ is permuted according to decomposition $\mathcal{D}$, a \emph{(bordered) block diagonal form}~\cite{Borndoerfer1998} is obtained.
A detailed description of decompositions and their handling in \scip can be found in the release report for version 7.0~\cite{SCIP7}.

\subsubsection{Improvement of Penalty Alternating Direction Method}
\label{subsubsec:padm}

Since version 7.0, \scip includes the decomposition heuristic Penalty Alternating Direction Method (\method{PADM}).
For the current version \method{PADM} has been extended by the option
to improve a found solution by reintroducing the original objective function.

This heuristic splits a \MINLP as listed in~\eqref{eq:generalminlp} into several subproblems
according to a given decomposition $\mathcal{D}$ with linking variables only,
whereby the linking variables get copied and the differences are penalized.
Then the subproblems are solved by an alternating procedure.
A detailed description of penalty alternating direction methods
and their practical application
can be found in Geißler et al.~\cite{Geissleretal2017}
and Schewe et al.~\cite{ScheweSchmidtWeninger2019}.

To converge faster to a feasible solution, the original objective function of each subproblem
has been completely replaced by a penalty term.
Since this can lead to arbitrarily bad solutions,
the heuristic was extended in the following way:
Initially, the original version of \method{PADM} runs
and tries to find a feasible solution.
If a feasible solution was found the linking variables are fixed to the values of this solution
and each independent subproblem is solved again but now with the original objective function.
In order to accelerate the solving of the reoptimization step,
the already found solution is used as a warm start solution and very small solving limits are imposed.
The additional reoptimization step must not take more time than was already used
by the heuristic in the first step and the node limit is set to one.
By setting the parameter \param{heuristics/padm/reoptimize} the feature of using
a second reoptimization step in \method{PADM} can be turned on/off (default: on).

The new feature was tested on the \miplib2017~\cite{miplib2017} benchmark instances, for which decompositions are provided on the web page.
If \method{PADM} could get called, preliminary results show that \method{PADM} finds a feasible solution in 15 of 31 cases.
The new reoptimization step successfully improves the solution of \method{PADM} in 33\% of these cases by 42\% on average.

\subsubsection{Dynamic Partition Search}

With \scipv the new decomposition heuristic Dynamic Partition Search (\method{DPS}) was added.
It is a primal construction heuristic which requires a decomposition with linking constraints only.

The \method{DPS} heuristic splits a \MILP as listed in~\eqref{eq:generalmip} into several subproblems according to a decomposition $\mathcal{D}$.
Thereby the linking constraints and their right/left-hand sides are also split
by introducing new parameters $p_q \in \R^{|L^{\text{row}}|}$ for each block $q \in \{1,\dots,k\}$ and
requiring that
\begin{equation}\label{eq:DecompPartition}
  \sum_{q=1}^{k} p_q = b_{[L^{\text{row}}]}
\end{equation}
holds.
To obtain information about the infeasibility of one subproblem and to speed up the solving process,
slack variables $z_q \in \R^{|L^{\text{row}}|}_+$ are added and the objective function is replaced by a weighted sum of these slack variables.
In detail, for penalty parameter $\lambda \in \R^{|L^{\text{row}}|}_{> 0}$
each subproblem $q$ has the form
\begin{equation}
\label{eq:dpsblock}
\begin{aligned}
\min\; & \lambda^\T z_q, \\
\mathrm{s.t.}\; & A_{[D^{\text{row}}_{q},D^{\text{col}}_{q}]}x_{[D^{\text{col}}_{q}]} \geq \rhs_{[D^{\text{row}}_{q}]}, \\
& \lb_{i} \leq x_{i} \leq \ub_{i} && \fa i \in \varindex \cap D^{\text{col}}_{q}, \\
& x_{i} \in \Z && \fa i \in \intvarindex \cap D^{\text{col}}_{q},\\
& A_{[L^{\text{row}},D^{\text{col}}_{q}]}x_{[D^{\text{col}}_{q}]} + z_q \geq p_{q},\\
& z_q \in \R^{|L^{\text{row}}|}_+.
\end{aligned}
\end{equation}

From~\eqref{eq:dpsblock}, it is immediately apparent that the correct choice of $p_{q}$ plays the central role.
Because if $p_{q}$ is chosen for each subproblem $q$ such that the slack variables $z_q$ take the value zero,
one immediately obtains a feasible solution.
For this reason, we refer to $(p_q)_{q \in \{1,\dots,k\}}$ as a \emph{partition} of $b_{[L^{\text{row}}]}$.
The goal of \method{DPS} is to find a feasible partition as fast as possible.

To get started, an initial partition $(p_q)_{q \in \{1,\dots,k\}}$ is chosen,
which fulfills~\eqref{eq:DecompPartition}.
Then it is checked whether this partition will lead to a feasible solution
by solving $k$ independent subproblems \eqref{eq:dpsblock}
with fixed $p_q$.
If all subproblems have an optimal objective value of zero,
a feasible solution was found and is given by the concatenation of the $k$ subsolutions.
Conversely, a lower bound on the objective function value of one subproblem
of greater than zero immediately provides evidence that the current partition does not lead to a feasible solution.

If the current partition does not correspond to a feasible solution, then partition $(p_q)_{q \in \{1,\dots,k\}}$ and penalty parameter $\lambda$
have to be updated:
For each single linking constraint $j \in L^{\text{row}}$ the value vector $z_j$
is subtracted from the current partition $p_j$
and the same amount is added to all blocks with $z_{jq} = 0$, so that~\eqref{eq:DecompPartition} still holds.
If at least one slack variable is positive, the corresponding penalty parameter is increased.
Then, the subproblems are solved again and the steps are repeated
until a feasible solution is found or until a maximum number of iterations
(controlled by parameter \param{heuristics/dps/maxiterations})
is reached.

To push the slack variables to zero and to speed up the algorithm,
the original objective function has been completely replaced by a penalty term.
Analogously to \method{PADM} (see Section~\ref{subsubsec:padm}) it is possible to improve the found solution
by reoptimizing with the original objective function.
In \method{DPS} the partition instead of the linking variables is fixed.
By setting parameter \param{heuristics/dps/reoptimize} this feature can be turned on/off (default: off).

The new decomposition heuristic was tested on the \miplib2017~\cite{miplib2017} benchmark instances, for which decompositions are provided on the web page.
If \method{DPS} could get called, preliminary results show that \method{DPS} finds a feasible solution in 17 of 80 cases.
A general performance improvement can not be shown.
The main reason for these slightly disappointing results is probably that \method{DPS} requires a well-decomposable problem structure.
The evaluated instances are general \MILPs which do not necessarily have such a structure.
However, on two instances (\texttt{proteindesign121hz512p9} and \texttt{30n20b8})  \method{DPS} is successful and reduces the time until the first found primal solution highly,
since no other heuristic is able to construct a feasible solution at or before the root node.
It is noticeable that in both instances the linking constraints contain only bounded integer variables.
The heuristic probably benefits from this, since the number of usable partitions is thus countable and finite.

\subsection{Benders' Decomposition}
\label{subsect:benders}

The work on the Benders' decomposition framework has moved into a research phase.
As such, only minor updates and bug fixes have been completed for the framework since the release of \scip 7.0.
The most important update for the Benders' decomposition framework is the option to apply the mixed integer rounding (MIR) procedure, as described by Achterberg~\cite{Achterberg2007a}, when generating optimality cuts.
The aim of applying the MIR procedure to the generated optimality cut is to potentially compute a stronger inequality.

Strengthening the classical Benders' optimality cut using the MIR procedure involves the following steps:
\begin{itemize}
  \item Generate a classical optimality cut from the solution of the Benders' decomposition subproblem.
  \item Attempt to compute a flow cover cut for the generated optimality cut. This is achieved by calling \texttt{SCIPcalcFlowCover}. If this process is successful, replace the optimality cut with the computed flow cover cut.
  \item Attempt to perform the MIR procedure on the optimality cut (this could have been updated in the previous step). The MIR procedure is performed by calling \texttt{SCIPcalcMIR}. If the MIR procedure is successful, the optimality cut is replaced with the resulting inequality.
  \item Finally, \texttt{SCIPcutsTightenCoefficients} is executed in an attempt to tighten the coefficients of the optimality cut.
\end{itemize}

The MIR procedure is active by default.
A new parameter 
\begin{center}
  \texttt{benders/<bendersname>/benderscut/optimality/mir}, 
\end{center}
where \texttt{<bendersname>} is the name of the Benders' decomposition plugin, has been added to enable/disable the MIR procedure for strengthening the Benders' optimality cuts.

\subsection{Cut Selectors}
\label{subsect:cutselectors}


The new cut selector plugin is introduced in SCIP 8.0. Users now have the ability to create their own cut selection rules and include them into SCIP. For a current summary on the state of cut selection in the literature, see Dey and Molinaro~\cite{DeyMolinaro2018}, and for an overview of cutting plane measures and the improvements provided by intelligent selection, see Wesselmann and Suhl~\cite{Wesselmann2012}. The existing rule used since SCIP 6.0~\cite{SCIP6} has been moved to \texttt{cutselection/hybrid}. The ability to include cut selectors has also been implemented through PySCIPOpt.

\subsection{Technical Improvements}
\label{subsect:further}


\paragraph{Thread Safety}

In previous versions, \scip contained the argument \texttt{PARASCIP} for
the Make and CMake build system to make it thread-safe. This has been replaced by
\texttt{THREADSAFE}, which is now true by default (\texttt{PARASCIP} still
exists for backward compatibility).

Most parts of \scip are in fact always thread-safe, but interfaces to
external programs are sometimes not. For instance, for the LP-solver
\gurobi, the thread-safe mode opens a new LP-environment for each thread.
Other interfaces to external software may use parallelization that has to
be controlled in order not to mix data from different threads, e.g., \cppad
and \filtersqp. The change to thread-safe mode should not significantly
affect performance.

\paragraph{Revision of External Memory Estimation}

\scip usually uses its own internal memory functions. This allows to keep
track of the used memory. If it approaches the memory limit, \scip can
switch to a memory saving mode, which, for instance, uses
depth-first-search. However, memory used by external software, in
particular, NLP and LP-solvers cannot easily be determined in a portable
way. Therefore, the estimation of used memory in \scip has been improved
for version 8 with data-fitting as follows. The memory
consumption by LP-solvers was measured using a stand-alone version on a
testset of LP-relaxations. Then a linear regression with the number of
constraints, variables, and nonzeros as features was computed. This current
estimation uses the weights $\num[round-precision=1]{8.5e-04}$,
$\num[round-precision=1]{7.6e-04}$, $\num[round-precision=1]{3.5e-05}$,
respectively, and works quite well (for \soplex we got $R^2 = 0.99$). If
NLPs are solved, the estimation is doubled.

\paragraph{Option to Forbid Variable Aggregation}

Similar to multi-aggregation, one can now forbid aggregation of a variable
by calling the function \texttt{SCIPdoNotAggrVar()}. This is sometimes
useful, for example, if certain constraint handlers cannot handle aggregated
variables. Note, however, that this can slow down the solving process since
the relaxations tend to be larger.

\paragraph{Debugging of Variable Uses}

SCIP counts the number of uses of a variable and frees a variable when its uses count reaches zero.
It is therefore important to capture a variable to prevent it from being freed too early and to release a variable when it is no longer used.
To assist on finding a missing or excessive capture or release of a variable, code has been added to \texttt{var.c} to print the function call stack when a variable of a specified name, optionally in a SCIP problem of a specified name, is captured or released.
The code requires GCC Gnulib (\texttt{execinfo.h} in particular) and will not work on every platform.

To activate this feature, define \texttt{DEBUGUSES\_VARNAME} and \texttt{DEBUGUSES\_PROBNAME} in \texttt{var.c}.
If the tool \texttt{addr2line} is available on the system, the printed call stacks provide more information, but its use causes a significant slowdown.
Defining \texttt{DEBUGUSES\_{\allowbreak}NOADDR2LINE} disables the call of this tool.

\paragraph{Improving Numerical Properties of Linear Inequalities}

When a constraint handler or separator computes a cutting plane, often its numerical properties need to be checked and possibly improved before it is added to a relaxation.
Further, changing coefficients or sides of a \texttt{SCIP\_ROW} may round numbers that are very close to integral values, which may invalidate a previously valid cut.
To assist on carefully improving the numerical properties of an inequality, the \texttt{SCIP\_ROWPREP} datastructure has been made available, see \texttt{pub\_misc\_rowprep.h}.
Routines are available to relax or scale linear inequalities to improve the range of coefficients and to avoid almost-integral numbers, see the paragraph ``Cut cleanup'' in Section~\ref{sect:nlsepa} and Step~\ref{enfosepa} in Section~\ref{sect:nlenfo} for more details.
Note, that ranged linear constraints (both left-hand-side and right-hand-side being finite) cannot be handled.

\paragraph{Reader for AMPL \texttt{.nl} Files}

The reader for \texttt{.nl} files has been rewritten and is now included with SCIP's default plugins.
See Section~\ref{sect:ampl} for more details.

\section{SCIP's New MINLP Framework}\label{sect:minlp}

\subsection{New Expressions Framework}
\label{sect:expr}

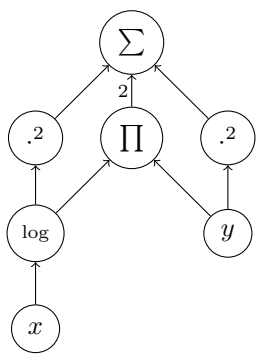
\begin{figure}[ht]
  \centering
  \begin{tikzpicture}[scale=0.6, sibling distance=6em, level distance=6em]
    \tikzstyle{expr} += [shape=circle, draw, align=center, minimum size = 18pt] 

    \node[expr] (plus) {$\sum$}
      child[<-] {
        node[expr] (sqr1) {$\cdot^2$}
        child[<-] {
          node[expr] (log) {\scriptsize $\log$}
          child[<-] {
            node[expr] (x) {$x$}
          }
        }
      }
      child[<-] {
        node[expr] (prod) {$\prod$}
        edge from parent node[left,inner sep=.1em] {\scriptsize 2}
      }
      child[<-] {
        node[expr] (sqr2) {$\cdot^2$}
        child[<-] {
          node[expr] (y) {$y$}
        }
      };

    \draw[<-] (prod) to[out=-135,in=45] (log);
    \draw[<-] (prod) to[out=-45,in=135] (y);

    \node[yshift=+0.5cm] (t3) at (plus) {};

    \node[yshift=-0.5cm]               (t1) at (x) {};
    \node[yshift=-0.5cm]               (t2) at (y) {};
  \end{tikzpicture}
  \caption{Expression graph for algebraic expression $\log(x)^2 + 2\log(x)y+y^2$.}
  \label{fig:exprgraph}
\end{figure}

Algebraic expressions are well-formed combinations of constants, variables, and various algebraic operations such as addition, multiplication, exponentiation, that are used to describe mathematical functions.
They are often represented by a directed acyclic graph with nodes representing variables, constants, and operations and arcs indicating the flow of computation, see Figure~\ref{fig:exprgraph} for an example.

With SCIP 7.0 and before, the following node types were supported in algebraic expressions:
\begin{itemize}
\item constant, parameter,
\item variable (specified by integer index),
\item addition, subtraction, multiplication, division  (with two arguments),
\item square, square-root, power (rational exponent), power (integer exponent), signed power ($x\mapsto \mathrm{sign}(x)\abs{x}^p$),
\item exponentiation, natural logarithm,
\item minimum, maximum, absolute value,
\item sum, product, affine-linear, quadratic, signomial  (with arbitrarily many arguments),
\item user-defined.
\end{itemize}
The operand type ``user-defined'', which was introduced with SCIP 3.2~\cite{SCIP32}, brought some of the extensibility typical for SCIP plugins.
However, only the most essential callbacks (evaluation, differentiation, linear under/overestimation) were defined for user-defined expressions.
Thus, other routines that worked on expressions, such as simplification, had built-in treatment for operands integrated in SCIP, but defaulted to some simple conservative behavior when a user-defined operand had to be dealt with.

Another problem with this design of expressions was the ambiguity and additional complexity due to the presence of high-level operators such as affine-linear, quadratic, and others.
For example, code that did some operation on a sum had to implement the same routine for any operand that represents some form of summation (plus, sum, affine-linear, quadratic, signomial), each time dealing with a slightly different data structure.

With SCIP 8, the expression system has been completely rewritten.
Proper SCIP plugins, referred to as \textit{expression handlers}, are now used to define all semantics of an operand.
These expression handlers support more callbacks than what was available for the user-defined operator before.
Furthermore, much ambiguity and complexity is avoided by adding expression handlers for basic operations only.
High-level structures such as quadratic functions can still be recognized, but are no longer made explicit by a change in the expression type.
An expression (\texttt{SCIP\_EXPR}) comprises various data, in particular the arguments of the expression, further on denoted as \textit{children}.
It can also hold the data and additional callbacks of an \textit{expression owner}, if any.
A prominent example of an expression owner is the constraint handler for nonlinear constraints, see Section~\ref{sect:consnl}, which stores data associated with the enforcement of nonlinear constraints in expressions that are used to specify nonlinear constraints.
Further, due to their many use cases, a representation of the expression as a quadratic function can be stored, see Section~\ref{sect:nlhdlrquaddetect} for details.
Various methods are available in the SCIP core to manage expressions (create, modify, copy, free, parse, print), to evaluate and compute derivative information at a point, to evaluate over intervals, to simplify, to identify common subexpressions, to check curvature and integrality, and to iterate over it.
Many of these methods access callbacks that can be implemented by expression handlers.
Some additional callbacks are used by the constraint handler for nonlinear constraints (Section~\ref{sect:consnl}).
The expression handler callbacks are:
\begin{itemize}
\item \texttt{COPYHDLR}: include expression handler in another SCIP instance;
\item \texttt{FREEHDLR}: free expression handler data;
\item \texttt{COPYDATA}: copy expression data, for example, the coefficients of a linear sum;
\item \texttt{FREEDATA}: free expression data;
\item \texttt{PRINT}: print expression;
\item \texttt{PARSE}: parse expression from string;
\item \texttt{CURVATURE}: detect convexity or concavity;
\item \texttt{MONOTONICITY}: detect monotonicity;
\item \texttt{INTEGRALITY}: detect integrality (is value of operation integral if arguments have integral value?);
\item \texttt{HASH}: hash expression using hash values of arguments;
\item \texttt{COMPARE}: compare two expressions of same type;
\item \texttt{EVAL}: evaluate expression (implementation of this callback is mandatory);
\item \texttt{BWDIFF}: evaluate partial derivative of expression with respect to specified argument (backward derivative evaluation);
\item \texttt{FWDIFF}: evaluate directional derivative of expression (forward derivative evaluation);
\item \texttt{BWFWDIFF}: evaluate directional derivative of partial derivative with respect to specified argument (backward over forward derivative);
\item \texttt{INTEVAL}: evaluate expression over interval;
\item \texttt{ESTIMATE}: compute linear under- or overestimator of expression with respect to given bounds on arguments and a reference point;
\item \texttt{INITESTIMATES}: compute one or several linear under- or overestimators of expression with respect to given bounds on arguments;
\item \texttt{SIMPLIFY}: simplify expression by applying algebraic transformations;
\item \texttt{REVERSEPROP}: compute bounds on arguments of expression from given bounds on expression.
\end{itemize}
The SCIP documentation provides more details on these callbacks.

Finally, for the following operators, expression handlers are included in SCIP 8.0:
\begin{itemize}
\item \texttt{val}: scalar constant;
\item \texttt{var}: a SCIP variable (\texttt{SCIP\_VAR});
\item \texttt{varidx}: a variable represented by an index; this handler is only used for interfaces to NLP solvers (NLPI);
\item \texttt{sum}: an affine-linear function, $y\mapsto a_0 + \sum_{j=1}^k a_jy_j$ for $y\in\R^k$ with constant coefficients $a\in\R^{k+1}$;
\item \texttt{prod}: a product, $y\mapsto c\prod_{j=1}^ky_j$ for $y\in\R^k$ with constant factor $c\in\R$;
\item \texttt{pow}: a power with a constant exponent, $y\mapsto y^p$ for $y\in \R$ and exponent $p\in\R$ (if $p\not\in\Z$, then $y\geq 0$ is required);
\item \texttt{signpower}: a signed power, $y\mapsto \mathrm{sign}(y)\abs{y}^p$ for $y\in\R$ and constant exponent $p\in\R$, $p>1$;
\item \texttt{exp}: exponentiation, $y\mapsto \exp(y)$ for $y\in\R$;
\item \texttt{log}: natural logarithm, $y\mapsto \log(y)$ for $y\in\R_{>0}$;
\item \texttt{entropy}: entropy, $y\mapsto\begin{cases}-y\log(y), & \textrm{if }y > 0,\\0, & \textrm{if }y=0,\end{cases}$ for $y\in\R_{\geq 0}$;
\item \texttt{sin}: sine, $y\mapsto\sin(y)$ for $y\in\R$;
\item \texttt{cos}: cosine, $y\mapsto\cos(y)$ for $y\in\R$;
\item \texttt{abs}: absolute value, $y\mapsto \abs{y}$ for $y\in\R$.
\end{itemize}
When comparing with the list for SCIP~7.0 above, one observes that support for parameters (these behaved like constants but could not be simplified away and were modifiable) and operators ``min'' and ``max'' has been removed.
Further, support for sine, cosine, and the entropy function has been added.


\subsection{New Handler for Nonlinear Constraints}
\label{sect:consnl}


For SCIP 8, the constraint handler for general nonlinear constraints (\texttt{cons\_nonlinear}) has been rewritten and the specialized constraint handlers for quadratic, second-order cone, absolute power, and bivariate constraints have been removed.
Some of the unique functionalities of the removed constraint handlers has been reimplemented in other plugin types. 

\subsubsection{Motivation}
\label{sect:nlmotivation}

An initial motivation for the rewrite of \texttt{cons\_nonlinear} has been a numerical issue which is caused by the explicit reformulation of constraints in SCIP 7.0 and earlier versions.
For an example, consider the problem
\begin{equation}
 \label{eq:consnl_exorig}
 \begin{aligned}
  \min\; & z, \\
  \text{s.t.}\; & \exp(\ln(1000)+1+x\,y) \leq z, \\
   & x^2 + y^2 \leq 2,
 \end{aligned}
\end{equation}
with optimal solution $x=-1$, $y=1$, $z=1000$.
Previously, solving this problem with SCIP could end with the following solution report:
{\small
\begin{verbatim}
SCIP Status        : problem is solved [optimal solution found]
Solving Time (sec) : 0.08
Solving Nodes      : 5
Primal Bound       : +9.99999656552062e+02 (3 solutions)
Dual Bound         : +9.99999656552062e+02
Gap                : 0.00 %
  [nonlinear] <e1>: exp((7.9077552789821368 + (<x>*<y>)))-<z>[C] <= 0;
violation: right-hand side is violated by 0.000673453314561812
best solution is not feasible in original problem

x                                   -1.00057454873626 	(obj:0)
y                                   0.999425451364613 	(obj:0)
z                                    999.999656552061 	(obj:1)
\end{verbatim}}
The reason that SCIP initially determined this solution to be feasible is that, in presolve, the problem gets rewritten as
\begin{equation}
 \label{eq:consnl_extrans}
 \begin{aligned}
  \min\; & z, \\
  \text{s.t.}\; & \exp(w) \leq z, \\
   & \ln(1000)+1+x\,y = w, \\
   & x^2 + y^2 \leq 2.
 \end{aligned}
\end{equation}
The constraints in this transformed problem are violated by $0.4659\cdot 10^{-6}$, $0.6731\cdot 10^{-6}$, and $0.6602\cdot 10^{-6}$, thus are feasible with respect to \texttt{numerics/feastol}$=10^{-6}$, and therefore the solution is accepted by SCIP.
On the MINLPLib library, the problem that a final solution is feasible for the presolved problem but violates nonlinear constraints in the original problem occurred for 7\% of all instances.

Problem~\eqref{eq:consnl_exorig} gets rewritten as~\eqref{eq:consnl_extrans} for the purpose of constructing a linear relaxation.
In this process, nonlinear functions are approximated by linear under- and overestimators.
As the formulas that were used to compute these estimators are only available for ``simple'' functions (for example, convex functions, concave functions, bilinear terms), new variables and constraints were introduced to split more complex expressions into adequate form~\cite{SmithPantelides1999,VigerskeGleixner2017}.

A trivial attempt to solve the issue of solutions not being feasible in the original problem would have been to add a feasibility check before accepting a solution.
However, if a solution is not feasible, actions to resolve the violation of original constraint need to be taken, such as a separating hyperplane, a domain reduction, or a branching operation needs to be performed. Since the connection from the original to the presolved problem was not preserved, it would not have been clear which operations on the presolved problem would help best to remedy the violation in the original problem.

Thus, the new constraint handler aims to preserve the original constraints by applying only transformations (simplifications) that, in most situations, do not relax the feasible space when taking tolerances into account.
The reformulations that were necessary for the construction of a linear relaxation are not applied explicitly anymore, but handled implicitly by annotating the expressions that define the nonlinear constraints (here, the mysterious ``data of an expression owner'', see Section~\ref{sect:expr}, comes into play).
Another advantage of this approach is a clear distinction between the variables that were present in the original problem and the variables added for the reformulation.
With this information, branching is avoided on variables of the latter type.
Finally, it is now possible to exploit overlapping structures in an expression simultaneously.

\subsubsection{Extended Formulations}
\label{sect:nlextform}


To explain the functionality of the new \texttt{cons\_nonlinear}, consider MINLPs of the form
\begin{equation}
  \label{eq:minlp}
  \tag{MINLP}
  \begin{aligned}
    \min\; & \linobj^\T x, \\
    \mathrm{s.t.}\; & \low{g} \leq g(x) \leq \upp{g}, \\
    & \low{b} \leq Ax \leq \upp{b}, \\
    & \low{x} \leq x \leq \upp{x}, \\
    & x_\intvarindex \in \Z^{\intvarindex},
  \end{aligned}
\end{equation}
with $\linobj\in\R^n$, $g:\R^n\to\overline{\R}^m$,
$\low{g}$, $\upp{g}\in\overline{\R}^m$,
$A\in\R^{\tilde m\times n}$,
$\low{b}$, $\upp{b}\in\smash{\overline{\R}^{\tilde m}}$,
$\low{x}$, $\upp{x}\in\overline{\R}^n$,
$\intvarindex\subseteq \{1,\ldots,n\}$,
$\overline{\R} \coloneqq \R \cup \{\pm\infty\}$.
Further, assume that $g_i(\cdot)$ is nonlinear and specified by an expression (see Section\ \ref{sect:expr}), $i=1,\ldots,m$,
$\low{g}\leq\upp{g}$,
$\low{g}_i \in\R$ or $\upp{g}_i \in\R$ for all $i=1,\ldots,m$,
$\low{b}\leq\upp{b}$,
$\low{b}_i \in\R$ or $\upp{b}_i \in\R$ for all $i=1,\ldots,\tilde m$, and
$\low{x}\leq\upp{x}$.
All nonlinear constraints $\low{g}\leq g(x)\leq \upp{g}$ are handled by \texttt{cons\_nonlinear}, while the linear constraints are handled by \texttt{cons\_linear} or its specializations.
(Of course, in general, any kind of constraint that SCIP supports is allowed, but for this section only linear and nonlinear constraints are considered.)
In comparison to SCIP 7.0, the specialized nonlinear constraint handlers and the distinction into a linear and a nonlinear part of a nonlinear constraint have been removed. 
As a consequence, all algorithms for nonlinear constraints (checking feasibility, domain propagation, separation, etc) work on expressions now.

SCIP solves problems like \eqref{eq:minlp} to global optimality via a spatial branch-and-bound algorithm that mixes branch-and-infer and branch-and-cut~\cite{BeKiLeLiLuMa12}.
Important parts of the solution algorithm are presolving, domain propagation (that is, tightening of variable bounds), linear relaxation, and branching.
For the domain propagation and linear relaxation aspects, two extended formulations of \eqref{eq:minlp} that are obtained by introducing \emph{slack variables} and replacing sub-trees of the expressions that define nonlinear constraints by \emph{auxiliary variables} are considered.

For domain propagation, the following extended formulation is considered:
\begin{equation}
  \label{eq:minlp_extdp}
  \tag{$\text{MINLP}_\text{ext}^\text{dp}$}
  \begin{aligned}
    \min\; & \linobj^\T x, \\
    \mathrm{s.t.}\; & h^\text{dp}_i(x,w^\text{dp}_{i+1},\ldots,w^\text{dp}_{m^\text{dp}}) = w_i^\text{dp}, & i=1,\ldots,m^\text{dp}, \\
    & \low{b} \leq Ax \leq \upp{b}, \\
    & \low{x} \leq x \leq \upp{x}, \\
    & \low{w}^\text{dp} \leq w^\text{dp} \leq \upp{w}^\text{dp}, \\
    & x_\intvarindex \in \Z^{\intvarindex}.
  \end{aligned}
\end{equation}
%
Initially, slack variables $w_1^\text{dp},\ldots,w_m^\text{dp}$ are introduced and $h_i^\text{dp}\defi g_i$ for $i=1,\ldots,m$.
Next, for each function $h_i^\text{dp}(x)$, subexpressions $f(x)$ may be replaced by new auxiliary variables $w_{i'}$, $i'>m$, and new constraints $h_{i'}^\text{dp}(x) = w_{i'}^\text{dp}$ with $h_{i'}^\text{dp} \defi f$ are added.
For the latter, subexpressions may be replaced again.
Since auxiliary variables that replace subexpression of $h_i^\text{dp}(x)$ always receive an index larger than $\max(m,i)$, the result is referred to by $h_i^\text{dp}(x,w_{i+1}^\text{dp},\ldots,w_{m^\text{dp}}^\text{dp})$ for any $i=1,\ldots, m^\text{dp}$.
That is, to simplify notation, $w_{i+1}^\text{dp}$ is used instead of $w_{\max(i,m)+1}^\text{dp}$.
If a subexpressions that is replaced by an auxiliary variable appears in several places, then only one auxiliary variable and one constraint is added to the extended formulation.
Reindexing may be necessary to have $h_i^\text{dp}$ depend on $x$ and $w_{i+1}^\text{dp},\ldots$ only.

The details of how subexpressions are chosen to be replaced by auxiliary variables will be discussed in Section~\ref{sect:detect}.
For the moment it is sufficient to assume that algorithms are available to compute interval enclosures of
\begin{align}
  \{ h_i^\text{dp}(x,w_{i+1}^\text{dp},\ldots,w_{m^\text{dp}}^\text{dp}) :
      \low{x} \leq x \leq \upp{x},\; \low{w}^\text{dp} \leq w^\text{dp} \leq \upp{w}^\text{dp} \},  \label{eq:inteval} \\
  \{ x_j : h_i^\text{dp}(x,w_{i+1}^\text{dp},\ldots,w_{m^\text{dp}}^\text{dp}) = w_i^\text{dp} :
      \low{x} \leq x \leq \upp{x},\; \low{w}^\text{dp} \leq w^\text{dp} \leq \upp{w}^\text{dp} \}, &\; j=1,\ldots,n, \label{eq:revproporig} \\
  \{ w^\text{dp}_j : h_i^\text{dp}(x,w_{i+1}^\text{dp},\ldots,w_{m^\text{dp}}^\text{dp}) = w_i^\text{dp} :
      \low{x} \leq x \leq \upp{x}, \low{w}^\text{dp} \leq w^\text{dp} \leq \upp{w}^\text{dp} \}, \label{eq:revpropaux} \\
  & \hspace{-2em}j=i+1,\ldots,m^\text{dp},  \nonumber
\end{align}
for $i=1,\ldots,m^\text{dp}$.
The variable bounds $\low{w}^\text{dp}$, $\upp{w}^\text{dp}\in\overline{\R}^{m^\text{dp}}$ are initially set to
$\low{w}^\text{dp}_i=\low{g}_i$, $\upp{w}^\text{dp}_i=\upp{g}_i$, $i=1,\ldots,m$, and
$\low{w}^\text{dp}_i=-\infty$, $\upp{w}^\text{dp}_i=\infty$, $i=m+1,\ldots,m^\text{dp}$.

It is worth noting here that the variables $w^\text{dp}$ are not actually added as SCIP variables, although this has been suggested, but merely serve notational purposes.
In the context of domain propagation, only the bounds $\low{w}^\text{dp}$ and $\upp{w}^\text{dp}$ are relevant and stored in the expression.

For the construction of a linear relaxation, a similar extended formulation is considered:
\begin{equation}
  \label{eq:minlp_extlp}
  \tag{$\text{MINLP}_\text{ext}^\text{lp}$}
  \begin{aligned}
    \min\; & \linobj^\T x, \\
    \mathrm{s.t.}\; & h^\text{lp}_i(x,w^\text{lp}_{i+1},\ldots,w^\text{lp}_{m^\text{lp}}) \lesseqgtr_i w_i^\text{lp}, & i=1,\ldots,m^\text{lp}, \\
    & \low{b} \leq Ax \leq \upp{b}, \\
    & \low{x} \leq x \leq \upp{x}, \\
    & \low{w}^\text{lp} \leq w^\text{lp} \leq \upp{w}^\text{lp}, \\
    & x_\intvarindex \in \Z^{\intvarindex}.
  \end{aligned}
\end{equation}
Functions $h_i^\text{lp}(\cdot)$ are again obtained from the expressions that define functions $g_i(\cdot)$ by recursively replacing subexpressions by auxiliary variables $w_{i+1}^\text{lp},\ldots,w_{m^\text{lp}}^\text{lp}$.
However, it is important to note that different subexpressions may be replaced when setting up $h^\text{lp}(\cdot)$ compared to setting up $h^\text{dp}(\cdot)$.
In fact, in contrast to \eqref{eq:minlp_extdp}, it is assumed that algorithms are available to compute a linear outer-approximation of the sets
\begin{equation}
  \{ (x,w^\text{lp}) : h_i^\text{lp}(x,w_{i+1}^\text{lp},\ldots,w_{m^\text{lp}}^\text{lp}) \lesseqgtr_i w_i^\text{lp}, x\in[\low{x},\upp{x}], w^\text{lp}\in[\low{w}^\text{lp},\upp{w}^\text{lp}] \},\;  i=1,\ldots,m^\text{lp}. \label{eq:extlpcons}
\end{equation}
Thus, the auxiliary variables $w_i^\text{lp}$, $i=m+1,\ldots,m^\text{lp}$, can be different from $w_i^\text{dp}$, $i=m+1,\ldots,m^\text{dp}$.
However, the slack variables $w_i^\text{lp}$, $i=1,\ldots,m$, can be considered as identical to $w_i^\text{dp}$.
Similarly to \eqref{eq:minlp_extdp}, the variable bounds $\low{w}^\text{lp}$, $\upp{w}^\text{lp}\in\smash{\overline{\R}^{m^\text{lp}}}$ are initially set to
$\low{w}^\text{lp}_i=\low{g}_i$, $\upp{w}^\text{lp}_i=\upp{g}_i$, $i=1,\ldots,m$, and
$\low{w}^\text{lp}_i=-\infty$, $\upp{w}^\text{lp}_i=\infty$, $i=m+1,\ldots,m^\text{lp}$.
Regarding the (in)equality sense $\lesseqgtr_i$, a valid simplification would be to assume equality everywhere.
For performance reasons, though, it can be beneficial to relax certain equalities to inequalities if that does not change the feasible space of~\eqref{eq:minlp_extlp} when projected onto $x$.
Therefore,
\[
  \lesseqgtr_i\; \coloneqq \begin{cases}
  =,    & \text{if } \low{g}_i > -\infty,\; \upp{g}_i < \infty, \\
  \leq, & \text{if } \low{g}_i = -\infty,\; \upp{g}_i < \infty, \\
  \geq, & \text{if } \low{g}_i > -\infty,\; \upp{g}_i = \infty,
\end{cases}
\quad\text{ for }i=1,\ldots,m.
\]
For $i>m$, monotonicity of expressions needs to be taken into account.
This is discussed in Section~\ref{sect:nllocks}.

In difference to~\eqref{eq:minlp_extdp}, the variables $w^\text{lp}$ are added to SCIP as variables when the LP is initialized. They are marked as \emph{relaxation-only}~\cite{SCIP7}, that is, are not copied when the SCIP problem is copied and are fixed or deleted when restarting (new auxiliary variables are added for the next SCIP round).

To decide for which constraints in \eqref{eq:minlp_extlp} it can make sense to try to improve their linear relaxation, the value of a subexpression needs to be compared with the value for $h_i^\text{lp}(\cdot)$.
Thus, define $\hat h_i^\text{lp}(x)$ to be the value of the subexpression that $h_i^\text{lp}(\cdot)$ represents if evaluated at $x$.
Formally, for $i=1,\ldots,m^\text{lp}$,
\[
\hat h_i^\text{lp}(x) \defi h_i^\text{lp}(x,w_{i+1}^\text{lp},\ldots,w_{m^\text{lp}}^\text{lp}) 
\;\mathrm{where}\; w_j^\text{lp} \defi h_j^\text{lp}(x,w_{j+1}^\text{lp},\ldots,w_{m^\text{lp}}^\text{lp}), j = i+1,\ldots,m^\text{lp}.
\]
Hence, $\hat h_i \equiv g_i$ for $i=1,\ldots,m$.

\paragraph{Example}
Recall Figure~\ref{fig:exprgraph} and the constraint
\[
  \log(x)^2 + 2\log(x)\,y+y^2 \leq 4.
\]
Using structural detection algorithms (discussed in Section~\ref{sect:detect} below), SCIP may replace $\log(x)$ by an auxiliary variable $w_2$, since that results in a quadratic form $w_2^2+2w_2y+y^2$, which is both bivariate and convex, the former being well suited for domain propagation and the latter being beneficial for linearization.
Therefore, the following extended formulation \eqref{eq:minlp_extdp} may be constructed:
\begin{align*}
h_1^\text{dp}(x,y,w_2^\text{dp}) \defi (w_2^\text{dp})^2 + 2w_2^\text{dp}y + y^2 & = w_1^\text{dp},
 \\
h_2^\text{dp}(x,y) \defi \log(x) & = w_2^\text{dp},
 \\
w_1^\text{dp} & \leq 4.  
\end{align*}
\eqref{eq:minlp_extlp} could be very similar,
\begin{align*}
h_1^\text{lp}(x,y,w_2^\text{lp}) \defi (w_2^\text{lp})^2 + 2w_2^\text{lp}y + y^2 & \leq w_1^\text{lp}, \\
h_2^\text{lp}(x,y) \defi \log(x) & = w_2^\text{lp}, \\
w_1^\text{lp} & \leq 4,
\end{align*}
where equality has been chosen for $h_2^\text{lp}(x,y) = w_2^\text{lp}$ because $(w_2^\text{lp})^2 + 2w_2^\text{lp}y + y^2$ is neither monotonically increasing nor monotonically decreasing in $w_2^\text{lp}$. If, however, $y\geq 0$ and $x\geq 1$, then one may relax to $\log(x) \leq w_2^\text{lp}$.

Next, consider the following slight modification:
\[
  \log(x)^2 + 4\log(x)\,y+y^2 \leq 4.
\]
SCIP may again replace $\log(x)$ by an auxiliary variable $w_2$, since that results in a bivariate quadratic form, but the expression is not convex anymore.
SCIP may therefore decide to introduce additional auxiliary variables to disaggregate the quadratic form for the purpose of constructing a linear relaxation.
Therefore, while~\eqref{eq:minlp_extdp} would be the same as above (with coefficient 2 changed to 4), \eqref{eq:minlp_extlp} would be the result of associating an auxiliary variable with every node of the expression graph:
\begin{align*}
w_2^\text{lp} + 4w_3^\text{lp} + w_4^\text{lp} & \leq w_1^\text{lp}, \\
(w_5^\text{lp})^2 & \leq w_2^\text{lp}, \\
w_5^\text{lp}\,y & \leq w_3^\text{lp}, \\
y^2 & \leq w_4^\text{lp}, \\
\log(x) & = w_5^\text{lp},  \\
w_1^\text{lp} & \leq 4.
\end{align*}



\subsubsection{Variable and Expression Locks}
\label{sect:nllocks}


For constraints that are checked for feasibility,
SCIP asks the constraint handler to add down- and uplocks to the variables in the constraint.
A downlock (uplock) indicates whether decreasing (increasing) the variable could render the constraint infeasible.
While it would be valid to add both down- and uplocks for each variable, more precise information can be useful, for example, for the effectiveness of primal heuristic or dual presolving routines.

For constraints as in \eqref{eq:minlp}, the monotonicity of $g(x)$ and $Ax$ with respect to a specific variable and the finiteness of left- and right-hand sides ($\low{g}$, $\upp{g}$, $\low{b}$, $\upp{b}$) decides which locks should be added.
While for $Ax$ it is sufficient to check the sign of matrix entries, the monotonicity of $g(x)$ can sometimes be deduced by analyzing the expression that defines $g(x)$.
Since monotonicity of $g(x)$ may depend on variable values, variable bounds should be taken into account when deriving monotonicity information and variable locks.

To derive locks for variables, the down- and uplocks for variables are generalized to expressions.
That is, in each expression $e$ a number of down- and uplocks (although they are referred to as negative and positive locks in the code) are stored, which indicate the number of constraints that could become infeasible when the value of $e$ is decreased or increased.
For variable-expressions, these down- and uplocks are then exactly the required down- and uplocks of the corresponding variables.

To start, take a constraint $\smash{\low{g}}_j \leq g_j(x) \leq \smash{\upp{g}}_j$ and assume that the expression that defines $g_j(x)$ is given as $\tilde g(f_1(x),f_2(x),\ldots)$ for some operand $\tilde g$ and (sub)expressions $f_1,f_2,\ldots$.
If $\upp{g}_j < \infty$, then increasing the value of $\tilde g$ could render the constraint infeasible, so an uplock is added to $\tilde g$.
Analogously, if $\low{g}_j > -\infty$, then decreasing the value of $\tilde g$ could render the constraint infeasible, so a downlock is added to $\tilde g$.

Next, these locks are ``propagated'' to the children $f_1,f_2,\ldots$.
First, the monotonicity of $\tilde g$ with respect to a child $f_k$ is checked by use of the \texttt{MONOTONICITY} callback of the expression handler for $\tilde g$.
If $\tilde g$ is monotonically increasing in $f_k$, then increasing $f_k$ could render those constraints infeasible that could become infeasible if $\tilde g$ is increased and decreasing $f_k$ could render those constraints infeasible that could become infeasible when $\tilde g$ is decreased.
Therefore, down- and uplocks stored for $\tilde g$ are added to the down- and uplocks, respectively, of $f_k$.
If $\tilde g$ is monotonically decreasing in $f_k$, then increasing $f_k$ would decrease $\tilde g$ and decreasing $f_k$ would increase $\tilde g$.
Therefore, the downlocks of $\tilde g$ are added to the uplocks of $f_k$ and the uplocks of $\tilde g$ are added to the downlocks of $f_k$.
Finally, if no monotonicity of $\tilde g$ in $f_k$ could be concluded, then the sum of down- and uplocks of $\tilde g$ are added to both the down- and uplocks of $f_k$.

This procedure is applied for all expressions $f_1,f_2,\ldots$ and recursively to their successors.
When a variable expression is encountered, then the down- and uplocks in the variable expression are added to the down- and uplocks of the variable.
Therefore, in difference to linear and many other types of constraints in SCIP, a variable in a single constraint can get several down- or uplocks if it appears several times.

When constraints need to be ``unlocked'', the same procedure is run, but down- and uplocks are subtracted instead of added.
To avoid that, due to tightened variable bounds, different monotonicity information is used when removing locks, the calculated monotonicity information is stored (removed) in an expression when it is locked the first time (unlocked the last time).

\begin{figure}[ht]
  \centering
  \begin{tikzpicture}[scale=0.6, sibling distance=12em, level distance=8em]
    \tikzstyle{expr} += [shape=ellipse, draw, align=center, minimum size = 18pt]

    \node[expr] (plus) {$\sum$ \\ \scriptsize up:1 down:0}
      child[<-] {
        node[expr] (sqr1) {$(\cdot)^2$ \\ \scriptsize up:1 down:0}
        child[<-] {
          node[expr] (log) {$\log$\\ \scriptsize up:0 down:2}
          child[<-] {
            node[expr] (x) {$x$\\ \scriptsize up:0 down:2}
            edge from parent node[left,inner sep=.2em] {\scriptsize INC}
          }
          edge from parent node[left,inner sep=.2em] {\scriptsize DEC}
        }
        edge from parent node[left,inner sep=.2em] {\scriptsize INC}
      }
      child[<-] {
        node[expr] (prod) {$\prod$ \\ \scriptsize up:1 down:0}
        edge from parent node[right,inner sep=.2em] {\scriptsize INC}
      }
      child[<-] {
        node[expr] (sqr2) {$(\cdot)^2$ \\ \scriptsize up:1 down:0}
        child[<-] {
          node[expr] (y) {$y$\\ \scriptsize up:0 down:2}
          edge from parent node[right,inner sep=.2em] {\scriptsize DEC}
        }
        edge from parent node[right,inner sep=.1em] {\scriptsize INC}
      };

    \draw[<-] (prod) to[out=-150,in=40] node[right, pos=0.5] {\scriptsize DEC} (log);
    \draw[<-] (prod) to[out=-30,in=140] node[left, pos=0.5] {\scriptsize DEC} (y);

    \node[yshift=+0.5cm] (t3) at (plus) {};

    \node[yshift=-0.5cm]               (t1) at (x) {};
    \node[yshift=-0.5cm]               (t2) at (y) {};
  \end{tikzpicture}
  \caption{Propagation of locks through expression graph. INC/DEC specifies the monotonicity of parent w.r.t.\ child.}
  \label{fig:locks}
\end{figure}
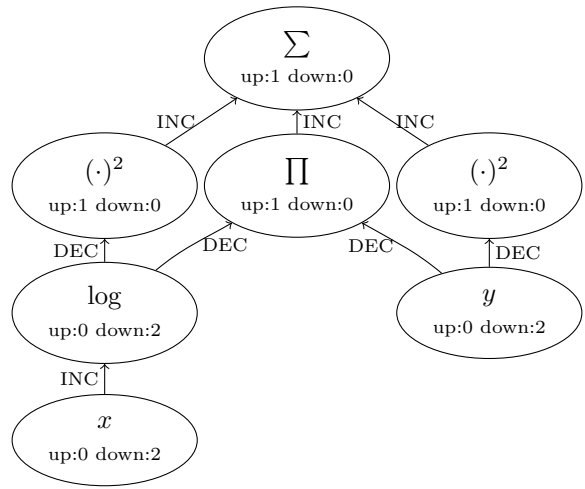

For an example, consider again the expression from Figure~\ref{fig:exprgraph} and the constraint $\log(x)^2 + 2\log(x)\,y+y^2\leq 4$.
Assume further that $\low{x} = 0$, $\upp{x} = 1$, and $\upp{y} = 0$.
The locks for $x$ and $y$ are deduced as follows, see also Figure~\ref{fig:locks}:
\begin{enumerate}
 \item One uplock and no downlock are assigned to the sum-node, because the constraint has a finite right-hand side and no left-hand side.
 \item Since every coefficient in the sum is nonnegative, every child of the sum is assigned one uplock and no downlock.
 \item $\log(x)^2$ is monotonically decreasing in $\log(x)$ because $\log(x)\leq 0$, so that the uplock of $\log(x)^2$ is added to the downlocks of $\log(x)$.
 \item $2\log(x)\,y$ is monotonically decreasing in $\log(x)$ because $2y\leq 0$, so the uplock of $2\log(x)\,y$ is added to the downlocks of $\log(x)$.
 \item $2\log(x)\,y$ is monotonically decreasing in $y$ because $2\log(x)\leq 0$, so the uplock of $2\log(x)\,y$ is added to the downlocks of $y$.
 \item $y^2$ is monotonically decreasing in $y$ because $y\leq 0$, so the uplock of $y^2$ is added to the downlocks of $y$.
 \item $\log(x)$ is monotonically increasing in $x$, so the downlocks of $\log(x)$ are added to the downlocks of $x$.
\end{enumerate}
Thus, eventually both $x$ and $y$ receive 2 downlocks, one for each appearance of the variables in the expression.
Presolvers, primal heuristics, or other plugins of SCIP may now use the information that increasing the value of these variables in any feasible solution does not render this constraint infeasible.

\subsubsection{Nonlinear Handler}
\label{sect:nlhdlr}


The construction of the extended formulations requires algorithms that analyze an expression for specific structures, for instance, quadratic or convex subexpressions as in the previous example.
Following the spirit of the plugin-oriented design of SCIP, these algorithms are not hardcoded into \texttt{cons\_nonlinear}, but are added as separate plugins, referred to as \emph{nonlinear handlers}.
Next to detecting structures in expressions, nonlinear handlers can also provide domain propagation and linear relaxation algorithms that act on these structures.
These plugins have to interact tightly with \texttt{cons\_nonlinear} and nonlinear constraints.
Therefore, in difference to other plugins in SCIP, nonlinear handlers are managed by \texttt{cons\_nonlinear} and not the SCIP core.

In fact, \texttt{cons\_nonlinear} acts both as a handler for nonlinear constraints and as a ``core'' for the management and enforcement of the extended formulations \eqref{eq:minlp_extdp} and \eqref{eq:minlp_extlp}.
As a constraint handler, it checks nonlinear constraints for feasibility, adds them to the NLP relaxation, applies various presolving operations (see Section~\ref{sect:nlpresol}), handles variable locks, and more.
When it comes to domain propagation, separation, and enforcement of nonlinear constraints (see Sections~\ref{sect:nlprop}--\ref{sect:nlenfo}), the constraint handler decides for which constraints in the extended formulations domain propagation or separation should be tried and calls corresponding routines in nonlinear handlers.
When separation fails in enforcement, the constraint handler also selects a branching variable from a list of candidates that has been assembled by nonlinear handlers.

Since domain propagation, separation, and enforcement is partially ``outsourced'' into nonlinear handlers, a certain similarity of nonlinear handler callbacks to constraint handler callbacks is not surprising.
A nonlinear handler can provide the following callbacks:
\begin{itemize}
\item \texttt{COPYHDLR}: include nonlinear handler in another SCIP instance;
\item \texttt{FREEHDLRDATA}: free nonlinear handler data;
\item \texttt{FREEEXPRDATA}: free expression-specific data of nonlinear handler;
\item \texttt{INIT}: initialization;
\item \texttt{EXIT}: deinitialization;
\item \texttt{DETECT}: analyze a given expression ($h_i^\text{dp}(\cdot)$ and/or $h_i^\text{lp}(\cdot)$) for a specific structure and decide whether to contribute in domain propagation for $h_i^\text{dp}(\cdot)=w_i^\text{dp}$ or linear relaxation of $h_i^\text{lp}(\cdot)\lesseqgtr_i w_i^\text{lp}$ (implementation of this callback is mandatory);
\item \texttt{EVALAUX}: evaluate expression with respect to auxiliary variables in descendants, that is, compute $h_i^\text{lp}(x,w_{i+1}^\text{lp},\ldots,w_{m^\text{lp}}^\text{lp})$;
\item \texttt{INTEVAL}: evaluate expression with respect to current bounds on variables, that is, compute an interval enclosure of \eqref{eq:inteval};
\item \texttt{REVERSEPROP}: tighten bounds on descendants, that is, compute interval enclosures of \eqref{eq:revproporig} and \eqref{eq:revpropaux} and update bounds $\low{x}_j$, $\upp{x}_j$, $\low{w}_j$, $\upp{w}_j$ accordingly;
\item \texttt{INITSEPA}: initialize separation data and add initial linearization of \eqref{eq:extlpcons} 
to the LP relaxation;
\item \texttt{EXITSEPA}: deinitialize separation data;
\item \texttt{ENFO}: given a point $(\hat x,\hat w)$, create a bound change or add a cutting plane that separates this point from the feasible set; usually, this routine tries to improve the linear relaxation of $h^\text{lp}_i(x,w^\text{lp}_{i+1},\ldots,w^\text{lp}_{m^\text{lp}}) \lesseqgtr_i w_i^\text{lp}$; if neither a bound change nor a cutting plane was found, register variables for which reducing their domain might help to make separation succeed;
\item \texttt{ESTIMATE}: given a point $(\hat x,\hat w)$, compute a linear under- or overestimator of function $h^\text{lp}_i(x,w^\text{lp}_{i+1},\ldots,w^\text{lp}_{m^\text{lp}})$ that is as tight as possible in $(\hat x,\hat w)$ and valid with respect to either the local or global bounds on $x$ and $w^\text{lp}$; further, register variables for which reducing their domain might help to produce a tighter estimator.
\end{itemize}
More details on the exact input and output of these callbacks is given in the SCIP documentation.

\subsubsection{Constructing Extended Formulations}
\label{sect:detect}


The extended formulations~\eqref{eq:minlp_extdp} and \eqref{eq:minlp_extlp} are constructed simultaneously by processing one nonlinear constraint of \eqref{eq:minlp} at a time (however, common subexpressions that are shared among different constraints are processed only once).
For a constraint $\smash{\low{g}}_i \leq g(x) \leq \smash{\upp{g}}_i$, $i\in\{1,\ldots,m\}$, 
for which domain propagation is enabled (which it is by default), a slack 
variable $w_i^\text{dp}$ and a constraint $h^\text{dp}_i(x) = w_i^\text{dp}$ 
with $h^\text{dp}_i \equiv g_i$ are added to \eqref{eq:minlp_extdp}.
If, additionally, separation or enforcement is enabled (which they are by
default) and SCIP is not in presolve, then the slack variable $w_i^\text{lp}$
and the constraint $h^\text{lp}_i(x) \lesseqgtr_i w_i^\text{lp}$ with 
$h^\text{lp}_i \equiv g_i$ are added\footnote{To be exact, extended 
formulations are not created explicitly and slack variables are not created 
at this state, but the top of the expression $g_i$ in the nonlinear 
constraints is marked for propagation and/or separation. Variables 
$w_i^{\text{lp}}$ are added in SCIP when the LP relaxation is initialized. For 
simplicity, these technicalities are omitted here.}  
to \eqref{eq:minlp_extlp}.
Thereby, $\lesseqgtr_i$ is decided 
according to the finiteness of $\low{g}_i$ and $\upp{g}_i$ as shown above.

Next, the \texttt{DETECT} callback of each nonlinear handler is called for the expression that defines $h^\text{dp}_i$ and $h^\text{lp}_i$.
The callback is informed if domain propagation and/or separation is 
required or if it was already provided by some other nonlinear handler.
The nonlinear handler then analyzes the expression and returns whether it 
wants to participate in (i) domain propagation for $h^\text{dp}_i(x) = 
w_i^\text{dp}$, (ii) separation for $h^\text{lp}_i(x) \leq w_i^\text{lp}$, 
and/or (iii) separation for $h^\text{lp}_i(x) \geq w_i^\text{lp}$.
Furthermore, the nonlinear handler can also introduce auxiliary variables for 
subexpressions, that is, transform $h^\text{dp}_i(x)$ into 
$h^\text{dp}_i(x,w_{i+1}^\text{dp},\ldots)$, $h^\text{lp}_i(x)$ into 
$h^\text{lp}_i(x,w_{i+1}^\text{lp},\ldots)$, and add constraints like 
$h_{i+1}^\text{dp}(x) = w_{i+1}^\text{dp}$ and $h_{i+1}^\text{lp}(x) 
\lesseqgtr_{i+1} w_{i+1}^\text{lp}$.
The sense $\lesseqgtr_{i+1}$ is decided according to the expression locks (see 
Section~\ref{sect:nllocks}) that have been deduced for $h_{i+1}^\text{dp}$:
\[\lesseqgtr_{i+1} \defi
  \begin{cases}
    =, & \text{if there are both down- and uplocks for } h_{i+1}^\text{dp}, \\
    \leq, & \text{if there are no downlocks for } h_{i+1}^\text{dp}, \\
    \geq, & \text{if there are no uplocks for } h_{i+1}^\text{dp}.
  \end{cases}
\]
Note that this definition is consistent with the definition of $\lesseqgtr_j$, $j=1,\ldots,m$, in Section~\ref{sect:nlextform} as downlocks (uplocks) are added to $g_j$ if and only if $\low{g}_j>-\infty$ ($\upp{g}_j<\infty$), see Section~\ref{sect:nllocks}.

After processing $h_i^\text{dp}$ and $h_i^\text{lp}$, subexpressions thereof are processed in a depth-first manner.
Whenever a subexpression is visited that is associated with a constraint $h_{i'}^\text{dp}(x) = w_{i'}^\text{dp}$ or $h_{i'}^\text{lp}(x) \lesseqgtr_{i'} w_{i'}^\text{lp}$, 
the \texttt{DETECT} callback of each nonlinear handler is called again, this time for the expression that defines $h^\text{dp}_{i'}$ or $h^\text{lp}_{i'}$.
This way, extended formulations are built for each constraint $\smash{\low{g}}_i \leq 
g_i(x) \leq \smash{\upp{g}}_i$ in a recursive manner by utilizing the structure 
detection algorithms of all available nonlinear handlers.

For an example, recall again the expression from Figure~\ref{fig:exprgraph} and assume that the initial extended formulations are $\log(x)^2+2\log(x)\,y+y^2 = w_1^\text{dp}$ and $\log(x)^2+2\log(x)\,y+y^2 \geq w_1^\text{lp}$.
A nonlinear handler that inspects the sum-expression may decide that it can 
provide domain propagation for the expression if the $\log$-expression was 
replaced by a variable. Thus, an auxiliary variable $w^\text{dp}_2$ as well as 
a constraint $w^\text{dp}_2 = \log(x)$ are introduced and $h_1^\text{dp}(x) = 
\log(x)^2+2\log(x)\,y+y^2$ is changed to $h_1^\text{dp}(x,w_2^\text{dp}) = 
(w_2^\text{dp})^2+2w_2^\text{dp}y+y^2$.
Similarly, the same or another nonlinear handler may decide that it can provide 
a linear relaxation for the inequality if $\log(x)$ was replaced by a variable. 
It will introduce an auxiliary variable $w^\text{lp}_2$ and a constraint 
$w^\text{lp}_2 = \log(x)$. Then it will change $\log(x)^2+2\log(x)\,y+y^2 \geq 
w_1^\text{lp}$ to $(w_2^\text{lp})^2+2w_2^\text{lp}\,y+y^2 \geq w_1^\text{lp}$.
In addition, the nonlinear handler then indicates that tight bounds for 
$w_2^\text{lp}$ and $y$ are required to compute the linear relaxation. This 
again initiates the introduction of an auxiliary variable $w^\text{dp}_2$ and a 
constraint $w^\text{dp}_2 = \log(x)$, given that they were not existing already.

\subsubsection{Nonlinear Handler ``Default''}
\label{sect:nlhdlrdefault}


To ensure that there always exist routines that can provide domain propagation 
and linear relaxation for an expression, the ``fallback'' nonlinear handler 
\texttt{default} is available.
This nonlinear handler resorts to callbacks of expression handlers (see 
Section~\ref{sect:expr}) to provide the necessary functionalities.
However, while nonlinear handlers are usually meant to handle larger parts of an expression, the methods implemented by the expression handlers are limited to the immediate children of an expression and thus have a rather myopic view on the expression.
Therefore, the \texttt{DETECT} callback of the default nonlinear handler is called with a low priority.
It then decides whether it contributes domain propagation or linear relaxation depending on what other nonlinear handler have declared before. 
If the nonlinear handler decides to contribute, it will introduce auxiliary variables $w^\text{dp}$ and/or $w^\text{lp}$ for all immediate children of the current expression.
The other callbacks of the default nonlinear handler, in particular \texttt{EVALAUX}, \texttt{INTEVAL}, \texttt{REVERSEPROP}, \texttt{ESTIMATE}, can then utilize the corresponding ``myopic'' callbacks of the expression handlers.


\subsubsection{Presolve}
\label{sect:nlpresol}



\paragraph{Simplify}
  The simplify callbacks of expression handlers are called to bring the expressions into a canonical form.
  For example, recursive sums and products are flattened and fixed or aggregated variables are replaced by constants or sums of active variables.
  See the documentation of function \texttt{SCIPsimplifyExpr()} for a more exhaustive list of applied simplifications.

\paragraph{Forbid Multiaggregation}
  For variables that appear nonlinearly, multiaggregation is forbidden.   This aims to prevent that a simple term like $x^3$ will be expanded into a possible long polynomial when $x$ is multiaggregated.
  The linear constraint handler currently does not account for such effects when deciding whether a variable should be multiaggregated.

\paragraph{Common Subexpressions}
  Subexpressions that appear several times are identified and replaced by a single expression.
  This also ensures that every variable is represented by only one variable-expression across all constraints and that for expressions that appear in several nonlinear constraints at most one auxiliary variable is introduced in the extended formulations.
  However, sums that are part of other sums are currently not identified, since in the canonical form no sum can have a sum as a child. The same holds for products.
  The \texttt{HASH} and \texttt{COMPARE} callback of the expression handlers are used to identify common subexpressions.

\paragraph{Scaling}
  For constraints for which the expression is a sum (which it always is if there is a constant factor different from $1.0$), it is ensured that the number of terms with positive coefficients is at least the number of terms with negative coefficients by scaling the constraint with $-1$.
  If there are as many positive as negative coefficients, then it is ensured that the right-hand side is not $+\infty$.
  This canonicalization step can be useful for the next point.

\paragraph{Merge Constraints}
  Nonlinear constraints that share the same expression are merged.

\paragraph{Constraint Upgrading}
  Upgrades to other constraint types are checked.
  Most importantly, nonlinear constraints that are linear after simplification are replaced by constraints that are handled by \texttt{cons\_linear}.
  Further, constraints that can be written as $(x - a_x)\cdot(y - a_y)=0$ with $x$ and $y$ binary variables and $a_x,a_y\in\{0,1\}$ are replaced by setpacking constraints.

\paragraph{Linearization of Binary Products}
  Products of binary variables are linearized.
  This is done in a way that is similar to previous SCIP versions~\cite{VigerskeGleixner2017}, but the consideration of cliques is new:
  \begin{itemize}
  \item In the simplest case, a product $\prod_i x_i$ is replaced by a new variable $z$ and a constraint of type ``and'' is added that models $z = \bigwedge_i x_i$.
    The ``and''-constraint handler will then separate a linearization of this product~\cite{BerHP09}.
  \item Optionally, for a product of only two binary variables, $xy$, the linearization can be added directly as linear constraints ($x\geq z$, $y\geq z$, $x+y \leq 1+z$).
  \item For a product of two binary variables, $xy$, it is checked whether $x$ (or its negation) and $y$ (or its negation) are contained in a common clique.
    Taking this information into account allows for simpler linearizations of $xy$.
    For example, $x$ and $y$ being in a common clique implies $x+y\leq 1$ and thus $xy=0$.
    Analogously, $x + (1-y) \leq 1$ gives $xy=x$, $(1-x)+y\leq 1$ gives $xy=y$, and $(1-x)+(1-y)\leq 1$ gives $xy=x+y-1$.
  \item Replacing every product in a large quadratic term $\sum_{i,j}Q_{ij}x_ix_j$ by a new variable and constraint can increase the problem size enormously. 
    SCIP therefore checks whether there exist sums of the form $x_i \sum_j Q_{ij} x_j$ ($Q_{ij}\neq 0$) with at least 50 terms and replaces them by a single variable $z_i$ and the linearization
    \begin{align*}
      \low{Q} x_i & \leq z_i, \\
      z_i &\leq \upp{Q} x_i, \\
      \low{Q} &\leq \sum_j Q_{ij} x_j - z_i + \low{Q} x_i, \\
      \upp{Q} &\geq \sum_j Q_{ij} x_j - z_i + \upp{Q} x_i,
   \end{align*}
   where $\low{Q} \defi \sum_j \min(0, Q_{ij})$, $\upp{Q} \defi \sum_j \max(0, Q_{ij})$.
   This usually gives a looser LP relaxation as when each product $x_ix_j$ would be replaced individually, but has the advantage that less variables and constraints need to be introduced.
   Variable $z_i$ is marked to be implicit integer if all coefficients $Q_{ij}$ are integer.
   Variables $x_i$ that appear in the highest number of bilinear terms are prioritized.
  \end{itemize}

\paragraph{Identification of Integrality}
  For constraints that can be written as $\sum_i a_if_i(x) + b y = c$, $b\neq 0$, it is checked whether the variable type of $y$ can be changed to \emph{implicitly integer}.
  Storing the information that a continuous variable can take only integer values in a feasible solution can be useful in the solving process, for example, when branching on $y$.
  To change the type of $y$, the following conditions need to be satisfied:
  $y$ is of continuous type, $\frac{a_i}{b}\in\Z$, $\frac{c}{b}\in\Z$, and $f_i(x)\in\Z$ for solutions that satisfy integrality requirements of \eqref{eq:minlp} ($x_\intvarindex \in \Z^{\intvarindex}$).
  To determine the latter, the \texttt{INTEGRALITY} callback of expression handlers is used.

\paragraph{Implicit Discreteness~\cite{HansenJaumardRuizXiong1993}}
  It is checked whether some variables can be restricted to be at one of their bounds.
  At the moment, the method looks for a non-binary variable $x$ with finite bounds, without coefficient in the objective function, and that appears in only one constraint.
  Furthermore, this constraint needs to be a polynomial inequality where $x$ appears only in monomials of the form $c_k x^{2k}$ with $k\in\mathbb{N}$ or other monomials where $x$ has exponent $1$.
  If all coefficients $c_k$ have the same sign, then the constraint function is convex ($c_k>0$) or concave ($c_k<0$) in~$x$.
  Finally, if in addition, the constraint has an infinite right-hand side (when $c_k>0$) or an infinite left-hand side (when $c_k<0$), then $x$ can be restricted to be in $\{\low{x},\upp{x}\}$.
  This is valid because any feasible solution with $x\in(\low{x},\upp{x})$ can be transformed into another feasible solution with same objective function value by moving $x$ to one of its bounds.

  If $\low{x}=0$ and $\upp{x}=1$, then $x$ is transformed into a binary variable.
  Otherwise, a bound disjunction constraint $(x\leq \low{x}) \vee (x\geq \upp{x})$ is added.
  This ``upgrade'' of continuous variables to discrete ones has been shown to be particularly effective for box-QP instances.

\paragraph{Identification of Unlocked Linear Variables}
  Since SCIP supports linear objective functions only, problems with a nonlinear objective function are reformulated by the readers of and interfaces to SCIP into one with a linear objective function ($\min f(x)$ becomes $\min z \,\mathrm{s.t.}\, f(x)\leq z$), 
  To ensure feasibility of such artificial constraints, nonlinear constraints are checked for a variable $x_i$, $i\in\{1,\ldots,n\}$, that appears linearly and which value could be increased or decreased in a solution without the risk of violating other constraints (see also Section~\ref{sect:nllocks}).
  When a solution candidate violates a nonlinear constraint where such a variable $x_i$ has been identified, the constraint handler postprocesses this solution by adjusting the value of $x_i$ such that the constraint becomes feasible.
  This modified solution is then passed on to primal heuristic ``trysol'', which will suggest it to the SCIP core the next time this primal heuristic is run.

\paragraph{Bound Tightening}
  Domain propagation is run (see Section~\ref{sect:nlprop}) to tighten variable bounds and identify redundant or always-infeasible constraints ($g([\low{x},\upp{x}]) \subseteq [\low{g},\upp{g}]$ or $g([\low{x},\upp{x}]) \cap [\low{g},\upp{g}] = \emptyset$).
  The extended formulation \eqref{eq:minlp_extdp} for domain propagation is constructed to make use of the \texttt{INTEVAL} and \texttt{REVERSEPROP} callbacks of nonlinear handlers.
  Further, bounds that are implied by the domain of expressions are enforced, if possible, such as the lower bound for arguments of $\log(x)$ or $x^p$ with $p\not\in\Z$ are set to a small positive value.


\subsubsection{Domain Propagation}
\label{sect:nlprop}


As in previous SCIP versions, \texttt{cons\_nonlinear} implements a feasibility-based bound tightening (FBBT) procedure.
For that, interval arithmetic is used to bound the preimage of each constraint function with respect to the constraint sides, i.e., interval over-estimates are computed for
\[
  \{ x \in [\low{x},\upp{x}] : \smash{\low{g}}_j \leq g_j(x) \leq \smash{\upp{g}}_j \}
\]  
for each $j=1,\ldots,m$.
As it is nontrivial to do so for an arbitrary function, the expression graph and extended formulation \eqref{eq:minlp_extdp} are utilized.
Recall from Section~\ref{sect:nlhdlr} that it is assumed that the nonlinear handlers provide methods to compute interval enclosures of $h_i^\text{dp}(\cdot)$ and its inverse, see~\eqref{eq:inteval}--\eqref{eq:revpropaux}.

On the implementation side, domain propagation consists of one or several forward and backward passes through the expression graph.
In the forward pass, interval enclosures of \eqref{eq:inteval} are computed using the current local bounds on variables $x$. 
The interval enclosures of \eqref{eq:inteval} are used to update the bounds $\low{w}_i^\text{dp}$, $\upp{w}_i^\text{dp}$.
In the backward pass, interval enclosures of \eqref{eq:revproporig} and \eqref{eq:revpropaux} are computed and used to update the bounds $\low{x}_j$, $\upp{x}_j$, $\low{w}_j$, $\upp{w}_j$.
Note that constraint sides are taken into account because, initially, $\low{w}_i^\text{dp} = \low{g}_i$ and $\upp{w}_i^\text{dp} = \upp{g}_i$, $i=1,\ldots,m$.

To only recalculate intervals that may result in a bound tightening, the constraint handler gets notified when the local bounds on an original variable $x_i$ that appears in a nonlinear constraint is changed.
If a bound is tightened, then all nonlinear constraints that contain this variable are marked for propagation.
If a bound is relaxed, however, then the previously computed bounds on auxiliary variables $\low{w}^\text{dp}$, $\upp{w}^\text{dp}$ are marked as invalid.
When the domain propagation routine of the constraint handler is called, the expressions of all constraints that were marked for propagation are processed in a depth-first manner (forward pass).
If for some $i=1,\ldots,m$, the interval enclosure of \eqref{eq:inteval} is not a subset of $[\smash{\low{g}}_i,\smash{\upp{g}}_i]$ (that is, the constraint is not redundant with respect to current variable bounds), then $h_i^\text{dp}$ is queued for backward propagation.
The backward propagation queue is then processed in a breadth-first-order.
Each time the interval enclosure of \eqref{eq:revpropaux} provides a sufficient tightening for $[\low{w}_j^\text{dp},\upp{w}_j^\text{dp}]$, $h_j^\text{dp}$ is appended to the backward propagation queue.
If a bound tightening for some $x_j$ is derived from \eqref{eq:revproporig}, constraints that contain $x_j$ are marked for propagation again, so that another forward pass may start after the current backward pass.


The following mentions a few more subtleties.
\paragraph{Auxiliary Variables in \eqref{eq:minlp_extlp}}
  Recall that the \texttt{DETECT} callback of a nonlinear handler can request bound updates for auxiliary variables in \eqref{eq:minlp_extlp}, see Section~\ref{sect:detect}.
  Thus, if $h_i^\text{dp}$ is not only associated with $w_i^\text{dp}$ but also an auxiliary variable $w_{i'}^\text{lp}$, then the bounds on $w_{i'}^\text{lp}$ are tightened, too.

\paragraph{Reducing Side Effects}
  The bounds computed in a backward pass are stored separately from those computed by the forward pass.
  That is, tightened bounds on $w^\text{dp}$ are not immediately used to compute the bounds on functions that use $w^\text{dp}$.
  Instead, a bound tightening on an auxiliary variable in the backward pass first has to result in a bound change on an original variable $x$, which should then result in tighter bounds computed by the forward pass.
  A reason for this implementation detail is that it is tried to reduce side-effects from the backward propagation in a node of the branch-and-bound tree on the domain propagation in another node of the tree.
  With the current implementation, the domain propagation in a node only depends on the bounds of SCIP variables $x$ and $w^\text{lp}$, but not bounds on $w^\text{dp}$ that were computed by backward propagation in a different part of the tree.

\paragraph{Integrality}
  Integrality information on expressions is taken into account to tighten intervals with fractional bounds to integral values.

\paragraph{Handling Rounding Errors in Variable Bounds and Constraint Sides}
  While the domain propagation in the constraint handler and expression and nonlinear handlers are implemented by using interval arithmetics with outward-rounding, this is not the case for many other parts of SCIP.
  For this reason, variable bounds are relaxed by a small amount when entering the forward pass.
  Since these small relaxations result in overestimates for the intervals of all following computations, several cases deserve a special treatment:
  \begin{itemize}
  \item By default, a bound $b$ is relaxed by $10^{-9}\max(1,\abs{b})$.
  \item If, however, the domain width is small, but the bound itself is large, then relaxing by $10^{-9} \abs{b}$ can have a large impact.
    Therefore, bound relaxation is additionally restricted to $10^{-3}$ times the width of the domain.
  \item Bounds on integer variables (including implicit integer) are not relaxed.
  \item Since integral values, especially 0, often have a special meaning, bounds are not relaxed beyond the next integer value.
  \end{itemize}
  Constraint sides are relaxed by a small amount, too.
  Here, an absolute relaxation of $10^{-9}$ is applied.
  
  Finally, also when updating existing bounds in original or auxiliary variables with newly computed ones, the latter are slightly relaxed if the new interval has a nonzero distance of less than \texttt{numerics/epsilon}=$10^{-9}$ to the existing domain.
  That is, instead of concluding infeasibility for the current subproblem, the variable is fixed to the bound that is closest to the new interval.

\paragraph{Special Case: Redundancy Check}
  When checking whether a constraint can be deleted because it is redundant, it needs to be ensured that the constraint is also satisfied for a solution that violates variable bounds by a small amount.
  Otherwise, a feasibility check for the solution in the original problem can fail.
  Hence, when doing a forward pass for the redundancy check, bound for all unfixed variables are relaxed by the feasibility tolerance of SCIP, independent of the variable type.
  Further, constraint sides are relaxed by the feasibility tolerance as well.

\paragraph{Stopping Criterion}
  In the backward pass, only tightenings that do sufficient progress on the bounds of variables $w^\text{dp}$ and $x$ are usually applied.
  This is to avoid many rounds of bound tightening that do only little progress.
  New bounds are considered sufficiently better than previous ones if the variable gets fixed, the relative improvement on a bound is at least \texttt{numerics/boundstreps}=5\%, or a bound changes sign, i.e., is moved to or beyond zero.

\subsubsection{Initialization of Solve and Relaxations}
\label{sect:nlinitlp}

After presolve, SCIP calls the constraint handlers to initialize their data structures for the branch-and-bound process and to initialize the linear and nonlinear relaxations of the problem.
For \texttt{cons\_nonlinear}, the following operations are performed.

\paragraph{Nonlinear Relaxation}
  For each constraint of \eqref{eq:minlp}, a simple check for convexity and concavity of function $g_i(x)$ on $[\low{x},\upp{x}]$ is done. This uses the \texttt{CURVATURE} callback of the expression handlers.
  The constraint is added to the NLP relaxation of SCIP and the row in the NLP is marked as convex or concave, if possible.
  This information is picked up by other plugins that work on a convex nonlinear relaxation of the problem, for example, \texttt{sepa\_convexproj} and \texttt{prop\_nlobbt}.
  
\paragraph{Extended Formulations}
  The extended formulations \eqref{eq:minlp_extdp} and \eqref{eq:minlp_extlp} are setup.
  That is, the \texttt{DETECT} callback of nonlinear handlers are called on the expressions in the nonlinear constraints to identify structure that can be exploited for domain propagation and linear relaxation, see also Section~\ref{sect:detect}.
  
  Afterwards, the slack- and auxiliary variables $w^\text{lp}$ are added to SCIP and are marked as relaxation-only~\cite{SCIP7}.
  For expressions $h_i^\text{lp}(\cdot)$ that were identified to always have an integral value in a feasible solution (see also Section~\ref{sect:nlpresol}), the type of variable $w_i^\text{lp}$ is set to be implicitly integer instead of continuous.
  
  Finally, bounds $\low{w}^\text{lp},\upp{w}^\text{lp}$ are tightened by running a specialized variant of domain propagation (see Section~\ref{sect:nlprop}).
  In this variant, backward propagation is called for all functions $h_i^\text{dp}(\cdot)$, $i=1,\ldots,m^\text{dp}$.
  This is to ensure that domain information that can not be inferred from bounds on the original variables $x$ is stored in the variable bounds $\low{w}^\text{lp},\upp{w}^\text{lp}$.
  For example, for $\log(w_1^\text{lp})$ with $w_1^\text{lp}=xy$ and $x,y\in[-1,1]$, the bound $w_1^\text{lp}>0$ is implied by the expression itself.
  However, bounds on $x$ and $y$ cannot be tightened such that $xy\geq 0$ is ensured.
  
\paragraph{Linear Relaxation}
  An initial linear relaxations of nonlinear constraints is constructed by calling the \texttt{INITSEPA} of all nonlinear handlers that participate in \eqref{eq:minlp_extlp}.
  
  The ``default'' nonlinear handler computes an initial linear relaxation of~\eqref{eq:extlpcons} by calling the \texttt{ESTIMATE} callback of the expression handler.
  This results in a set of linear under- or overestimators of $h_i^\text{lp}(x,w_{i+1}^\text{lp},\ldots,w_{m^\text{lp}}^\text{lp})$, which are completed to valid hyperplanes by adding auxiliary variable $w_i^\text{lp}$.

\paragraph{Collect Square and Bilinear Terms}
  All expressions in \eqref{eq:minlp_extlp} that are of the form $xy$ or $x^2$ (where $x$ and $y$ can be either original or auxiliary variables) are collected in a data structure that is easy to traverse and search.
  This is used by some plugins that work on bilinear terms (Sections~\ref{sect:nlhdlrbilin}, \ref{sect:rlt}). 

\subsubsection{Separation}
\label{sect:nlsepa}

After SCIP solved the LP relaxation for a node of the branch-and-bound tree, it calls the separator callback of the constraint handlers and separators to check whether a cutting plane that separates the current LP solution $(\hat x, \hat w)$ is available.
For a constraint $\smash{\low{g}}_i \leq g_i(x) \leq \smash{\upp{g}}_i$ of \eqref{eq:minlp} that is violated by $\hat x$, the corresponding extended formulation is checked for separating cutting planes.
During separation only ``strong'' cuts desired, by what cutting planes that are more than just barely violated by $(\hat x,\hat w)$ are meant.
The quantification of ``more than just barely'' is left to the separation algorithm (discussed below and in the following sections).

First, if $h_i^\text{lp}(x,w_{i+1}^\text{lp},\ldots,w_{m^\text{lp}}^\text{lp}) \lesseqgtr_i w_i$ is violated by $(\hat x, \hat w)$, then the nonlinear handlers that registered to contribute to the linear relaxation of this constraint are called.
For a nonlinear handler that implements the \texttt{ENFO} callback, it is left completely to the nonlinear handler to decide how to separate $(\hat x,\hat w)$ from \eqref{eq:extlpcons}.
The callback is also informed that only ``strong'' cuts are desired and candidates for branching are not collected.
If the \texttt{ENFO} callback is not implemented, then the \texttt{ESTIMATE} callback must be implemented.
Thus, a linear under- or overestimator of $h_i^\text{lp}(\cdot)$ is requested from the nonlinear handler and completed to a cutting plane.
The cutting plane is deemed as ``strong'' if the estimator is sufficiently close to the value of $h_i^\text{lp}(\cdot)$ in $(\hat x,\hat w)$.

Formally, assume that $h_i^\text{lp}(\hat x,\hat w_{i+1}^\text{lp},\ldots,\hat w_{m^\text{lp}}^\text{lp}) > \hat w_i^\text{lp}$ and that a nonlinear handler provides a linear underestimator $\ell(x,w_{i+1}^\text{lp},\ldots,w_{m^\text{lp}}^\text{lp})$ of $h_i^\text{lp}(\cdot)$ with respect to current local variable bounds.
If $\ell(\hat x,\hat w_{i+1}^\text{lp},\ldots,\hat w_{m^\text{lp}}^\text{lp}) > \hat w^\text{lp}_{i}$, then $\ell(x,w_{i+1}^\text{lp},\ldots,w_{m^\text{lp}}^\text{lp}) \leq w^\text{lp}_i$ is a cutting plane that separates $(\hat x,\hat w)$ and is valid for the current branch-and-bound node (and it is valid globally if $\ell(\cdot)$ does not depend on local variable bounds).
Further, the cut is regarded as \emph{strong} if
\begin{equation}
\ell(\hat x,\hat w_{i+1}^\text{lp},\ldots,\hat w_{m^\text{lp}}^\text{lp}) \geq \hat w^\text{lp}_i + \alpha
(h_i^\text{lp}(\hat x,\hat w_{i+1}^\text{lp},\ldots,\hat w_{m^\text{lp}}^\text{lp}) - \hat w^\text{lp}_i),
\label{eq:strongcut}
\end{equation}
where $\alpha$ is given by parameter \texttt{constraints/nonlinear/weakcutthreshold} and currently set to $0.2$.
That is, for a strong cut, it is required that it closes at least 20\% of the convexification gap.
Note that if $h_i^\text{lp}(\cdot)$ is convex (and this is also detected by SCIP), then the linear underestimator is typically a linearization of $h_i^\text{lp}(\cdot)$ at $(\hat x,\hat w)$ and thus
$\ell(\hat x,\hat w_{i+1}^\text{lp},\ldots,\hat w_{m^\text{lp}}^\text{lp}) = h_i^\text{lp}(\hat x,\hat w_{i+1}^\text{lp},\ldots,\hat w_{m^\text{lp}}^\text{lp})$.

Once $h_i^\text{lp}(x,w_{i+1}^\text{lp},\ldots,w_{m^\text{lp}}^\text{lp}) \lesseqgtr_i w_i$ has been processed, subexpressions $h_{i'}^\text{lp}$ of $h_i^\text{lp}$ are inspected and the nonlinear handler associated with $h_{i'}^\text{lp}$ are called for separation or linear under-/overestimation.

Again, some more subtleties are discussed next.

\paragraph{Constraints to Separate}
  Separation is not called for every violated nonlinear constraint of \eqref{eq:minlp_extlp}.
  For a subexpression $h_{i'}^\text{lp}(\cdot)$ of $h_i^\text{lp}(\cdot)$ (including $h_i^\text{lp}(\cdot)$ itself), $i\in\{1,\ldots,m\}$, separation is only called if the absolute violation of
$$h_{i'}^\text{lp}(x,w_{i+1}^\text{lp},\ldots,w_{m^\text{lp}}^\text{lp})\lesseqgtr_i w_i$$
is at least a certain factor of the violation of original constraint $\smash{\low{g}}_i\leq g(x)\leq \smash{\upp{g}}_i$.
  This factor is controlled by parameter \texttt{constraints/nonlinear/enfoauxviolfactor} and currently set to 0.01.
  This threshold has been added to prevent the separation for constraints whose violation does not contribute significantly to the violation of original constraints.
  In terms of the first example from Section~\ref{sect:nlextform} ($\log(x)^2+2\log(x)y+y^2\leq 4$), this means that if the violation of $(w_2^\text{lp})^2+2w_2^\text{lp}y+y^2 = w_1^\text{lp}$ is very small in comparison to the violation of $\log(x)^2+2\log(x)y+y^2\leq 4$, then separation for the quadratic equation is suspended until the violation of $\log(x)=w_2^\text{lp}$ has been sufficiently reduced.
  
  Further, separation is skipped for nonlinear constraints of \eqref{eq:minlp_extlp} if their absolute violations is below the feasibility tolerance of SCIP, as no strong cuts are expected in this case.

\paragraph{Cut Cleanup}
  Before a cut is passed to the separation storage of SCIP, its numerical properties are checked and improved, if possible.
  For a cut $\sum_{j=1}^k a_jx_j \leq b$ ($x_j$ refers here to either original or auxiliary variables of \eqref{eq:minlp_extlp}) with $\abs{a_j} > \abs{a_{j+1}}$, $j=1,\ldots,k-1$, $a_k\neq 0$, the following operations are performed:
  \begin{enumerate}
  \item Ensure that the coefficient range $\abs{\frac{a_1}{a_k}}$ is below a certain threshold, which by default is set to $10^7$ (parameter \texttt{separating/maxcoefratiofacrowprep}).
    To achieve this, the procedure tries to eliminate either the term $a_1x_1$ or $a_kx_k$ from the cut by adding the inequality $-a_jx_j \leq -a_j\smash{\low{x}}_j$ if $a_j>0$, or $-a_jx_j \leq -a_j\upp{x}_j$ if $a_j<0$, for $j=1$ or $j=k$. Whether the first or the last term is chosen depends on finiteness of variable bounds and the amount that the cut is relaxed in $\hat x$ by this operation.
    Since local variable bounds are used here, a cut that was previously globally valid may now be locally valid only.
  \item \label{strongcutmalcoef} If the absolute value of the maximal coefficient, $\abs{a_1}$, is below $10^{-4}$ or above $10^4$ (parameter \texttt{constraints/nonlinear/strongcutmaxcoef}), then the cut is scaled by a factor $2^p$, $p\in\mathbb{Z}$, that suffices to ensure $\abs{a_1} \in[10^{-4},10^4]$.
  Thus, together with the previous criterion, this ensures that the absolute value of all coefficients is within $[10^{-4},10^4]$ and aims on sorting out cutting planes whose absolute violation is very small or large due to bad scaling only.
  \item Ensure that fractional coefficients are not within $\varepsilon=10^{-9}$ (\texttt{numerics/epsilon}) of an integer value. That is, for $a_j$ with $0<\abs{a_j-\lfloor a_j\rceil}\leq \varepsilon$, where $\lfloor a_j\rceil$ denotes the integer-rounding of $a_j$, the cut is relaxed by adding the inequality $(\lfloor a_j\rceil-a_j)x_j \leq (\lfloor a_j\rceil-a_j)\upp{x}_j$ if $\lfloor a_j\rceil > a_j$ or $(\lfloor a_j\rceil-a_j)x_j \leq (\lfloor a_j\rceil-a_j)\low{x}_j$ if $\lfloor a_j\rceil < a_j$.
  Due to the use of bound information, a previously globally valid cut may be only locally valid now.
  The main motivation for this operation is to prevent the replacement of $a_j$ by $\lfloor a_j\rceil$ that would occur when the cut is stored in a \texttt{SCIP\_ROW}, since that could make the cut invalid.
  \item Similar to the previous point, a right-hand side $b$ that is very close to $0$ is relaxed. If $b\in[-10^{-9},0]$, then $b$ is relaxed to 0. If $b\in(0,10^{-9}]$, then $b$ is relaxed to $1.1\cdot 10^{-9}$.
  This is done to prevent the replacement of $b\in(0,10^{-9}]$ by $0$ that would occur when a \texttt{SCIP\_ROW} is formed.
  \end{enumerate}
  If the cleanup failed, for example, because finite bounds were not available to relax the cut, the relaxed cut is no longer violated in $\hat x$, or is no longer ``strong'' (\eqref{eq:strongcut} is one way of defining what a ``strong'' cut is), then it is discarded.

\paragraph{Linearization in Incumbents}
  In the last decades, solvers for convex MINLPs have demonstrated that the choice of the reference point in which to linearize convex nonlinear constraints is essential.
  While using the solution of the LP relaxation still leads to a convergent algorithm~\cite{Kelley1960}, better performance is achieved by using a reference point that is close to or at the boundary of the feasible region~\cite{DuranGrossmann1986,SerranoSchwarzGleixner2020}.
  Therefore, also the new implementation of \texttt{cons\_nonlinear} includes a feature where feasible solutions are used as reference points to generate cutting planes.
  
  That is, whenever a primal heuristic finds a new feasible solution $x^*$, SCIP iterates through the nonlinear constraints of~\eqref{eq:minlp_extlp} in reverse order, sets $(w_i^\text{lp})^* \defi h_i^\text{lp}(x^*,(w_{i+1}^\text{lp})^*,\ldots,(w_{m^\text{lp}}^\text{lp})^*)$ and calls the \texttt{ESTIMATE} callback of the registered nonlinear handler (if it implements this callback) with $(x^*,w^*)$ as reference point.
  If a globally valid underestimator $\ell(x,w_{i+1}^\text{lp},\ldots,w_{m^\text{lp}}^\text{lp})$ is returned with
  $\ell(x^*,(w_{i+1}^\text{lp})^*,\ldots,(w_{m^\text{lp}}^\text{lp})^*) = (w_i^\text{lp})^*$ (that is, it is supporting $h_i^\text{lp}(\cdot)$ at $(x^*,w^*)$), then the cut $\ell(x,w_{i+1}^\text{lp},\ldots,w_{m^\text{lp}}^\text{lp}) \leq w_i^\text{lp}$ is added to the cutpool of SCIP.
  Overestimators are handled analogously.
  However, since this feature gave mixed computational results when it was added, it is currently disabled by default (parameter \texttt{constraints/nonlinear/linearizeheursol}).

\subsubsection{Enforcement}
\label{sect:nlenfo}


The enforcement callbacks of constraint handlers are the ones where resolving infeasibility of solutions has to be taken most seriously.
While domain propagation and separation callbacks are allowed to return empty-handed, the enforcement for nonlinear constraint needs to find some action to enforce violated nonlinear constraints in a given solution point.
Especially when points are almost feasible, i.e., when violations in \eqref{eq:minlp_extlp} are small (reconsider also the motivating example from Section~\ref{sect:nlmotivation}), enforcing constraints can be difficult and some measures taken may appear desperate.

In summary, the constraint handler attempts to enforce constraints of \eqref{eq:minlp} by separation on \eqref{eq:minlp_extlp}, domain propagation on \eqref{eq:minlp_extdp}, or branching on a variable $x_i$, $i\in\{1,\ldots,n\}$.

In the unlikely case that no relaxation has been solved, then the constraint handler is asked to enforce the pseudo-solution (\texttt{ENFOPS} callback), that is, a vertex of the variables domain with best objective function value.
In this case, domain propagation is called (see Section~\ref{sect:nlprop}).
If no bound change is found and infeasibility of the node is not concluded, then all variables in all violated nonlinear constraints and with domain width larger than $\varepsilon$ (\texttt{numerics/epsilon}) are registered as branching candidates.
The branching rules of SCIP for external branching candidates will then take care of selecting a variable for branching.
If no branching candidate could be found, then it is not clear whether there is no feasible solution left in the current node (though relevant domains are tiny).
In this case, the constraint handler instructs SCIP to solve the LP relaxation.

When the constraint handler has to enforce a solution $(\hat x,\hat w)$ of the LP relaxation (\texttt{ENFOLP} callback), then the following steps are taken:
\begin{enumerate}
 \item The violation of the solution in \eqref{eq:minlp} and \eqref{eq:minlp_extlp} is analyzed.
   Let $v^g$, $v^h$, and $v^b$, be the maximal absolute violation of the nonlinear constraints in \eqref{eq:minlp}, the nonlinear constraints in \eqref{eq:minlp_extlp}, and the bounds of variables in nonlinear constraints ($\low{x},\upp{x},\low{w},\upp{w}$), respectively.
   Further, let $\text{tol}^\text{feas}$ be the feasibility tolerance of SCIP (\texttt{numerics/feastol}), $\text{tol}^\text{lp}$ be the current primal feasibility tolerance of the LP solver, and $\varepsilon$ be the value of \texttt{numerics/epsilon}.
   By default, $\text{tol}^\text{feas} = \text{tol}^\text{lp} = 10^{-6}$ and $\varepsilon = 10^{-9}$.
   Thus, if $v^g \leq \text{tol}^\text{feas}$, then all nonlinear constraints are satisfied with respect to SCIPs feasibility tolerance and no enforcement is necessary.
   Further, note that SCIP itself already ensures $v^b \leq \text{tol}^\text{lp}$ and $\text{tol}^\text{lp} \leq \text{tol}^\text{feas}$.

 \item \label{enfobndviol} If $v^b > v^h$, that is, the violation of variable bounds is larger than violations of the nonlinear constraints in \eqref{eq:minlp_extlp}, then chances to derive cutting planes from \eqref{eq:minlp_extlp} that separate $(\hat x,\hat w)$ are low.
  This is because methods that work on nonconvex constraints often take variable bounds into account and do not work well when the reference point is outside these bounds.
  Hence, if $v^b > v^h$ and $\text{tol}^\text{lp} > \varepsilon$, then $\text{tol}^\text{lp}$ is reduced to $\max (\varepsilon, v^b/2)$ and a resolve of the LP is triggered.
  
  \item If $v^h < \text{tol}^\text{lp}$, that is, violations of the nonlinear constraints in \eqref{eq:minlp_extlp} are below the feasibility tolerance of the LP solver, then deriving a valid cut that is violated in the current LP solution by more than $\text{tol}^\text{lp}$ can be very difficult.
  Therefore, if also $\text{tol}^\text{lp} > \varepsilon$, then $\text{tol}^\text{lp}$ is reduced to $\max (\varepsilon, v^h/2)$ and a resolve of the LP is triggered.
  
  \item \label{enfosepa} The separation algorithm from Section~\ref{sect:nlsepa} is called with some additional flags that indicate that it is called from the enforcement callback.
  These additional flags extend the separation algorithm as follows.
  \begin{itemize}
    \item When the \texttt{ENFO} or \texttt{ESTIMATE} callbacks of a nonlinear handler are called, then they are instructed to register variables $x_j$ or $w^\text{lp}_i$ for branching, if useful.
    A variable should be registered as a branching candidate if branching on that variable could result in finding tighter cutting planes on the resulting subproblems.
    Usually, this is the case when a convexification gap was introduced due to convexification of a nonconvex function with respect to the current variable domain.
    Thus, nonlinear handler that underestimate convex expressions usually do not register branching candidates.
    
    \item A forward pass of domain propagation in \eqref{eq:minlp_extdp} (see Section~\ref{sect:nlprop}) is run to ensure that recent bound tightenings are taken into account.

    \item Recall that for a violated constraint $\smash{\low{g}}_i\leq g_i(x)\leq \smash{\upp{g}}_i$ with $i \in \{1, \dots, m\}$, constraint $h_i^\text{lp}(x,w_{i+1}^\text{lp},\ldots,w_{m^\text{lp}}^\text{lp}) \lesseqgtr_i w_i$ and subexpressions of $h_i^\text{lp}$ are tried for separation.
      If for none of them a ``strong'' cut could be found, no branching candidate was registered, and the violation of constraint $\smash{\low{g}}_i\leq g_i(x)\leq \smash{\upp{g}}_i$ is at least $0.5\,v^g$ (parameter \texttt{constraints/nonlinear/weakcutminviolfactor}), then separation is repeated without the requirement that cutting planes need to be ``strong''.
      
    \item Dropping the requirement for ``strong'' cuts has various consequences on the separation algorithm described in Section~\ref{sect:nlsepa}:
      The requirement that the absolute violation of constraints of \eqref{eq:minlp_extlp} is at least $\text{tol}^\text{feas}$ is dropped (recall again the motivating example from Section~\ref{sect:nlmotivation} where a solution was feasible with respect to $\text{tol}^\text{feas}$ for \eqref{eq:minlp_extlp} but not feasible for \eqref{eq:minlp}).
      
      Instead of \eqref{eq:strongcut}, it is now sufficient that the violation of the cutting plane in $(\hat x,\hat w)$ is at least $\text{tol}^\text{lp}$.
      
      The cleanup of the cut is modified to take the minimal violation $\text{tol}^\text{lp}$ into account.
      That is, if the violation is in $[10\varepsilon,\text{tol}^\text{lp}]$, then it is scaled up to reach a violation of $10^{-4}$ (parameter \texttt{separating/minefficacy(root)}), if possible\footnote{See implementation of \texttt{SCIPcleanupRowprep()} for more details.}, or at least $\text{tol}^\text{lp}$.
      Step~\ref{strongcutmalcoef} in the original cut cleanup (scale to get coefficients into $[10^{-4},10^4]$) is replaced by scaling down the cut to achieve $\abs{a_1} < 10/\text{tol}^\text{feas}$, if this is possible without the violation to drop under $\text{tol}^\text{lp}$.
      
      Since cuts with violations that are just barely above feasibility tolerances are allowed, it is tried to ensure that floating-point rounding-off errors do not falsify the magnitude of the calculated violation.
      For that, it is required that the violation of the cut $\sum_ja_jx_j\leq b$ is sufficiently large when compared to the terms of the cut, i.e., at least $2^{-50} \max (\abs{b}, \max_j \abs{a_j\hat x_j})$ is required.
      The value $50$ has been chosen because the mantissa of a floating-point number in double precision has $52$ bits.
            
     The cut cleanup procedure is instructed to record for which variables it has modified coefficients in order to achieve the desired coefficient range or to avoid coefficients to be within $\varepsilon$ of an integral value.
     If the cleanup failed to produce a violated cut, then these variables are registered as branching candidates (auxiliary variables may be mapped onto original variables, though, see Section~\ref{sect:nlbranch}).
     The motivation is that since a bound of these variables was used to relax the cut, having a smaller domain may result in less relaxation and thus a higher chance to find a violated cut.

     \paragraph{\small Case Study} {\small
       Most of the ``cut cleanup'' routines have been added to improve numerical stability on test instances.
       One of the more peculiar cases is detailed in the following.
       On instance \texttt{ex1252} from MINLPLib, constraint \texttt{e4} is originally given as $-6.52 (0.00034 x_6)^3 - 0.102 (0.00034x_6)^2 x_{12} + 7.86\cdot 10^{-8}x_{12}^2 x_6 + x_3 = 0$ (coefficients have been rounded).
       After simplification of expressions, this is represented in SCIP as $x_3-2.54\cdot 10^{-10}x_6^3-1.17\cdot 10^{-8}x_6^2x_{12}+7.86\cdot 10^{-8}x_6x_{12}^2 = 0$.
       When~\eqref{eq:minlp_extlp} is constructed, none of the specialized nonlinear handlers detect a structure, so nonlinear handler ``default'' introduces an auxiliary variable for each nonlinear term.
       The resulting constraint\footnote{The attentive reader observes that the constraint handler is partially responsible for its own misery here by naively replacing each monomial by an auxiliary variable. Adding an automated scaling for newly introduced variables or being more considerate in the simplification step may help here.} in~\eqref{eq:minlp_extlp} is
       $x_3 -2.54\cdot 10^{-10}w_{13} -1.17\cdot 10^{-8}w_{14} +7.86\cdot 10^{-8}w_{16} = w_{12}$.
       Assume that in a solution the value in the left-hand side is below the one on the right-hand side.
       Though the constraint is actually linear, enforcing it uses the separation procedures of the nonlinear handler.
       Therefore, the cut that is generated via the help of the expression handler ``sum'' is, not surprisingly,
       $x_3 -2.54\cdot 10^{-10}w_{13} -1.17\cdot 10^{-8}w_{14} +7.86\cdot 10^{-8}w_{16} \geq w_{12}$.
       This cut is marked as globally valid.
       Next, the cut cleanup procedure is run and recognizes that the coefficient range is $\approx 10^{10} > 10^7$.
       It then uses the variable bounds at the current node to eliminate variables from the cut until the coefficient range is sufficiently reduced.
       Apparently, the least relaxation is necessary if the terms for $x_3$, $w_{12}$, and $w_{13}$ are removed.
       The resulting cut, now only valid for the current node, is $7.86\cdot 10^{-8}w_{16} -1.17\cdot 10^{-8}w_{14}\geq -13.94$, which turns out to be no longer violated by the solution to be separated.
       Since the relaxation of the cut used the bounds of $x_3$, $w_{12}$, and $w_{13}$, the only choice left to resolve the violation is to tighten these bounds.
       Therefore, variables $x_3$ and $x_6$ (due to $w_{13}=x_6^3$) are registered as branching candidates (in the current implementation, only the left-hand side of constraints in~\eqref{eq:minlp_extlp} are considered).
       In a later node, the whole procedure repeats, but since variable bounds are tighter, cut cleanup results in the cut $7.86\cdot 10^{-8}w_{16} -1.17\cdot 10^{-8}w_{14}\geq -12.68$, which has a higher chance to be violated.
       Eventually the instance can be solved to a gap below 1\%, but the challenging numerical properties and the costly way they are currently handled take their toll on the performance.}
    \end{itemize}
  
  \item If the separation algorithms were not successful, but branching candidates have been collected, then these candidates are either passed on to the SCIP core as external branching candidates or the branching rules of the nonlinear constraint handler are employed.
    The latter is currently the default (parameter \texttt{constraints/nonlinear/{\allowbreak}branching/{\allowbreak}external}) and described in more detail in Section~\ref{sect:nlbranch} below.
\end{enumerate}
  
In most situations, it is either possible to separate an infeasible solution or to find a variable in a nonconvex term such that branching on that variable should reduce the convexification gap, which would allow for a tighter linear relaxation.
However, enforcement needs to handle also the less likely situations where neither separation nor branching was successful.
This leads to the following (less strategical) attempts.

\begin{enumerate}
  \setcounter{enumi}{5}
  \item If $v^b>\varepsilon$, then $\text{tol}^\text{lp}$ is reduced to $\max(\varepsilon, v^b/2)$ and a resolve of the LP is triggered.
    As in Step~\ref{enfobndviol}, the hope is that separation methods will work better if the LP solution is less outside the variable bounds.

  \item If $v^h>\varepsilon$ and $\text{tol}^\text{lp}>\varepsilon$, then $\text{tol}^\text{lp}$ is reduced to $\max(\varepsilon, v^b/2)$ and a resolve of the LP is triggered.
    The hope here is that i) less tolerance on the feasibility for previously generated cuts may lead to a feasible solution, and ii) more cuts can be added if the minimal required violation is reduced.
  
  \item Domain propagation (Section~\ref{sect:nlprop}) is run in the hope that some bound change that hasn't previously been found in the separation-and-propagation loop for the current node is discovered now.
    This bound change may separate the current LP solution or have an influence on the next separation attempts.
  
  \item Any unfixed variable in violated nonlinear constraints is registered as external branching candidate.
    SCIP then branches on one of these variable and the hope is that infeasibilities in child nodes will be easier to resolve.
    Note that when the domain width of a variable is reduced to less than $\varepsilon$, then the variable is treated as if fixed to a single value.
  
  \item If all variables in violated constraints are fixed, then it may be the overestimation of variable bounds that prevented domain propagation to conclude that the current node is infeasible.
    The node will be cut off and a message issues to the log.
\end{enumerate}

The constraint handler collects statistics on how often it added ``weak'' cuts, tightened the LP feasibility tolerance ($\text{tol}^\text{lp}$ is reset to $\text{tol}^\text{feas}$ whenever processing of a new node starts), branched on any unfixed variable, etc.
The occurrence of such behavior is an indication that SCIP has numerical problems to solve this instance.
To see these statistics, enable parameter \texttt{table/cons\_nonlinear/active}.

Finally, if the constraint handler has to enforce a solution of a relaxation other than the LP (\texttt{ENFORELAX}), then almost the same algorithm is run as for enforcing LP solutions.
The only differences are that i) $\text{tol}^\text{feas}$ is used instead of $\text{tol}^\text{lp}$ as minimal required cut violation and ii) reduction of $\text{tol}^\text{lp}$ is omitted.
Note, that the enforcement of relaxation solutions has not been tested and would probably require some patching up to work reliably.

\subsubsection{Branching}
\label{sect:nlbranch}

The handler for nonlinear constraints now includes its own branching rule to select a variable for branching among a number of candidates.
The candidates are variables that usually appear in nonconvex expressions of violated nonlinear constraints and are collected while trying to find a cutting plane that separates a given relaxation solution (Step~\ref{enfosepa} in the previous section).
Branching on such a variable should reduce the gap that is introduced by convexifying the nonconvex expression in both children because this gap is typically proportional to the domain width.

\paragraph{Mapping Constraint Violation onto Variables}

Within the \texttt{ESTIMATE} and \texttt{ENFO} callbacks of a nonlinear handler, the handler should register with the constraint handler those variables of \eqref{eq:minlp_extlp} where branching could potentially help to produce tighter estimators or cutting planes.
With the branching candidates a ``violation score'' is enclosed, which typically is the relative violation of the nonlinear constraint in \eqref{eq:minlp_extlp} that is currently handled,
\begin{equation}
\label{eq:violscore}
s^v \defi \frac{\vert h_i^\text{lp}(\hat x,\hat w_{i+1}^\text{lp},\ldots,\hat w_{m^\text{lp}}^\text{lp}) - \hat w_i^\text{lp}\vert}{\max(1,\vert\hat w_i^\text{lp}\vert)}
\end{equation}
This value serves as a proxy for the convexification gap associated with $h_i^\text{lp}(\cdot)$.

For each branching candidate, the number of violation scores that have been added, the maximal score, and the sum of scores are stored.
If a nonlinear handler registers only one branching candidate for an expression, then the value $s^v$ can be added to the score of that variable immediately.
For a multivariate function $h_i^\text{lp}(\cdot)$, several candidates may be registered, which requires distributing the $s^v$ onto several variables.
Let $x_{i_1},\ldots,x_{i_k}$, $\{i_1,\ldots,i_k\}\subseteq\varindex$, be such a set of variables.
Let $k_u$ be the number of unbounded variables in this set,
\[
  k_u \defi \abs{\{j\in\{1,\ldots,k\} : \low{x}_{i_j} = -\infty \text{ or } \upp{x}_{i_j} = \infty \}},
\]
If $k_u>0$, then to each unbounded variable an equal part of the violation score is assigned.
That is, variable $x_{i_j}$ is assigned the score
\begin{equation}
\label{eq:varviolscoreunbnd}
\begin{cases}
  \frac{s^v}{k_u}, & \text{if } \low{x}_{i_j} = -\infty \text{ or } \upp{x}_{i_j} = \infty, \\
  0, & \text{otherwise}.
\end{cases}
\end{equation}
Hence, only unbounded variables are considered for branching.
This is because the computation of a linear outer-approximation of \eqref{eq:extlpcons} often depends on the presence of variable bounds.
If all variables are bounded, the following variable weights are considered instead:
\begin{equation}
  \label{eq:midnessweight}
  \lambda_j \defi
    \begin{cases}
      \max\left( 0.05, \frac{\min(\hat x_{i_j} - \low{x}_{i_j}, \upp{x}_{i_j} - \hat x_{i_j})}{\upp{x}_{i_j} - \low{x}_{i_j}} \right), & \text{if } \low{x}_{i_j} \neq \upp{x}_{i_j}, \\
      0, & \text{otherwise},
    \end{cases}
      \qquad j=1,\ldots,k.
\end{equation}
Value $\lambda_j\in[0.05,0.5]$ measures the ``midness'' of the current solution point with respect to the variables domain.
Larger shares of the violation score are then assigned to variables that are closer to the middle of the domain:
\begin{equation}
 \label{eq:varviolscorebnd}
 \frac{\lambda_j}{\sum_{j'=1}^k\lambda_{j'}}s^v.
\end{equation}
This choice is inspired by the observation that the convexification gap is typically smallest at the boundary of the domain.
Further, since a value close to $\hat x_{i_j}$ is typically selected as branching point, this choice prefers variables that lead to children in the branch-and-bound tree with similar domain sizes.
The following alternatives to the weights \eqref{eq:midnessweight} for a bounded unfixed variable $x_{i_j}$ can be chosen (parameter \texttt{constraints/nonlinear/{\allowbreak}branching/violsplit}):
\begin{align*}
  \text{uniform:} \;& 1.0 \\
  \text{domain width:} \;& \upp{x}_{i_j}-\low{x}_{i_j} \\
  \text{logarithmic scale of domain width:} \;&
    \begin{cases}
      10 \log_{10}(\upp{x}_{_j} - \low{x}_{i_j}), & \text{if } \upp{x}_{_j} - \low{x}_{i_j} \geq 10, \\
      \frac{1}{-10\log_{10}(\upp{x}_{_j} - \low{x}_{i_j})}, & \text{if } \upp{x}_{_j} - \low{x}_{i_j} \leq 0.1, \\
      \upp{x}_{_j} - \low{x}_{i_j}, & \text{otherwise}.
    \end{cases}
\end{align*}

\paragraph{Auxiliary Variables}

While the choice of notation in the previous section implied that violation scores would only be distributed onto original variables $x_i$, $i\in\varindex$, it is clear that the same formulas can be used if some or all variables are auxiliary variables of the extended formulation~\eqref{eq:minlp_extlp}.
However, recall that an auxiliary variable $w_i^\text{lp}$ is essentially just a proxy for a subexpression that is defined with respect to original variables and other auxiliary variables $w_{i'}^\text{lp}$, $i'>i$.
Due to this construction, branching on original variables $x_i$ is usually preferred, as this tightens not only the bounds on $x_i$ directly but also the bounds on one or several auxiliary variables implicitly (see Section~\ref{sect:nlprop}).
On the other hand, there may be situations where branching on auxiliary variables could be preferable (after all, such branching could tighten bounds on original variables via domain propagation, too) as it has a more direct effect on the bounds on auxiliary variables.
As we have not come up with an intuitive criterion on when to allow branching on auxiliary variables, currently only the minimal depth required for nodes in the branch-and-bound tree to allow branching on auxiliary variables can be specified (parameter \texttt{constraints/nonlinear/branching/aux}). 
The default is to never branch on auxiliary variables, though.
Therefore, when a nonlinear handler registers a set of variables and a violation score for branching, each auxiliary variable $w_i^\text{lp}$ in this set is replaced by the variables that are appear in $h_i^\text{lp}(x,w_{i+1}^\text{lp},\ldots,w_{m^\text{lp}}^\text{lp})$.
This is repeated until only original variables are left.
The violation score is then distributed among this set of original variables.

If during enforcement, separation failed to find a cut (recall Step~\ref{enfosepa} in Section~\ref{sect:nlenfo}), all variables with an assigned violation score are collected by default.
Optionally, only candidates from constraints which violation is a certain factor of $v^g$ are considered for branching.
However, this factor is by default 0 (parameter \texttt{constraints/nonlinear/{\allowbreak}branching/{\allowbreak}highviolfactor}).

\paragraph{Branching Candidate Scores}

Let $x_{i_1},\ldots,x_{i_k}$, $\{i_1,\ldots,i_k\}\subseteq\varindex$, be the set of branching candidates.
With each candidate, up to five different scores are associated.

The \emph{violation score} was already introduced in the previous paragraph.
When for a variable several violation scores have been added, then currently the sum of these values is used:
\[
s_j^v \defi\, \text{sum of violations scores } \eqref{eq:varviolscoreunbnd} \text{ or } \eqref{eq:varviolscorebnd} \text{ assigned to } x_{i_j},\quad j=1,\ldots,k.
\]
Alternatively, also the average or maximum can be used (parameter \texttt{constraints/{\allowbreak}nonlinear/{\allowbreak}branching/scoreagg}).
Summing up has the effect that variables appearing in several violated constraints are likely to be assigned a larger violation score.

The \emph{pseudo-costs score} mimics the adaptation of pseudo costs to spatial branching as it was introduced by Belotti et al.~\cite{BeLeLiMaWa08} and is implemented also in the SCIP core and \texttt{pscost} branching rule plugin.
Pseudo costs associated with a variable $x_i$ are estimates on the change in the dual bound that is provided by the LP relaxation and that result from branching on variable $x_i$.
For brevity, only a simplified explanation of the calculation of pseudo costs and pseudo-cost branching scores for continuous variables and the default setting for \texttt{branching/lpgainnormalize} is presented here.
Assume that the LP relaxation has been solved in a node of the branch-and-bound tree after branching on variable $x_i$ at branching point $\tilde x_i$.
The pseudo costs $\psi_i^+$ and $\psi_i^-$ store the average change in the LP objective function value normalized by the change in the domain width of $x_i$.
That is, if $\Delta$ is the absolute change in the LP relaxations objective function value in the node $x_i\leq \tilde x_i$ with respect to the parent node, then value $\Delta/(\tilde x_i - \low{x}_i)$ is used to update the pseudo cost $\psi_i^+$.
Alternatively, for the node $x_i\geq \tilde x_i$, $\psi_i^-$ is updated with $\Delta/(\upp{x}_i-\tilde x_i)$.

The pseudo-cost score is a prediction of the dual bound gain that can be expected from branching on variable $x_{i_j}$, $j\in\{1,\ldots,k\}$.
Let $\tilde x_{i_j}\in(\low{x}_{i_j},\upp{x}_{i_j})$ be the branching point (see also end of this section) that would be chosen if branching on $x_{i_j}$.
Then the quantities $\psi_{i_j}^+ (\tilde{x}_{i_j} - \low{x}_{i_j})$ and $\psi_{i_j}^- (\upp{x}_{i_j} - \tilde{x}_{i_j})$ are used to define the pseudo-cost branching score:
\[
s_j^p \defi\;
  \begin{cases}
    \text{n/a}, & \text{if } \low{x}_{i_j} = -\infty \text{ or } \upp{x}_{i_j} = \infty, \text { otherwise} \\
    \psi_{i_j}^+ (\tilde{x}_{i_j} - \low{x}_{i_j}) \cdot \psi_{i_j}^- (\upp{x}_{i_j} - \tilde{x}_{i_j}), & \text{if both } \psi_{i_j}^+ \text{ and } \psi_{i_j}^- \text{ are deemed reliable},\\
    \psi_{i_j}^+ (\tilde{x}_{i_j} - \low{x}_{i_j}), & \text{if only } \psi_{i_j}^+ \text{ is deemed reliable},\\
    \psi_{i_j}^- (\upp{x}_{i_j} - \tilde{x}_{i_j}), & \text{if only } \psi_{i_j}^- \text{ is deemed reliable},\\
    \text{n/a}, & \text{otherwise}
  \end{cases}
\]
A value $\psi_{i_j}^+/\psi_{i_j}^-$ is deemed reliable if it has been updated at least twice (\texttt{constraints/{\allowbreak}nonlinear/{\allowbreak}branching/pscostreliable}).
The pseudo-cost score is not computed for problems with constant objective function ($c=0$ in \eqref{eq:minlp}).

The \emph{domain score} aims on giving preference to variables with a domain that is not very large or very small.
The motivation is that relatively large domains may require many branching operations until their domain is small enough to allow for a useful linear relaxation and branching on relatively small domains may not reduce the convexification gap considerably anymore.
The domain score is therefore largest for domains of width~1 and slowly decreases for larger and smaller domains:
\[
s_j^b \defi\;
  \begin{cases}
    \log_{10}(2\cdot 10^{20} / (\upp{x}_{i_j} - \low{x}_{i_j})), & \text{if } \upp{x}_{i_j} - \low{x}_{i_j} \geq 1, \\
    \log_{10}(2\cdot 10^{20} \max(\varepsilon,\upp{x}_{i_j} - \low{x}_{i_j})), & \text{otherwise},
  \end{cases}
\]
The appearance of $10^{20}$ in this formula is due to the implicit bound of $10^{20}$ (\texttt{numerics/{\allowbreak}infinity}) that SCIP applies to unbounded variables.
Thus, in this formula, $\low{x}_{i_j}$ and $\upp{x}_{i_j}$ should be understood as $\pm 10^{20}$ if at infinity.

The \emph{integrality score} aims on giving preference to variables that are of integral type because the domains of integer branching variables in child nodes does not overlap.
Further, binary variables are preferred over integer variables as branching on a binary variable will fix it in both children.
The score is defined as
\[
s_j^i \defi\;
  \begin{cases}
    1.0, & \text{if } i_j \in \intvarindex, \low{x}_{i_j} = 0, \upp{x}_{i_j} = 1, \\
    0.1, & \text{if } i_j \in \intvarindex, \low{x}_{i_j} \neq 0 \text{ or } \upp{x}_{i_j} \neq 1, \\
    0.01, & \text{if } i_j \in \varindex \setminus \intvarindex, x_{i_j} \text{ has been marked to be implicitly integer}, \\
    0.0, & \text{otherwise}.
  \end{cases}
\]

Finally, the \emph{dual score} is a coarse idea that tries to evaluate the importance of violation scores~\eqref{eq:violscore} from the perspective of the dual bound that the LP relaxation provides.
Assume that for a constraint $h_i^\text{lp}(x,w_{i+1}^\text{lp},\ldots,w_{m^\text{lp}}^\text{lp}) \leq w_i^\text{lp}$ of \eqref{eq:minlp_extlp} a cut $\ell(x,w_{i+1}^\text{lp},\ldots,w_{m^\text{lp}}^\text{lp}) \leq w_i^\text{lp}$, where $\ell(\cdot)$ is a linear underestimator of $h_i^\text{lp}(\cdot)$, was added to the LP.
If $\mu$ denotes the dual variable associated with this cut in the LP, then this cut contributes $\mu (\ell(x,w_{i+1}^\text{lp},\ldots,w_{m^\text{lp}}^\text{lp}) - w_i^\text{lp})$ to the Lagrangian function of the LP relaxation.
If instead of the cut the function $h_i^\text{lp}(\cdot)$ could have been used in the LP, then this would change the value of the Lagrangian function by $\mu (h_i^\text{lp}(x,w_{i+1}^\text{lp},\ldots,w_{m^\text{lp}}^\text{lp}) - \ell(x,w_{i+1}^\text{lp},\ldots,w_{m^\text{lp}}^\text{lp}))$.
Therefore, this product of dual variable and convexification gap is used to evaluate the importance that the linear relaxation of this nonlinear constraint has on the dual bound that is provided by the LP.

In the current experimental implementation, the convexification gap
\[
\abs{h_i^\text{lp}(\hat{\hat x},\hat{\hat w}_{i+1}^\text{lp},\ldots,\allowbreak\hat{\hat w}_{m^\text{lp}}^\text{lp}) - \ell(\hat{\hat x},\hat{\hat w}_{i+1}^\text{lp},\ldots,\hat{\hat w}_{m^\text{lp}}^\text{lp})}
\]
in the LP solution $(\hat{\hat x}, \hat{\hat w})$ at the time the cut is generated is stored together with the cut. (The use of the absolute value is to accommodate overestimators from the case where $\lesseqgtr_i$ is $\geq$).
To compute the dual score $s_j^d$ of a variable $x_{i_j}$, for all rows in the LP that contain $x_{i_j}$ and that were generated from a nonlinear constraint of~\eqref{eq:minlp_extlp}, the quantities $\abs{\hat\mu (h_i^\text{lp}(\hat{\hat x},\hat{\hat w}_{i+1}^\text{lp},\ldots,\hat{\hat w}_{m^\text{lp}}^\text{lp}) - \ell(\hat{\hat x},\hat{\hat w}_{i+1}^\text{lp},\ldots,\hat{\hat w}_{m^\text{lp}}^\text{lp}))}$ are added.
Here, $\hat\mu$ refers to the dual value of the cut in the current LP solution.
The current implementation has a number of disadvantages that will need to be addressed before the dual score could be usable by default.
For example, it would obviously be better to use the convexification gap in the current LP solution instead of $\hat{\hat x}$.
Further, cuts may be defined in terms of auxiliary variables, but branching is done on original variables only.
Thus, the replacement of auxiliary variables by original variables (see paragraph ``Auxiliary Variables'' above) would need to be considered here as well.

In a final step, the scores $s_j^v$, $s_j^p$, $s_j^b$, $s_j^i$, $s_j^d$ are aggregated into a single score for each variable.
For that, weights $\gamma^v$, $\gamma^p$, $\gamma^b$, $\gamma^i$, $\gamma^d$ are used, which can be set by parameters \texttt{constraints/nonlinear/branching/*weight} and default to $\gamma^v = 1.0$, $\gamma^p = 1.0$, $\gamma^b = 0$, $\gamma^i=0.5$, $\gamma^d=0$.
Since the scores can be of different magnitudes, they are scaled by the maximal score in each category.
Thus, let $s_{\max}^v \defi \max_{j=1,\ldots,k}s_j^v$ and similar for $s_{\max}^p$, $s_{\max}^b$, $s_{\max}^i$, $s_{\max}^d$.
Further, the case that pseudo-cost scores may not be available for each variable needs to be considered.
Therefore, for a variable where pseudo-cost scores are available, the final score is computed as
\[
  s_j^f \defi\; \frac{\gamma^v \frac{s_j^v}{s_{\max}^v}
                   +\gamma^p \frac{s_j^p}{s_{\max}^p}
                   +\gamma^b \frac{s_j^b}{s_{\max}^b}
                   +\gamma^i \frac{s_j^i}{s_{\max}^i}
                   +\gamma^d \frac{s_j^d}{s_{\max}^d}}{\gamma^v + \gamma^p + \gamma^b + \gamma^i + \gamma^d}.
\]
If a pseudo-cost score is not available, then the other scores are magnified:
\[
  s_j^f \defi\; \frac{\gamma^v \frac{s_j^v}{s_{\max}^v}
                   +\gamma^b \frac{s_j^b}{s_{\max}^b}
                   +\gamma^i \frac{s_j^i}{s_{\max}^i}
                   +\gamma^d \frac{s_j^d}{s_{\max}^d}}{\gamma^v + \gamma^b + \gamma^i + \gamma^d}.
\]

\paragraph{Branching Variable and Coordinate}

Since the variable scores are rather a heuristic guideline than a clear indication which variable is ``best'', the code chooses from all variable with final score at least $0.9\max_{j=1,\ldots,k}s_j^f$ (\texttt{constraints/nonlinear/branching/{\allowbreak}highscorefactor}) uniformly at random.
This allows to exploits performance variability due to branching decisions by changing the seed for the random number generator (\texttt{randomization/randomseedshift}).

The branching point selection rule has not been changed since the last SCIP release.
For a bounded variable $x_j$, a value $\tilde x_j$ between $\hat x_j$ and $\frac{1}{2}(\low{x}_j+\upp{x}_j)$ is chosen, see also Section~4.4.5 of the SCIP Optimization Suite 7.0 release report~\cite{SCIP7}.
Two child nodes are created, one with $x_j \leq \tilde x_j$ and another with $x_j\geq \tilde x_j$, if $j\not\in\intvarindex$.
For $j\in \intvarindex$, domains are ensured to be disjoint ($x_j \leq \lfloor x_j\rfloor$, $x_j \geq \lfloor x_j\rfloor+1$).


\subsection{Nonlinear Handler for Quadratic Expressions}
\label{sect:nlhdlrquad}

The quadratic nonlinear handler detects quadratic expressions, provides specialized domain propagation,
and generates intersection cuts.

\subsubsection{Detection of Quadratic Expressions}
\label{sect:nlhdlrquaddetect}

An expression in a constraint of \eqref{eq:minlp_extdp} or \eqref{eq:minlp_extlp} is recognized as quadratic if it is a sum of terms where at least one term is either a product expression of two expressions or a power expression with exponent $2$.
Formally, the detection routine checks whether the expression can be written as
\begin{equation}
  \label{eq:quadexpr}
  q(y) = \sum_{i=1}^k q_i(y) \quad\text{with}\quad q_i(y) = a_iy_i^2 + c_iy_i + \sum_{j\in P_i} b_{i,j}y_iy_j
\end{equation}
where $y_i$ is either an original variable ($x$) or another expression, $a_i,c_i\in\R$, $b_{i,j}\in\R\setminus\{0\}$, $j\in P_i \Rightarrow i\not\in P_j$ for all $j\in P_i$, $P_i\subset\{1,\ldots,k\}$, $i=1,\ldots,k$.
A bilinear term $y_iy_j$ is associated with $y_i$ (i.e., $j\in P_i$) if $y_i$ appears more often in $q(y)$ than $y_j$.
This is a heuristic choice that should be beneficial for domain propagation.
In case of a tie, the order of expressions (see Section~\ref{sect:expr}) is used as tie-breaker.
If $q(y)$ is linear in $y$ or consists of one power or product term only, then detection is aborted.

After a quadratic structure~\eqref{eq:quadexpr} has been established, the construction of the extended formulations for \eqref{eq:minlp_extdp} and \eqref{eq:minlp_extlp} can differ.

For domain propagation (\eqref{eq:minlp_extdp}), the nonlinear handler checks further whether the quadratic expression is \textit{propagable}, by what is meant that the termwise domain propagation does not necessarily yield the best possible results due to suffering from the so-called \textit{dependency problem} of interval arithmetics.
Specifically, \eqref{eq:quadexpr} is propagable if at least one argument ($y_i$) appears at least twice.
For instance, $x^2 + y^2$ is not propagable, but $x^2 + x$ is.
Only if \eqref{eq:quadexpr} is propagable, the nonlinear handler registers itself as responsible for domain propagation.
Otherwise, the default nonlinear and the expression handlers for sum, product, and power will take care of a termwise propagation of domain.

To construct~\eqref{eq:minlp_extdp}, the nonlinear handler requests an auxiliary variable ($w^\text{dp}$) for any $y_i$ that is an expression and not yet a variable, but with two notable exceptions.
If a variable $y_i$ appears only in a square term of~\eqref{eq:quadexpr} ($a_i\neq 0$, $c_i=0$, $i\not \in P_{i'}$ for all $i'=1,\ldots,k$), then an auxiliary variable is introduced for $y_i^2$ instead of $y_i$.
Similarly, if two variable $y_i$ and $y_j$ appear only as one bilinear term $y_iy_j$ ($a_i=0$, $a_j=0$, $c_i=0$, $c_j=0$, $P_i=\{j\}$ or $P_j=\{i\}$), then an auxiliary variable is introduced for $y_iy_j$ instead of $y_i$ and $y_j$.
That is, a non-propagable part of~\eqref{eq:quadexpr} is split off and treated as if linear since this part does not suffer from the dependency-problem and sometimes better domain propagation routines are available for the single terms $y_i^2$ or $y_iy_j$ (for an example see the bilinear nonlinear handler described in Section~\ref{sect:nlhdlrbilin}).
For an example, consider $xy + z^2 + z$, which is propagable because $z$ appears twice.
However, for~\eqref{eq:minlp_extdp}, the reformulation $w + z^2 + z$, $w=xy$, is applied.
The quadratic nonlinear handler then handles domain propagation for $w + z^2 + z$, while either the default or the bilinear nonlinear handler handle domain propagation for $w=xy$.
An additional advantage of this division of work is that for other expression where $xy$ appears, variable $w$ and its domain information can be reused.

For separation (\eqref{eq:minlp_extlp}), the nonlinear handler registers itself for participation if intersection cuts are enabled (\texttt{nlhdlr/quadratic/useintersectioncuts}, currently disabled by default), no other nonlinear handler (for example the SOC nonlinear handler) handles separation yet, and the corresponding constraint in \eqref{eq:minlp_extlp} is nonconvex.
To decide the latter, the eigenvalues and eigenvectors of the quadratic coefficients matrix (defined by $a_i$ and $b_{i,j}$ of \eqref{eq:quadexpr}) are calculated via LAPACK and stored for later use.
To construct~\eqref{eq:minlp_extlp}, the nonlinear handler requests an auxiliary variable ($w^\text{lp}$) for any $y_i$ that is an expression and not yet a variable.
Thus, even when the quadratic nonlinear handler participates in both domain propagation and separation, the created extended formulations may differ if parts of~\eqref{eq:quadexpr} are not propagable.
This flexibility is a feature of the current design.

A nonlinear handler can choose whether it solely wants to be responsible for domain propagation or separation, or only wants to participate in addition to other routines.
The separation by the quadratic nonlinear handler is such a case, i.e., the nonlinear handler informs to the constraint handler that other possible nonlinear handlers should also be requested for separation.
Currently, this means that the default and bilinear nonlinear handler will become active, auxiliary variables will be introduced for each square and bilinear term, and corresponding under- and overestimators will be computed by these routines if an intersection cut was not generated.
These nonlinear handlers are also the only ones that register branching candidates.
For intersection cuts, bound information is not used explicitly by default
and the quadratic nonlinear handler does not register 
variables for branching.

\subsubsection{Propagation of Quadratic Expressions}

The goal of domain propagation is to use existing bounds on $y$ and $q(y)$ in \eqref{eq:quadexpr} to derive possibly tighter bounds on $q(y)$ and $y$, respectively.
The implementation is similar to the one of \texttt{cons\_quadratic} in SCIP 7 and before~\cite{VigerskeGleixner2017}, but backward propagation has been extended.
For simplicity, the special treatment for some square or bilinear terms as mentioned in the previous section is disregarded here.

\paragraph{Forward Propagation (\texttt{INTEVAL}, \eqref{eq:inteval})}
Here, the goal is to propagate bounds on $y$ to the quadratic expression $q(y)$.
To obtain the best bounds, it would be necessary to maximize/minimize $q(y)$ for $y\in[\low{y}, \upp{y}]$.
Since this is too expensive in general, bounds are overestimated.
For that, each quadratic term $q_i(y)$, $i=1,\ldots,k$, is considered separately.
If $q_i(y)$ is not univariate ($P_i\neq\emptyset$), then each $y_j$, $j\in P_i$, is replaced by its current bounds and the min/max of $a_iy_i^2+(c_i+\sum_{j\in P_i}b_{i,j}[\low{y}_j,\upp{y}_j])y_i$ are calculated~\cite{DomesNeumaier2010}.

\paragraph{Backward Propagation (\texttt{REVERSEPROP}, \eqref{eq:revproporig}, \eqref{eq:revpropaux})}
The goal of backward propagation is to derive bounds on each $y_i$ from given bounds $[\low{q},\upp{q}]$ on $q(y)$ and current bounds on $y$. Similar to forward propagation, the tightest bounds can be derived by maximizing/minimizing each $y_i$ w.r.t.\ $y\in [\low{y},\upp{y}]$ and $q(y)\in [\low{q},\upp{q}]$.
Since this is again too expensive to do in general, bounds are overestimated again.

Let $[\low{q}_i, \upp{q}_i]$, $i=1,\ldots,k$, be bounds on each $q_i(y)$ for $y\in[\low{y},\upp{y}]$.
These are computed in forward propagation.
Similar to forward propagation, bounds for each $y_i$ can be computed by reduction to and solving of a univariate quadratic interval equation~\cite{DomesNeumaier2010}:
\begin{equation}
  \label{eq:quadrevprop}
  a_iy_i^2 + (c_i+\sum_{j\in P_i}b_{i,j}[\low{y}_j,\upp{y}_j])y_i \in [\low{q},\upp{q}] - \sum_{i'=1,i'\neq i}^k [\low{q}_{i'},\upp{q}_{i'}].
\end{equation}

A downside of this approach is that bounds for variables that appear less often may not be deduced.
For example, consider $y_1^2 + y_1y_2 + y_1y_3 + y_2y_3 + y_3$.
As $y_2$ has less appearances than $y_1$ and $y_3$, this quadratic gets partitioned into $q_1(y) = y_1^2 + y_1y_2 + y_1y_3$, $q_2(y)=0$, and $q_3(y) = (y_3 + y_3y_2)$.
Therefore, no bounds are computed for $y_2$ in backward propagation.
The quadratic constraint handler of SCIP 7 handled the case of $q_i \equiv 0$ in certain situations where $y_i$ appeared in only one bilinear term.
For SCIP 8, this has been generalized.
In the example, a bound on $y_2$ is obtained by rewriting as $y_2 + y_3\in ([\low{q},\upp{q}] - [\low{q}_3,\upp{q}_3])/y_1 - y_1$, finding the min/max of the function on the right-hand side, and using this interval for backward propagation on $y_2+y_3$.
In general, after solving~\eqref{eq:quadrevprop}, the quadratic equation $q(y)\in [\low{q},\upp{q}]$ is interpreted as
\[
  c_i + \sum_{j\in P_i}b_{i,j}y_j \in \frac{1}{y_i}\left([\low{q},\upp{q}] - \sum_{i'=1,i'\neq i}^k [\low{q}_{i'},\upp{q}_{i'}]\right) - a_iy_i,
\]
the min/max of the univariate interval function on the right-hand side are calculated, and the resulting interval is used for backward propagation on $c_i + \sum_{j\in P_i}b_{i,j}y_j$.

\subsubsection{Intersection Cuts for Quadratic Constraints}
\label{sect:intercuts}

For separation, assume the constraint of~\eqref{eq:minlp_extlp} is $q(y)\leq w$ with $q(y)$ as in \eqref{eq:quadexpr} and $w$ an auxiliary variable of \eqref{eq:minlp_extlp}.
Further, assume that $q(y)$ is nonconvex ($q(y)$ being convex is handled by the nonlinear handler for convex expressions, see Section~\ref{sect:nlhdlrconvex}).
The quadratic nonlinear handler implements the separation of intersection cuts~\cite{Tuy1964,Balas1971,Glover1973} for the set $S \defi \{ (y,w) \in \R^k : q(y) \leq w \}$ that is defined by this constraint.

Let $(\hat{y},\hat w)$ be a basic feasible LP solution violating $q(y) \leq w$.
First, a convex inequality $g(y,w) < 0$ is build that is satisfied by $(\hat{y},\hat w)$, but by no point of $S$.
This defines a so-called \emph{$S$-free set} $C = \{ (y,w) \in \R^{k+1} : g(y,w) \leq 0 \}$, that is, a convex set 
with $(\hat{y},\hat w) \in \text{int}(C)$ containing no point of $S$ in its interior. 
The quality of the resulting cut highly depends on which $S$-free set is used. 
The tightest possible intersection cuts are obtained by using
\textit{maximal} $S$-free sets as proposed by Muñoz and 
Serrano~\cite{MunozSerrano2020}.

By using the conic relaxation $K$ of the LP-feasible region defined by the 
nonbasic variables at $(\hat{y},\hat w)$, the intersection points between 
the extreme rays of $K$ and the boundary of $C$ are computed. The intersection cut is then 
defined by the hyperplane going through these points and successfully separates 
$(\hat{x},\hat w)$ and $S$. Adding this cut to the LP relaxation excludes the violating 
point $(\hat{x},\hat w)$ from the LP-feasible region and thus enforces the quadratic 
constraint $q(y) \leq w$. To obtain even better cuts, there is also a 
strengthening procedure implemented that uses the idea of negative edge 
extension of the cone $K$~\cite{Glover1974}.
A detailed description of how the (strengthened) intersection cuts are
implemented can be found in the paper by Chmiela et al.~\cite{ChmielaMunozSerrano2021}.

\subsection{Nonlinear Handler for Second-Order Cones}
\label{sect:nlhdlrsoc}


The nonlinear handler for second-order cone (SOC) structures replaces and extends the previous constraint handler for second-order cone constraints.
It detects second-order cone constraints in the original or extended formulation and provides separation by means of a disaggregated cone reformulation.


\subsubsection{Detection}
Given a constraint $h^\text{lp}_i(x) \leq w_i^\text{lp}$ ($\lesseqgtr_i$ being $\geq$ is handled similarly) of the extended formulation~\eqref{eq:minlp_extlp}, the nonlinear handler checks for four different structures.
In some cases, it distinguishes between constraints that are copies of an original constraint with slack variable $w_i^\text{lp}$ added, that is, $i\leq m$, and constraints that only exist in the extended formulation due to the introduction of auxiliary variables, i.e., $i>m$.
Further, binary variables are treated as if they were squared, since this increases the likelihood of finding a SOC-structure.

\paragraph{Euclidian Norm}
If $i>m$, it is checked whether $h^\text{lp}_i(x)$ has the form
\begin{equation}
 \label{eq:socnorm1}
 \sqrt{\sum_{j=1}^k (a_j y_j^2 + b_j y_j) + c}
\end{equation}
for some coefficients $a_j,b_j,c\in\R$, $a_j>0$, and where $y_j$ is either an original variable ($x$) or some subexpression of $h^\text{lp}_i(\cdot)$, $j=1,\ldots,k$, for some $k\geq 2$.
Rewriting~\eqref{eq:socnorm1} reveals the constraint
\begin{equation}
 \label{eq:socnorm2}
 \sqrt{\sum_{j=1}^k \left(\left(\sqrt{a_j} y_j + \frac{b_j}{2\sqrt{a_j}}\right)^2 - \frac{b_j^2}{4a_j}\right) + c} \leq w_i^\text{lp}.
\end{equation}
If $c-\sum_{j=1}^k \frac{b_j^2}{4a_j}\geq 0$, then~\eqref{eq:socnorm2} has SOC-structure.
Thus, the nonlinear handlers requests auxiliary variables for each $y_j$, $j=1,\ldots,k$, and declares that it will provide separation.
In a future version, any positive-semidefinite quadratic expression should be allowed for the argument of $\sqrt{\cdot}$ in~\eqref{eq:socnorm1}.

If $i\leq m$, then $w_i^\text{lp}$ is just a slack-variable and the constraint is equivalent to $\sqrt{\sum_{j=1}^k (a_jy_j^2 + b_jy_j)+c} \leq \upp{g}_i$.
In that case, the nonlinear handler does not get active.
Assuming $a_j>0$ again, this will result in the extended formulation
$\sqrt{w_0} \leq \upp{g}_i$, $\sum_{j=1}^k(a_jw_j+b_jy_j+c) \leq w_0$, $y_j^2 \leq w_j$, $j=1,\ldots,k$, where $w_0,\ldots,w_k$ are new auxiliary variables.
We believe that separation for this formulation will be more efficient than for~\eqref{eq:socnorm2} (constraint $\sqrt{w_0}\leq \upp{g}_i$ is easily enforced by domain propagation).

\paragraph{Simple Quadratics}
Check whether $h^\text{lp}_i(x)$ has the form
\[
 \sum_{j=1}^k (a_jy_j^2) - a_{k+1}y_{k+1}^2 + c
\]
where $a_j,c\in\R$, $a_j>0$, and $y_j$ is either an original variable ($x$) or some subexpression of $h^\text{lp}_i(\cdot)$, $j=1,\ldots,k+1$, for some $k\geq 1$.
The relaxed constraint
\begin{equation}
 \label{eq:socsimple1}
 \sum_{j=1}^k (a_jy_j^2) - a_{k+1}y_{k+1}^2 + c \leq \upp{w}_i^\text{lp}
\end{equation}
has SOC-structure if $c-\upp{w}_i^\text{lp}\geq 0$.
Thus, in this case the nonlinear handler requests auxiliary variables for each $y_j$, $j=1,\ldots,k+1$, and declares that it will provide separation.

If $i\leq m$, then the replacement of the slack variable $w_i^\text{lp}$ by $\upp{w}_i^\text{lp}$ will not be problematic since it is sufficient to enforce the original constraint $g_i(x) \leq \upp{g}_i$ (recall $h_i^\text{lp}=g_i$, $\upp{w}_i^\text{lp} = \upp{g}_i$ initially).
However, if $i>m$, then relaxing $w_i^\text{lp}$ to $\upp{w}_i^\text{lp}$ could mean that infeasibility in~\eqref{eq:minlp} cannot be resolved by enforcing~\eqref{eq:socsimple1}.
Therefore, if $i>m$, then the nonlinear handler indicates to the constraint handler that separation should be requested from other nonlinear handlers as well.
In the current configuration, this introduces auxiliary variables for each square term in~\eqref{eq:socsimple1} by the default nonlinear handler.
The same distinction into $i\leq m$ and $i>m$ applies to the following two structures.

\paragraph{Simple Quadratics (Rotated SOC Variant)}
Check whether $h^\text{lp}_i(x)$ has the form
\[
 \sum_{j=1}^k (a_jy_j^2) - a_{k+1}y_{k+1}y_{k+2} + c
\]
where $a_j,c\in\R$, $a_j>0$, and $y_j$ is either an original variable ($x$) or some subexpression of $h^\text{lp}_i(\cdot)$, $j=1,\ldots,k+2$, for some $k\geq 0$.
The relaxed constraint
\[
 \sum_{j=1}^k (a_jy_j^2) + \frac{a_{k+1}}{4}(y_{k+1}-y_{k+2})^2 - \frac{a_{k+1}}{4}(y_{k+1}+y_{k+2})^2 + c \leq \upp{w}_i^\text{lp}
\]
has SOC-structure if $c-\upp{w}_i^\text{lp}\geq 0$.
Thus, in this case the nonlinear handler requests auxiliary variables for each $y_j$, $j=1,\ldots,k+2$, and declares that it will provide separation.

\paragraph{General Quadratics}
Check whether $h^\text{lp}_i(x) \leq \upp{w}_i^\text{lp}$ is a quadratic constraint that is SOC-representable.
As suggested by Mahajan and Munson~\cite{MahajanMunson2010}, this is done by computing an eigenvalue-decomposition of the quadratic coefficients matrix via LAPACK and attempting to rewrite $h^\text{lp}_i(x)$ as
\[
 \sum_{j=1}^{k+1} \lambda_j (v_j^\T y + \beta_j)^2 + c 
\]
where $\lambda_{1},\ldots,\lambda_{k+1}\in\R$ are the nonzero eigenvalues with corresponding eigenvectors $v_1,\ldots,v_{k+1}\in\R^{\ell}$, $\lambda_j>0$, $j=1,\ldots,k$, $\lambda_{k+1}<0$, $\beta_j,c\in\R$, $j=1,\ldots,{k+1}$, $k\geq 1$, and $y_j$ is either an original variable ($x$) or some subexpression of $h^\text{lp}_i(\cdot)$, $j=1,\ldots,\ell$.
If $0\not\in v_{k+1}^\T [\low{y},\upp{y}]+\beta_{k+1}$ (thus $\sqrt{(v_{k+1}^\T y + \beta_{k+1})^2} = \pm (v_{k+1}^\T y + \beta_{k+1})$ is linear) and $c-\upp{w}_i^\text{lp}\geq 0$, then an SOC-structure has been detected, the nonlinear handler requests auxiliary variables for each $y_j$, $j=1,\ldots,k$, and declares that it will provide separation.



\subsubsection{Separation}
\label{sect:nlhdlrsocsepa}

The SOC constraint that has been detected before is stored in the form
\begin{equation}
 \label{eq:soc}
 \sqrt{\sum_{j=1}^k (v_j^\T y + \beta_j)^2} \leq v_{k+1}^\T y + \beta_{k+1}
\end{equation}
with $v_j\in\R^\ell$, $j=1,\ldots,k+1$, where $y_1,\ldots,y_\ell$ are variables of~\eqref{eq:minlp_extlp}.

Since the left-hand side of~\eqref{eq:soc} is convex, a solution $\hat y$ that violates~\eqref{eq:soc} can be separated by linearization of the left-hand side of~\eqref{eq:soc} via a gradient cut:
\[
  \sqrt{\sum_{j=1}^k (v_j^\T \hat y + \beta_j)^2} +
  \sum_{i=1}^\ell \frac{\partial}{\partial y_i}\sqrt{\sum_{j=1}^k (v_j^\T y + \beta_j)^2}\, (y_i - \hat y_i)
  \leq v_{k+1}^\T y + \beta_{k+1}
\]

However, if there are many terms on the left-hand side of~\eqref{eq:soc} ($k$ being large), then it can require many cuts to provide a tight linear relaxation of~\eqref{eq:soc}.
Thus, as suggested by Vielma et al.~\cite{Vielma2016}, 
a disaggregation of the cone is used if $k\geq 3$:
\begin{align}
 (v_j^\T y + \beta_j)^2 & \leq z_j (v_{k+1}^\T y + \beta_{k+1}), \quad j=1,\ldots,k, \label{eq:socext1} \\
 \sum_{j=1}^k z_j & \leq v_{k+1}^\T y + \beta_{k+1}, \label{eq:socext2}
\end{align}
where variables $z_1,\ldots,z_k$ are new variables that are added to SCIP and marked as ``relaxation-only''.
A solution $(\hat y,\hat z)$ that violates~\eqref{eq:soc} needs to violate also~\eqref{eq:socext1} for some $j\in\{1,\ldots,k\}$ or \eqref{eq:socext2}.
The latter is already linear and can be added as cut.
If a rotated second-order cone constraint~\eqref{eq:socext1} is violated from some $j$, then it is transformed into the standard form
\[
 \sqrt{4(v_j^\T y + \beta_j)^2 + (v_{k+1}^\T y + \beta_{k+1} - z_j)^2} \leq v_{k+1}^\T y + \beta_{k+1} + z_j
\]
and a gradient cut is constructed by linearization of the left-hand side.

If the constraint handler requests only ``strong'' cuts (see Section~\ref{sect:nlsepa}), then gradient cuts are only added when their efficacy is at least $10^{-5}$ (\texttt{nlhdlr/soc/{\allowbreak}mincutefficacy}).
The efficacy is the violation of a cut divided by the Euclidian norm of its coefficient vector.

\subsection{Nonlinear Handler for Bilinear Expressions}
\label{sect:nlhdlrbilin}


The bilinear nonlinear handler identifies expressions of the form $y_1y_2$, where $y_1$ and $y_2$ are either non-binary variables of~\eqref{eq:minlp_extlp} or other expressions.
For a product $y_1y_2$, the expressions handler for products already provides linear under- and overestimators and domain propagation that is best possible when considering the bounds $[\low{y}_1,\upp{y}_1] \times [\low{y}_2,\upp{y}_2]$ only.
The nonlinear handler, however, can exploit linear inequalities over $y_1$ and $y_2$ to provide possibly tighter linear estimates and variable bounds.
These inequalities are found by projection of the LP relaxation onto variables $(y_1,y_2)$.
For more details, see M\"{u}ller et al.~\cite{Muller2020}.

\subsection{Nonlinear Handler for Convex and Concave Expressions}
\label{sect:nlhdlrconvex}


Two nonlinear handlers are available that try to detect convexity or concavity of a given expression $h_i^\text{lp}(x)$ and provide appropriate linear under- and overestimators.
The naming of the nonlinear handlers may be slightly confusing as the convex nonlinear handler checks for concavity of $h_i^\text{lp}(x)$ if overestimators are desired and the concave nonlinear handler checks for convexity if underestimators are desired.
After all, the detection algorithms of both nonlinear handlers are similar, so that they are discussed together here.
The linear estimators are computed differently, though.

In the following only the underestimating case ($\lesseqgtr_i$ being either~$\leq$ or~$=$) is considered.
The overestimating case is handled analogously.
The nonlinear handlers do not contribute to domain propagation so far.

\subsubsection{Detection}
\label{sect:nlhdlrconvexdetect}

Assume the constraint handler requests that underestimators of $h_i^\text{lp}(x)$ need to be found.
The convex nonlinear handler then seeks to find subexpressions of $h_i^\text{lp}(x)$ that need to be replaced by auxiliary variables $w_{i+1}^\text{lp},\ldots$ such that the remaining expression $h_i^\text{lp}(x,w_{i+1}^\text{lp},\ldots)$ is convex.
Similarly, the concave nonlinear handler seeks for $h_i^\text{lp}(x,w_{i+1}^\text{lp},\ldots)$ to be concave.
In both cases, the detection algorithm can aim for the remaining expression to be as large as possible.
This point will be revisited later.

To construct a maximal convex subexpression of $h_i^\text{lp}(x)$, the usual convexity and concavity detection rules are inverted and applied on $h_i^\text{lp}(x)$ in reverse order.
To do so, the expression is traversed in depth-first-search order, starting from the root of $h_i^\text{lp}(x)$.
With each subexpression the requirement of it being convex and/or concave is associated.
For the root, this will be convexity.
When a subexpression is considered, it is checked whether the subexpression can have the required curvature.
This is done by formulating requirements on convexity/concavity on the children of the subexpression.
If there are no conditions under which a subexpression can have the required curvature, then it is marked as to be replaced by an auxiliary variable.

For an example, consider the function $-\sqrt{\exp(x)}\sqrt{y} + \exp(x)$ with $\low{y}=0$.
First, it will be checked under which conditions on its arguments the sum will be convex.
This will create the requirements ``$\sqrt{\exp(x)}\sqrt{y}$ must be concave'' and ``$\exp(x)$ needs to be convex''.
Checking the former, the special structure $\sqrt{\cdot}\sqrt{\cdot}$ (a product of two power expressions, both exponents being $0.5$) may be detected and the requirements ``$\exp(x)$ must be concave'' and ``$y$ must be concave'' are created.
The check for ``$\exp(x)$ must be concave'' fails, i.e., there are no conditions on $x$ (other than $\low{x}=\upp{x}$) such that $\exp(x)$ is concave.
Therefore, this appearance of $\exp(x)$ is marked for replacement by an auxiliary variable.
The check for ``$y$ must be concave'' succeeds, since the function $y \mapsto y$ is both convex and concave.
The remaining check for ``$\exp(x)$ needs to be convex'' succeeds under the new condition ``$x$ needs to be convex'', which is satisfied.
Thus, the resulting maximal convex subexpression is $-\sqrt{w}\sqrt{y} + \exp(x)$, where $w$ is a new auxiliary variable and $w\leq\exp(x)$ is added to the extended formulation.
As this example has shown, it is possible that several appearances of the same subexpression ($\exp(x)$) are treated differently, depending on what requirements are imposed on the subexpression by its parents.

Four checks whether a subexpression can have a required curvature are currently implemented.
These are called in the given order.

\paragraph{Product Composition {\normalfont(\texttt{nlhdlr/\{convex,concave\}/cvxprodcomp})}}
Check whether a gi\-ven expression is a product of the form $af(b g(x)+c) g(x)$ with constants $a,b,c$ and repeating subexpression $g(x)$.
Considering the second derivative and by using available information on bounds and the \texttt{CURVATURE} and \texttt{MONOTONICITY} callbacks of the expression handlers, a condition on the curvature of $g(\cdot)$ can be derived that is sufficient for $af(b g(x)+c) g(x)$ to be convex or concave.

\paragraph{Signomial {\normalfont(\texttt{nlhdlr/\{convex,concave\}/cvxsignomial})}}
Check whether a given expression is a signomial, i.e., a product of power expressions, $c\prod_{j=1}^k f_j^{p_j}(x)$ with $c,p_j\in\R$ and subexpressions $f_j(x)$, $j=1,\ldots,k$.
If $\low{f}_j\geq 0$, then the product is convex if i) $p_j<0$ and $f_j(x)$ is concave for all $j=1,\ldots,k$ or ii) $\sum_{j=1}^kp_j\geq 1$ and there exists a $j^*\in\{1,\ldots,k\}$ such that $p_j<0$ and $f_j(x)$ is concave for all $j\neq j^*$ and $f_{j^*}(x)$ is convex.
Further, the product is concave if $p_j>0$ and $f_j(x)$ is concave for all $j=1,\ldots,k$ and $\sum_{j=1}^k p_j \leq 1$.
These conditions are proven by Maranas and Floudas~\cite{MaranasFloudas1995} and Chen and Huang~\cite{ChenHuang09} for the case that each $f_j(x)$ is equal to a variable.
If $\upp{f}_j<0$, then one can adapt by replacing $f_j$ by $-f_j$.
If the exponents satisfy the given conditions for convexity or concavity of the product, conditions on convexity and concavity of the subexpressions $f_j(x)$ are formulated.

\paragraph{Quadratics {\normalfont(\texttt{nlhdlr/\{convex,concave\}/cvxquadratic})}}
Check whether a given expression is quadratic, that is, is of the form $q(y)$ given by \eqref{eq:quadexpr}, where each $y_j$ is either an original variable or a subexpression.
With the same methods as in the nonlinear handler for quadratics, the sign of the eigenvalues of the quadratic coefficients matrix of $q(y)$ can be checked to decide whether $q(y)$ is convex or concave.
If $q(y)$ has the desired curvature, then it is required that every $y_j$ is linear.

The check on quadratics is currently disabled for the concave nonlinear handler.
It is not clear yet under which conditions it is beneficial to compute underestimators on a multivariate concave $q(y)$ via the methods of the concave nonlinear handler instead of handling each square and bilinear term of $q(y)$ separately.

\paragraph{Expression Handler}
For an expression $f(g_1(x),\ldots,g_k(x))$, call the \texttt{CURVATURE} callback of the expression handler for $f(\cdot)$.
If implemented and successful, then it provides convexity or concavity requirements for each $g_j(x)$.

\bigskip
As has been pointed out by Tawarmalani and Sahinidis~\cite{TawarmalaniSahinidis2005}, a tighter linear relaxation of a convex set (in the sense that less cuts are required to achieve the same outer-approximation) can usually be obtained when an extended formulation is used for function composition.
For instance, for $f(g(x))$ with both $f(\cdot)$ and $g(\cdot)$ being convex, $f(\cdot)$ being monotonically increasing, and $g(\cdot)$ being nonlinear, it is beneficial to consider the extended formulation $f(w)$, $w\geq g(x)$.
This is easily achieved in the detection algorithm by changing the requirement on a subexpression from convex or concave to linear (parameter \texttt{nlhdlr/convex/extendedform}).
Furthermore, the nonlinear handlers ignore expressions $h_i^\text{lp}(\cdot)$ that are a sum with more than one non-constant term (parameters \texttt{nlhdlr/\{convex,concave\}/detectsum}), unless the sum is a quadratic expression with at least one bilinear term, for example, $x^2 + 2xy + y^2$.

For the concave nonlinear handler, however, the observation by Tawarmalani and Sahinidis~\cite{TawarmalaniSahinidis2005} does not apply.
Instead, the number of variables in the expression for which estimators need to be computed can be an issue.
Therefore, here auxiliary variables are requested for multivariate linear subexpressions.
That is, even though concavity of $\log(x+y+z)$ can be recognized, the extended formulation $\log(w)$, $w\lesseqgtr x+y+z$, is used.
This way, only one- instead of three-dimensional underestimators need to be calculated.

Finally, if $h_i^\text{lp}(x)$ were transformed by the nonlinear handler into $h_i^\text{lp}(x,w_{i+1},\ldots,w_{m^\text{lp}})$ such that the corresponding expression has only original or auxiliary variables as children, then the detection of the nonlinear handler is reported as failed (parameter \texttt{nlhdlr/{\allowbreak}\{convex,concave\}/{\allowbreak}handletrivial}).
Instead, the default nonlinear handler will provide linear estimates via the \texttt{ESTIMATE} callback of the expression handlers.
We assume that these are more efficient than the generic implementation in the convex and concave nonlinear handlers.

\subsubsection{Underestimators for Convex Expressions}

For the convex function $h_i^\text{lp}(x,w_{i+1}^\text{lp},\ldots,w_{m^\text{lp}}^\text{lp})$ and reference point $(\hat x,\hat w^\text{lp})$, an underestimator is given by computing a tangent on the graph of $h_i^\text{lp}(\cdot)$ at $(\hat x,\hat w^\text{lp})$:
\[
  h_i^\text{lp}(\hat x,\hat w_{i+1}^\text{lp},\ldots,\hat w_{m^\text{lp}}^\text{lp}) + 
  \nabla h_i^\text{lp}(\hat x,\hat w_{i+1}^\text{lp},\ldots,\hat w_{m^\text{lp}}^\text{lp})
  \left(\begin{array}{c}x-\hat x\\ w^\text{lp}-\hat w^\text{lp}\end{array}\right).
\]
If, however, $h_i^\text{lp}(x,w_{i+1}^\text{lp},\ldots,w_{m^\text{lp}}^\text{lp})$ is univariate, that is, $h_i^\text{lp}(x,w_{i+1}^\text{lp},\ldots,w_{m^\text{lp}}^\text{lp})=f(y)$ for some variable $y$, and $y$ is integral, then taking the secant on the graph of $f(y)$ can give a tighter underestimator:
\[
 f(\lfloor \hat y\rfloor) + (f(\lfloor \hat y\rfloor+1) - f(\lfloor \hat y\rfloor)) \,(y-\lfloor \hat y\rfloor).
\]

\subsubsection{Underestimators for Concave Expressions}

To simplify notation, let $f(y)$, $y\in\R^k$, be the concave function $h_i^\text{lp}(x,w_{i+1}^\text{lp},\ldots,w_{m^\text{lp}}^\text{lp})$ for which a linear underestimator needs to be computed.
Assume further that variables that are currently fixed have been replaced by the corresponding constant.

Since $f(y)$ is concave, its convex envelope with respect to $[\low{y},\upp{y}]$ is vertex-polyhedral, that is, it is a polyhedral function which vertices correspond to the vertices of $[\low{y},\upp{y}]$.
Therefore, any hyperplane $\alpha y+\beta$ that underestimates $f(y)$ in all vertices of $[\low{y},\upp{y}]$ is a valid linear underestimator.
Maximizing $\alpha \hat y + \beta$ under these constraints gives an underestimator that is as tight as possible in the reference point $\hat y$.
For the frequent cases $k=1$ and $k=2$, routines that directly compute such an underestimator are available.
For $k>2$, a linear program is solved.

Since such underestimators need to be computed repeatedly for varying domains, updates to the cut generating linear program (CGLP) are kept to a minimum.
By use of the linear transformation $T:[0,1]^k \to [\low{y},\upp{y}]$ given by $T(\tilde y)_i = \low{y}_i + (\upp{y}_i - \low{y}_i)\tilde y_i$, $i=1,\ldots,k$, the matrix of the CGLP can be kept constant:
\begin{align*}
  \max_{\tilde\alpha,\tilde\beta}\; & \tilde\alpha^\T T^{-1}(\hat y) + \tilde\beta, \\
  \text{s.t.}\; & \tilde\alpha^\T \tilde y + \tilde\beta \leq f(T(\tilde y)), \qquad \forall \tilde y \in \{0,1\}^k.
\end{align*}
The cut in the original $y$-space is then recovered via $\alpha_j = \frac{\tilde\alpha_j}{\upp{y}_j-\low{y}_j}$, $j=1,\ldots,k$, and $\beta = \tilde\beta - \sum_{j=1}^k \alpha_j\low{y}_j$.
Since the CGLP typically has more rows than columns, the dual of CGLP is formulated and solved.
To increase the chance that $\alpha y+\beta$ is a facet of the convex envelope of $f(y)$, the reference point is perturbed and moved into the interior of $[\low{y},\upp{y}]$.

At the moment, underestimators for concave functions in more than 14 variables are not computed due to the size of the CGLP being exponential in $k$.
In fact, the detection algorithm in the concave nonlinear handler already returns unsuccessful if the recognized concave expression has more than 14 variables.
Dynamic row or column generation methods could be added to overcome this limit~\cite{BaoSahinidisTawarmalani2009}.

Since the underestimator may not be tight at $(\hat y,f(\hat y))$, all variables are registered as branching candidates by this nonlinear handler.

\subsection{Nonlinear Handler for Quotients}
\label{sect:nlhdlrquotient}


Note that the available expression handlers (see Section~\ref{sect:expr}) do not include a handler for quotients since they can equivalently be written using a product and a power expression.
However, the default extended formulation for an expression $y_1y_2^{-1}$ is given by replacing $y_2^{-1}$ by a new auxiliary variable $w$.
The linear outer-approximation is then obtained by estimating $y_1w$ and $y_2^{-1}$ separately.
The quotient nonlinear handler can provide tighter estimates by checking whether a given function $h_i^\text{lp}(x)$ can be cast as
\begin{equation}\label{quotient_constraint}
 f(y) = \frac{ay_1 + b}{cy_2 + d} + e
\end{equation}
with $a,b,c,d,e\in\R$, $a,c\neq 0$, and $y_1$ and $y_2$ being either original variables ($x$) or subexpressions of $h_i^\text{lp}(x)$.
At the moment, only expressions $h_i^\text{lp}(\cdot)$ that are of the form
$ay_1y_2^{-1}$ 
or $ay_1y_2^{-1} + by_2^{-1} + e$ 
are recognized as quotient and it is checked whether $y_j$ equals $a_j'y_j'+b_j'$ for some $a_j',b_j'$, $j=1,2$.
For estimation and domain propagation, the univariate ($y_1=y_2$) and bivariate ($y_1 \neq y_2$) cases are handled separately.

\subsubsection{Univariate Quotients ($y_1=y_2$)}

If $-d/c\not\in [\low{y}_2,\upp{y}_2]$, then $f(y)$ is either convex or concave on $[\low{y},\upp{y}]$.
Thus, under- and overestimators are computed via a tangent or a secant on the graph of $f(y)$.
%
If the singularity is in the domain of $y$, then no estimator can be computed.

For forward domain propagation, observe that the minimum and maximum of $f(y)$ is attained at $\low{y}$ or $\upp{y}$ if $-d/c\not\in [\low{y}_2,\upp{y}_2]$.
It is therefore sufficient for evaluate $f(y)$ at $\low{y}$ and $\upp{y}$ to obtain $f([\low{y},\upp{y}])$.
If the singularity is in the domain of $y$, then no finite bounds on $f(y)$ can be computed.

For backward domain propagation, let $[\low{f},\upp{f}]$ be the bounds given for $f(y)$.
Inverting~\eqref{quotient_constraint} yields
\[
  y = \frac{b-d [\low{f},\upp{f}]}{c [\low{f},\upp{f}]-a}.
\]
This interval can be evaluated as in forward propagation.

\subsubsection{Bivariate Quotients}

Let $\low{x}, \upp{x}, \low{y}, \upp{y}$ denote lower and upper bounds on variables $x$ and $y$,
respectively.
In the bivariate case, estimators are computed for $y_1'/y_2'$ and then transformed back into $(y_1,y_2)$-space.
Thus, for the following, $a=1$, $b=0$, $c=1$, $d=0$ is assumed.

If $0\in[\low{y}_2,\upp{y}_2]$, then no estimators can be computed.


If $\low{y}_1 \geq 0$ and $\low{y}_2 > 0$, then an underestimator is obtained by computing a tangent on the graph of the convex underestimator of $f(y)$, which is given by Zamora and Grossmann~\cite{ZaGr98b} as
\[
  \frac{1}{y_2} \left(\frac{y_1 + \sqrt{\low{y}_1\upp{y}_1}}{\sqrt{\low{y}_1} + \sqrt{\low{y}_1}}\right)^2.
\]
An overestimator is given by a hyperplane that passes either through the points
$(\low{y}_1, \low{y}_2$, $\low{y}_1/\low{y}_2)$,
$(\low{y}_1, \upp{y}_2, \low{y}_1/\upp{y}_2)$ and
$(\upp{y}_1, \low{y}_2, \upp{y}_1/\low{y}_2)$ or through the points
$(\low{y}_1, \upp{y}_2, \low{y}_1/\upp{y}_2)$,
$(\upp{y}_1, \low{y}_2$, $\upp{y}_1/\low{y}_2)$ and
$(\upp{x}_1,\upp{y}_2,\upp{y}_1/\upp{y}_2)$.
The choice of the three points depends on which combination yields the estimator that is tightest
at the given reference point, for example, the LP solution to be separated.

If $0\in[\low{y}_1,\upp{y}_1]$ and $\low{y}_2 > 0$, recall that the nonlinear handler works on the constraint $y_1/y_2 \lesseqgtr_i w_i^\text{lp}$ in~\eqref{eq:minlp_extlp} and that it is valid to replace $\lesseqgtr_i$ by $=$.
Rewriting as $y_1 = y_2w_i^\text{lp}$, linear McCormick envelopes~\cite{McCormick1976} are computed for $y_2w_i^\text{lp}$.
These inequalities are then rearranged to obtain linear estimators on $w_i^\text{lp}$.

The cases of $\upp{y}_1<0$ and/or $\upp{y}_2<0$ are handled similarly, up to swapping and negation of
bounds, coefficients, and inequality signs.

Since each variable appears only once in the expression, the default domain propagation for $f(y)$ does not suffer from the dependency problem of interval arithmetics.
Therefore, the nonlinear handler does not provide an extra implementation of domain propagation for the bivariate case.

\subsection{Nonlinear Handler for Perspective Reformulation}
\label{sect:nlhdlrpersp}


This nonlinear handler creates strengthened cutting planes for constraints that depend
on semi-continuous variables.
A variable $x_j$, $j\in\varindex$, is semi-continuous with respect to the binary indicator variable $x_{j'}$, $j'\in\intvarindex$, if it is restricted to the domain $[\low{x}^1_j, \upp{x}^1_j]$ when $x_{j'}=1$ and has a fixed value $x^0_j$ when $x_{j'}=0$.
In the rest of this subsection, 
the superscript $0$ denotes the value of a semi-continuous variable at $x_{j'}=0$.

Consider the constraint
\begin{equation}\label{semicontinuous_constraint}
 h_i^\text{lp}(x,w^\text{lp}_{i+1},\ldots) \lesseqgtr w_i^\text{lp}
\end{equation}
and write $h^\text{lp}_i(\cdot)$ as a sum of its nonlinear and linear parts:
\[
  h^\text{lp}_i(x,w^\text{lp}_{i+1},\ldots) = h_i^\text{nl}(x_\text{nl},w^\text{lp}_\text{nl}) + h_i^\text{l}(x_\text{l},w^\text{lp}_\text{l}),
\]
where $h_i^\text{nl}(\cdot)$ is a nonlinear function, $h_i^\text{l}(\cdot)$ is a linear function, $x_\text{nl}$ and $w^\text{lp}_\text{nl}$ are the vectors of variables $x$ and $w^\text{lp}$, respectively, that appear only in the nonlinear part of $h_i^\text{lp}$, and $x_\text{l}$ and $w^\text{lp}_\text{l}$ are the vectors of variables $x$ and $w^\text{lp}$, respectively, that appear only in the linear part of $h_i^\text{lp}(\cdot)$.

The perspective handler works on Constraint~\eqref{semicontinuous_constraint} if $x_\text{nl}$ and $w^\text{lp}_\text{nl}$ are semi-continuous with respect to the same indicator variable $x_{j'}$, and at least one other nonlinear handler provides estimation (\texttt{ESTIMATE} callback) for $h_i^\text{lp}(\cdot)$.
Thus, a nonlinear handler that implements only the \texttt{ENFO} callback, such as, for example, the SOC handler, is not suitable.

\subsubsection{Detection of Semi-continuous Variables}

To determine whether a variable $x_j$ is semi-continuous, the detection callback
searches for pairs of implied bounds on $x_j$ with the same indicator $x_{j'}$:
\begin{align*}
x_j &\leq \alpha^{(u)} x_{j'} + \beta^{(u)},\\
x_j &\geq \alpha^{(\ell)} x_{j'} + \beta^{(\ell)}.
\end{align*}
If $\beta^{(u)} = \beta^{(\ell)}$, then $x_j$ is a semi-continuous variable and $x_j^0 = \beta^{(u)}$,
$\low{x}^1_j = \alpha^{(\ell)} + \beta^{(\ell)}$, and $\upp{x}^1_j = \alpha^{(u)} + \beta^{(u)}$.

This information can be obtained either from linear constraints in $x_j$ and $x_{j'}$ or by
finding implicit relations between $x_j$ and $x_{j'}$.
Such relations can be detected by probing, which fixes $x_{j'}$ to its possible values and propagates
all constraints in the problem, thus detecting implications of $x_{j'} = 0$ and $x_{j'} = 1$.
SCIP stores the implied bounds in a globally available data structure.

The perspective nonlinear handler also detects semi-continuous auxiliary
variables. 
Given $h^\text{lp}_i(x, w_{i+1}^\text{lp},\ldots) \lesseqgtr_i w^\text{lp}_i$, where $x,w_{i+1}^\text{lp},\ldots$ are semi-continuous variables
depending on the same indicator $x_{j'}$, the auxiliary variable $w^\text{lp}_i$ can also be assumed to be semi-continuous since it is valid to replace $\lesseqgtr_i$ by $=$.

\subsubsection{Separation}

Suppose that the current relaxation solution violates Constraint~\eqref{semicontinuous_constraint}.
If the non-perspective nonlinear handlers claimed that estimators of $h^\text{lp}_i(\cdot)$ depend on variable bounds, probably because the functions is nonconvex, then probing is first
performed for $x_{j'}=1$ in order to tighten the implied bounds on variables $x, w_{i+1}^\text{lp}, \ldots$.
Linear underestimators (for ``$\leq$'' constraints) or overestimators
(for ``$\geq$'' constraints) that are valid when $x_{j'}=1$ are then obtained for the tightened bounds.
This estimator $\ell(\cdot)$ can be separated into parts corresponding to the nonlinear and linear variables of $h_i^\text{lp}(\cdot)$, respectively:
\[
  \ell(x,w_{i+1}^\text{lp},\ldots) =
  \ell^\text{nl}(x_\text{nl},w_\text{nl}^\text{lp}) + \ell^\text{l}(x_\text{l},w_\text{l}^\text{lp})
\]

%
An extension procedure is applied to the nonlinear part to ensure it is valid and tight for $x_{j'}=0$, while the linear part can remain unchanged since it shares none of the variables with the nonlinear part:
\[
  \ell^\text{nl}(x_\text{nl},w_\text{nl}^\text{lp}) +
  \left(h^\text{nl}_i(x^0_\text{nl},w^{\text{lp},0}_\text{nl}) - \ell^\text{nl}(x^0_\text{nl},w^{\text{lp},0}_\text{nl})\right)(1-x_{j'}) +
  \ell^\text{l}(x_\text{l},w_\text{l}^\text{lp}).
\]
This extension ensures that the estimator is equal to $h_i^\text{lp}(x,w^\text{lp}_{i+1},\ldots)$ for
$x_{j'}=0$, $x_\text{nl} = x^0_\text{nl}$, and $w_\text{nl} = w^0_\text{nl}$, and equal to $\ell(x,w_{i+1}^\text{lp},\ldots)$ for $x_{j'}=1$.
In the convex case, cuts thus obtained are equivalent to the classic perspective cuts~\cite{frangioni2006perspective}.
More details on the implementation in SCIP can be found in the paper by Bestuzheva et al.~\cite{BestuzhevaGleixnerVigerske2021}.

\subsection{Separator for Cuts from the Reformulation-Linearization Technique}
\label{sect:rlt}


A nonlinearity that appears frequently is a product between two variables and/or functions.
The separator for Reformulation-Linearization Technique (RLT) cuts~\cite{adams1986tight,adams1990linearization,adams1993mixed} for bilinear product relations in~\eqref{eq:minlp_extlp} and the separators discussed in the following two sections focus on enforcing the relationship between a product of two variables (original or auxiliary) and a corresponding auxiliary variable.
The RLT separator can additionally reveal linearized products between binary and continuous variables.

There exist variations of the RLT that can be applied to any (not necessarily quadratic) polynomials~\cite{sherali1992global}.
This separator, however, deals with bilinear products only.

In the following, $x$ refers to any variable of~\eqref{eq:minlp_extlp} and $X_{i,j}$ refers to the auxiliary variable ($w^\text{lp}$) that is associated with a constraint $x_ix_j \lesseqgtr X_{i,j}$ in~\eqref{eq:minlp_extlp}.
Note that $X_{i,j}$ may not exist in~\eqref{eq:minlp_extlp} for every pair of $x_i$ and $x_j$, even when $x_ix_j$ appears in some constraint of~\eqref{eq:minlp} (for example, auxiliary variables are not created for terms in convex quadratic constraints).
Both $X_{i,j}$ and $X_{j,i}$ refer to the same variable.

Given a product relation $X_{ij} = x_ix_j$, where $x_i \in [\low{x}_i,\upp{x}_i]$,
$x_j \in [\low{x}_j,\upp{x}_j]$ and a linear constraint $a^\T x \leq b$, RLT cuts
are derived by first multiplying the constraint by nonnegative bound factors $(x_i - \low{x}_i)$,
$(\upp{x}_i - x_i)$, $(x_j - \low{x}_j)$, and $(\upp{x}_j - x_j)$.
For instance, consider multiplication by the factor $(x_i - \low{x}_i)$, which yields a valid nonlinear inequality:
\begin{equation}\label{rlt_reformulated}
 a^\T x\, (x_i - \low{x}_i) \leq b\,(x_i - \low{x}_i).
\end{equation}
This is referred to as the reformulation step.

The linearization step is then performed for all terms $x_kx_i$ in~\eqref{rlt_reformulated}.
If a product relation $X_{ki} = x_kx_i$ exists, then the product is replaced with $X_{ki}$.
If $x_k$ and $x_i$ are contained in the same clique, the product is replaced with an equivalent
linear expression.
Otherwise, it is replaced by a linear under- or overestimator such that the inequality remains
valid.
By default, RLT cuts are constructed only for combinations of rows and bound factors where
all relations $X_{ki} = x_kx_i$ exist (parameter \texttt{separating/rlt/maxunknownterms}).

\subsubsection{Implicit Product Detection}

Bilinear product relations in which one of the multipliers is binary can equivalently be written via
mixed-integer linear constraints.
Likewise, MILP constraints representing such relations can be identified in order to derive
these implicit bilinear products.

Consider two linear constraints depending on the same three variables $x_i$, $x_j$ and $x_k$,
where $x_i$ is binary:
\begin{subequations}\label{rlt_detection_linrels}
\begin{gather}
a_1x_i + b_1x_k + c_1x_j \leq d_1,\label{rlt_detection_linrels1}\\
a_2x_i + b_2x_k + c_2x_j \leq d_2.\label{rlt_detection_linrels2}
\end{gather}
\end{subequations}
If $b_1b_2 > 0$ and $c_2b_1 - b_2c_1 \neq 0$, then these constraints imply a product relation,
\begin{equation}\label{rlt_derived_product}
  Ax_i + Bx_k + Cx_j + D \lesseqgtr x_ix_j,
\end{equation}
where the coefficients $A$, $B$, $C$, and $D$ and the inequality sign are obtained by:
\begin{itemize}
\item setting $x_i$ to $1$ in \eqref{rlt_detection_linrels1} and \eqref{rlt_derived_product}, and
requiring that the coefficients are similar for each variable, and the constants are equal;
\item setting $x_i$ to $0$ in \eqref{rlt_detection_linrels2} and \eqref{rlt_derived_product}, and
similarly requiring equivalence;
\item solving the linear system resulting from the first two steps.
\end{itemize}

SCIP analyses the linear constraints in the problem and stores all detected implicit products.
RLT cuts that use these products may strengthen the default continuous relaxation
$\{ (x_i,x_j,x_k) : x_i \in [0,1], \eqref{rlt_detection_linrels}\}$.

\subsubsection{Separation}

Let $(\hat x,\hat X)$ be the solution to be separated.
In order to reduce the computational cost of RLT cut separation, SCIP takes into account the signs
of coefficients of linear constraints and signs of product relation violations.
In particular, when multiplying a constraint $a^\T x \leq b$ by a bound factor, the resulting RLT cut
can only be violated if $a_k\hat x_k\hat x_j < a_k\hat X_{kj}$, that is, when replacing the product with the
corresponding variable increases the violation of the inequality.
This fact is used to ignore combinations of linear constraints and bound factors that will not
produce a violated cut, thus reducing the computational effort.

This is implemented via a row marking algorithm which, for every variable $x_i$ that participates in
bilinear products, iterates over all variables $x_j$ that appear in products together with $x_i$.
When it encounters a violated product, the algorithm iterates over all linear rows where $x_j$ has
a nonzero coefficient and stores them in a sparse sorted array together with the marks indicating
which bound factors of $x_i$ it should be multiplied with.
The cut generation algorithm then iterates over the array of marked rows and constructs RLT cuts
from the products of each row with the suitable bound factors.

More details on the algorithms and implementation will be included in the upcoming
paper~\cite{achterbergEfficient}.

\subsection{Separator for Principal Minors of $X\succeq xx^\T$}
\label{sect:minor}


Another new separator that enforces bilinear product relations in~\eqref{eq:minlp_extlp} is \texttt{sepa\_{\allowbreak}minor}.
The notation introduced in the previous section is used.

A convex relaxation of condition $X = xx^\T$ is given by requiring $X-xx^\T$ to be positive definite.
Separation for the set $\{(x,X) : X-xx^\T \succeq 0\}$ itself is possible, but cuts are typically dense and may include variables $X_{ij}$ for products that do not exist in the problem~\cite{QualizzaBelottiMargot2012}.
Therefore, \texttt{sepa\_minor} considers only (principle) $2 \times 2$ minors of $X-xx^\T$, which also need to be positive semi-definite.
By Schurs complement, this means that the condition
\begin{equation}
   \label{eq:minorposdef}
	A_{ij}(x,X) \defi \begin{bmatrix} 1 & x_i & x_j \\ x_i & X_{ii} & X_{ij} \\ x_j & 
	X_{ij} & X_{jj} \end{bmatrix} \succeq 0
\end{equation}
needs to hold.
The separator detects principle minors for which $X_{ii}$, $X_{jj}$, $X_{ij}$ exist and enforces $A_{ij}(x,X)\succeq 0$.

To identify which entries of the matrix $X$ exist, the separator iterates over 
the available nonlinear constraints. For each constraint, its expressions are 
explored and all expressions of the form $x_i^2$ and $x_ix_j$ are collected. 
Then, the separator iterates through the found bilinear terms $x_ix_j$ and if 
the corresponding expressions $x_i^2$ and $x_j^2$ exist, a minor is detected.

Let $(\hat{x},\hat{X})$ be a solution that violates~\eqref{eq:minorposdef}, i.e., there exists an eigenvector $v\in\R^3$ of $A_{ij}(\hat{x},\hat{X})$ with $v^\T A_{ij}(\hat{x},\hat{X})v<0$.
To separate $(\hat{x},\hat{X})$, \texttt{sepa\_minor} adds the globally valid linear inequality $v^\T A_{ij}(x,X)v \geq 0$ to the separation storage of SCIP.

For circle packing instances, the minor cuts are not really 
helpful~\cite{Khajavirad2017}. Since experiments showed that SCIP's overall 
performance was negatively affected, circle packing constraints are identified 
and their bilinear terms are ignored by \texttt{sepa\_minor} (parameter \texttt{separating/minor/ignorepackingconss}).

\subsection{Separator for Intersection Cuts on Rank-1 Constraint for $X$}
\label{sect:interminor}


Another new separator that enforces bilinear product relations in~\eqref{eq:minlp_extlp} is \texttt{sepa\_{\allowbreak}interminor}.
The notation introduced in Section~\ref{sect:rlt} is used.

Since $X = xx^\T$ has rank 1 in any feasible solution, any $2\times 2$ minor $\begin{pmatrix}
X_{i_1j_1} & X_{i_1j_2} \\
X_{i_2j_1} & X_{i_2j_2}
\end{pmatrix}$
of $X$ needs to have determinant 0.
That is, for any set of variable indices $i_1$, $i_2$, $j_1$, $j_2$ with $i_1\neq i_2$ and $j_1\neq j_2$, the condition
\begin{align}
  \label{eq:impliedquad}
  X_{i_1j_1}X_{i_2j_2} = X_{i_1j_2}X_{i_2j_1}
\end{align}
needs to hold.
If all variables in this constraint exist in the problem and the solution $(\hat x,\hat X)$ that is to be separated violates~\eqref{eq:impliedquad}, the separation strategy described in Section~\ref{sect:intercuts} is used to add (strengthened) intersection cuts that separate $(\hat x,\hat{X})$.
Additionally, it is also possible (parameter \texttt{separating/interminor/usebounds}) to use the bounds on $x_{i_1}$, $x_{i_2}$, $x_{j_1}$, $x_{j_2}$ to improve the cut by enlarging the corresponding $S$-free set~\cite{ChmielaMunozSerrano2021}.

The separator is currently disabled by default.

\subsection{Revised Primal Heuristic that Solves NLP Subproblem}
\label{sect:subnlp}

The primal heuristic \texttt{subnlp} targets problems like~\eqref{eq:minlp}, but runs on any CIP where the NLP relaxation is enabled.
Given a point $\tilde x$ that satisfies the integrality requirements ($\tilde x_i\in\Z$ for all $i\in\intvarindex$), the heuristic fixes all integer variables to the values given by $\tilde x$ in a copy of the CIP, presolves this copy, and triggers a solution of the NLP relaxation by an NLP solver using $\tilde x$ as starting point.
If the NLP solver, such as \ipopt, finds a solution that is feasible (and often also locally optimal) for the NLP relaxation, it is tried whether it is also feasible for the CIP.
If the CIP is a MINLP, then this should usually be the case.
The starting point $\tilde x$ can be the current solution of the LP relaxation if integer-feasible, can be a point that a primal heuristic that searches for feasible solutions of the MILP relaxation has computed, or can have been passed on by other primal heuristics that look for MINLP solutions, such as \texttt{undercover} or \texttt{mpec}.

The \texttt{subnlp} primal heuristic, which is implemented in virtually any global MINLP solver, had been added to SCIP together with the support for quadratic constraints (SCIP~1.2.0).
The rewrite of the algebraic expression system (Section~\ref{sect:expr}) and the handling of nonlinear constraints (Section~\ref{sect:consnl}) and the updates to the NLP solver interfaces and NLP relaxation (Section~\ref{sect:nlp}) were a good opportunity for a thorough revision of the heuristic.

\paragraph{Starting Condition and Iteration Limit}
By default, the heuristic is called in every node of the branch-and-bound tree, but invoking an NLP solver whenever a starting point $\tilde x$ is available would be too costly.
After the heuristic has been run, it therefore waits until a certain number of nodes have been processed.
How many nodes these are depends on the success of the heuristic in previous calls, the number of iterations the NLP solver used in previous calls, and the iteration limit that would be imposed for the following NLP solve.
Previously, the iteration limit was essentially static, which could mean that on problems with difficult NLPs a lot of effort was wasted on NLP solves that were interrupted by a too small iteration limit.

With SCIP 8, the heuristic tries to adapt the iteration limit to the NLPs to be solved.
For that, the heuristic counts how often an NLP solve stopped due to an iteration limit ($n^\text{iterlim}$) and how often it finished successfully, that is, stopped because convergence criteria were fulfilled ($n^\text{okay}$).
Let $i^\text{iterlim}$ be the highest iteration limit used among all NLP solves that stopped due to an iteration limit and let $i^\text{okay}$ be the total number of iterations used in all NLP solvers that finished successfully.
Further, let $i^\text{min}$ be a minimal number of iterations that should be granted to every NLP solve (parameter \texttt{heuristics/subnlp/itermin} = 20).
Finally, let $n^\text{init}$ be the number of initial NLP solves that should be granted $i^\text{init}$ many iterations (parameters \texttt{heuristics/subnlp/ninitsolves} = 2 and \texttt{heuristics/subnlp/iterinit} = 300).
The iteration limit $i^\text{next}$ for the next NLP solve is then decided as follows:
\begin{enumerate}
 \item If $n^\text{iterlim} > n^\text{okay}$, then $i^\text{next} \defi \max(i^\text{min}, 2i^\text{iterlim})$.
   That is, double the iteration limit if more solves ran into an iteration limit than were successful.
 \item Otherwise, if $n^\text{okay} > n^\text{init}$, then $i^\text{next} \defi \max(i^\text{min}, 2\frac{i^\text{okay}}{n^\text{okay}})$.
   That is, if there were a few successful solves to far, then use twice the average number of iterations spend in these solves as iteration limit.
 \item Otherwise, if $n^\text{okay}$, then $i^\text{next} \defi \max(i^\text{min}, i^\text{init}, 2\frac{i^\text{okay}}{n^\text{okay}})$.
   That is, consider also $i^\text{init}$ if there had not been enough successful solves so far.
 \item Otherwise, $i^\text{next} \defi \max(i^\text{min}, i^\text{init})$.
\end{enumerate}

To decide whether to execute the heuristic, an iteration contingent $i^\text{cont}$ is calculated and checked against $i^\text{next}$.
Compared to SCIP 7, this has received only minor updates:
\begin{enumerate}
 \item Initialize $i^\text{cont} \defi 0.3(\text{number of nodes processed} + 1600)$ (parameters
   \texttt{heuristics/{\allowbreak}subnlp/{\allowbreak}\{nodesfactor,nodesoffset\}}).
 \item Weigh by previous success of heuristic: Let $n^\text{tot}$ the total number of times the heuristic has run and $n^\text{sol}$ the number of solutions found by the heuristic.
   If the heuristic ran a few times and is no longer in a phase where it tries to find a suitable iteration limit, then weigh $i^\text{cont}$ by success of heuristics.
   That is, if $n^\text{tot} - n^\text{iter} > n^\text{init}$, then $i^\text{cont} \defi \frac{n^\text{sol}+1}{n^\text{tot}+1}i^\text{cont}$.
   Parameter $\beta\defi$\texttt{heuristics/subnlp/successrateexp} allows to replace $\frac{n^\text{sol}+1}{n^\text{tot}+1}$ by $(\frac{n^\text{sol}+1}{n^\text{tot}+1})^\beta$.
 \item Let $i^\text{tot}$ be the total number of iterations used in all NLP solves (successful or not) so far.
   Then $i^\text{cont} \defi i^\text{cont} - i^\text{tot}$.
 \item If $i^\text{cont} \geq i^\text{next}$, then the heuristic is run with $i^\text{next}$ as iteration limit for the NLP solver.
\end{enumerate}

\paragraph{Presolve}
The heuristic triggers a solve of the NLP relaxation of SCIP in a copy of the CIP.
When the heuristic is run for a starting point $\tilde x$, integer variables are fixed to the values given in $\tilde x$, the current primal bound is set as cutoff, and SCIP's presolve is run with presolve emphasis set to ``fast''. 
The aim of the presolve is to propagate the fixing of the integer variables in the problem since many NLP solvers, in particular those that are interfaced by SCIP, only implement a very limited presolve.
After presolve, if the problem is not empty or infeasible, SCIP is put into a state where its NLP relaxation can be solved.
If the original CIP is a MINLP, then solutions that are feasible to this NLP relaxation should also be feasible in the original problem.
Further, also solutions that are found during presolve are passed on to the original problem.

This process of fixing integer variables, setting a cutoff, and presolving the CIP repeats every time the heuristic is run.
If, however, there are no binary or integer variables, then setting a cutoff and presolve is skipped and the copied problem is kept in a state where its NLP relaxation can be solved.

\paragraph{NLP Solve}
The NLP relaxation in the presolved copied CIP instance is solved by a NLP solver that is interfaced by SCIP.
The solver is given $\tilde x$ as starting point and the iteration limit is set to $i^\text{next}$.
If the NLP solver is \ipopt, then also the ``expect infeasible problem'' heuristic of \ipopt is enabled.

If the solver claims to have found a feasible solution, then it is tried to add this solution to the original problem.
This can fail for three reasons: the objective function value is not good enough, the NLP relaxation is missing some constraints of the original CIP, or the solution is only slightly infeasible due to presolve reductions.
For example, due to tolerances, bounds of aggregated variables might be slightly violated.
To work around this case, if a solution is not accepted, its objective value is not worse than the current primal bound, its maximal constraint violation is close to the feasibility tolerance, and the copied problem has been presolved ($\intvarindex\neq\emptyset$), then the NLP is resolved with a tightened feasibility tolerance (parameter \texttt{heuristics/subnlp/feastolfactor}).
For this resolve, warmstart from the previous solution is enabled and the iteration count of the previous NLP solve is used as iteration limit.
If the NLP resolve succeeds and produces a solution that is accepted in the original problem, then the tightened feasibility tolerance is used for all following NLP solves by the heuristic.

\subsection{NLP Relaxation and Interfaces to NLP Solvers and Automatic Differentiation}
\label{sect:nlp}


The updated expressions framework (Section~\ref{sect:expr}) triggered a revision of the NLP relaxation and the interfaces to NLP solvers (\texttt{NLPI}) and automatic differentiation (\texttt{EXPRINT}).

\subsubsection{NLP Relaxation}
\label{sect:nlprelax}

The rows of the NLP relaxation (\texttt{SCIP\_NLROW}) no longer distinguish a quadratic part.
Therefore, rows now have the form
\[
  \text{left-hand side} \leq \text{linear terms} + \text{nonlinear term} \leq \text{right-hand side}
\]
where the nonlinear term is given as a SCIP expression.
When the nonlinear constraint handler (see Section~\ref{sect:consnl}) creates an NLP row for a constraint $\low{g}\leq g(x)\leq \upp{g}$ of \eqref{eq:minlp}, it separates linear terms from $g(x)$.
The constraint handlers \texttt{and}, \texttt{bounddisjunction}, \texttt{knapsack}, \texttt{linear}, \texttt{linking}, \texttt{logicor}, \texttt{setppc}, and \texttt{varbound} now add themselves to the NLP relaxation.
Previously, \texttt{and}, \texttt{bounddisjunction}, and \texttt{linking} constraints were not added.
For \texttt{bounddisjunction}, only univariate constraints are added.

Further, it is pointed out that the NLP relaxation of SCIP is no longer based on the extended formulation~\eqref{eq:minlp_extlp}, but is now closer to the continuous relaxation of the original problem~\eqref{eq:minlp}.

\subsubsection{Interfaces to NLP Solvers}
\label{sect:nlpi}

Since expression handlers are now proper SCIP plugins that require a SCIP pointer for many operations and since expressions are used to specify NLPs, also the NLP solver interfaces (\texttt{NLPI}) are now proper SCIP plugins that require a SCIP pointer.
However, as before, the NLPs that are specified via an NLPI can be independent of the problem that is solved by SCIP.
For the expressions in the objective and constraints of such an NLP this means that the ``var'' expression handler, which refers to a SCIP variable (\texttt{SCIP\_VAR*}), cannot be used.
Instead, the handler for ``varidx'' expressions, which refer to a variable index, needs to be used.
As a consequence, the evaluation and differentiation methods of expressions, which work with a SCIP solution (\texttt{SCIP\_SOL}), are not available (the \texttt{EVAL} callback of the ``varidx'' expression handler raises an error).
Instead, the NLP solver interfaces either implement their own evaluation and differentiation or resort to the helper functions implemented in \texttt{nlpioracle.\{h,c\}}.

Next to the adjustments to the new expressions framework, further updates and removals to the NLPI callbacks were implemented.
For a detailed list, see the \texttt{CHANGELOG}.
A notable change, though, is that parameter settings that specify the working limits and tolerances of an NLP solve are now passed directly to the \texttt{NLPISOLVE} callback and, thus, are used for the corresponding solve only.
The same applies to the NLP relaxation of SCIP and \texttt{SCIPsolveNLP()} (now a macro).
The default values for the NLP solve parameters are now uniform among all NLP solvers and some parameters were added, removed, or renamed.
The solve statistics now include information on the violation of constraints and variable bounds of the solution, if available.

The problem and optimization statistics that SCIP collects and prints on request (\texttt{display statistics}) now include a table for each used NLP solver, which prints the number of times the solver was used, the time spend, and how often each termination and solution status occurred.
Additionally, the time spend for evaluation and differentiation can be shown (parameter \texttt{timing/nlpieval}).

As before, SCIP includes interfaces to the NLP solvers \filtersqp, \ipopt, and \worhp.
In particular the interface to \ipopt has been improved. Only some points are mentioned here:
\begin{itemize}
\item Warmstarts from a primal/dual solution pair, either set via \texttt{NLPISETINITIALGUESS} or by using the solution from the previous solve, are now available.
Further, \ipopt is instructed to reinitialize less datastructures if the structure of the NLP did not change since the last solve.
\item When \ipopt requests an evaluation of the Jacobian or Hessian, function reevaluation is now skipped if possible.
\item When \ipopt stops at a point that it claims to be locally infeasible, it is now checked whether the solution proves infeasibility, see Berthold and Witzig~\cite[Theorem 1]{BertholdWitzig2021}.
If that is not the case, the solution status is changed to ``unknown''.
\item A few \ipopt parameters can now be set directly via SCIP parameters (\texttt{nlpi/ipopt/*}).
\item Due to changes in how the \ipopt output is redirected into the SCIP log, the \ipopt banner was no longer printed reliably for the first run of \ipopt anymore.
Therefore, the banner has now been disabled completely.
\end{itemize}

\subsubsection{Interface to Algorithmic Differentiation}
\label{sect:exprint}

For the computation of first and second derivatives, SCIP traditionally relied on a third-party automatic differentiation (AD) library.
With the new expressions framework (Section~\ref{sect:expr}), first derivatives and Hessian-vector products are available in SCIP itself.
Their implementation relies on the \texttt{BWDIFF}, \texttt{FWDIFF}, and \texttt{BWFWDIFF} callbacks of the expression handlers.
The latter two are not implemented for every expression handler so far.
However, some NLP solvers make use of full Hessians and their sparsity pattern, something that is not available in the expressions framework itself yet.
Further, the current datastructure for expressions with its many pointer-redirections does not perform too well when a fixed expression needs to be evaluated repeatedly in many points.
Therefore, a separate AD library is still used in the interfaces to NLP solvers.

Currently, the only library that is interfaced is \cppad\footnote{\url{https://github.com/coin-or/CppAD}}.
In the \cppad interface, a given expression is compiled into the serial datastructure (the ``tape'') that is used by \cppad.
Here, expression types (i.e., which handler is used) are checked and translated into a form that is native to \cppad when possible.
Since the \cppad interface is used by NLPIs only, it only supports the ``varidx'' expression and not the ``var'' expression (see begin of previous section).
With SCIP 8, \cppad's feature to optimize the tape has been enabled.

Mapping of expression handlers to \cppad's operator types is available for all expression handler that are included in SCIP.
For some expression types, such as signpower, this translation has been improved to avoid repeated recompilation of an expression.  
For expression handlers that are not known to the \cppad interface, the backward- and forward-differentiation callbacks of the expression handler are used to provide first derivatives.
However, second derivatives (Hessians) are not yet available. In the \ipopt interface, the Hessian approximation will be activated in this case.

With SCIP 7, quadratic functions, including their derivatives, were treated differently from other nonlinear function.
Further, the NLPs to be solved were build from the extended, thus sparse, formulation~\eqref{eq:minlp_extlp}.
Therefore, nonlinear functions typically depended on only a few variables and, thus, it was usually sufficient to work with dense Hessians.
With SCIP 8, though, also the derivatives of quadratics are computed by the AD library and the NLPs to be solved are closer to the original form~\eqref{eq:minlp}.
For these reasons, \cppad's routines to compute sparse Hessians are used now unless more than half of the Hessian entries are nonzero.


\subsection{Performance Impact of Updates for Nonlinear Constraints}
\label{sect:perfconsexpr}


While Section~\ref{sect:perfminlp} compared the performance of SCIP 7.0 and SCIP 8.0 on a set of MINLP instances, this section takes a closer look on the effect of replacing only the handling of nonlinear constraints in SCIP.
That is, here the following two versions of SCIP are compared:
\begin{description}
\item \textbf{classic}: the main development branch of SCIP as of 23th of August 2021; in this version, nonlinear constraints are handled as it has been in SCIP 7.0, with just a few bugfixes added;
\item \textbf{new}: as classic, but with the handling of nonlinear constraints replaced as detailed in this section and symmetry detection extended to handle nonlinear constraints (see Section~\ref{sect:symmetrynl}).
\end{description}

For this comparison, SCIP has been build with GCC 7.5.0 and uses \papilo~1.0.2 for MILP presolves, \bliss~0.73 to find graph automorphisms, \cplex~20.1.0.1 as LP solver, \ipopt~3.14.4 as NLP solver, \cppad~20180000.0 for automatic differentiation, and Intel MKL 2020.4.304 for linear algebra (LAPACK).
\ipopt uses the same LAPACK and HSL MA27 as linear solver.
All runs are carried out on identical machines with Intel Xeon CPUs E5-2660 v3 @ 2.60GHz and 128GB RAM in a single-threaded mode.
As working limits, a time limit of one hour, a memory limit of 100000MB, an absolute gap tolerance of $10^{-6}$, and a relative gap tolerance of $10^{-4}$ are set.
All 1678 instances of \minlplibtwo (version 66559cbc from 2021-03-11) that can be handled by both the classic and the new version are used.
It is noted that \minlplibtwo is not designed to be benchmark set, though, since, for example, some models are overrepresented with a large number of instance.
For each instance, two additional runs where the order of variables and constraints were permuted by SCIP were conducted.
Thus, in total 5034 jobs were run for each version of SCIP.

\begin{table}
 \centering
 \begin{tabular}{lrlrrr}
  \toprule
  Subset & instances & metric & classic & new & both \\ \midrule
  all & 5034 & solution infeasible & 481 & 49 & 20 \\
    & & failed & 143 & 70 & 18 \\
    & & solved & 2929 & 3131 & 2742 \\
    & & time limit & 1962 & 1833 & 1598 \\
    & & memory limit & 0 & 0 & 0 \\
  clean & 4839 & fastest & 3733 & 3637 & 2531 \\
    & & mean time & 75.9s & 70.3s \\
    & & mean nodes & 2543 & 2601 \\
  $[0,3600)$ & 2742 & fastest & 1990 & 1697 & 945 \\
    & & mean time & 4.7s & 5.4s \\
    & & mean nodes & 415 & 455 \\
  $[10,3600)$ & 985 & fastest & 618 & 554 & 187 \\
    & & mean time & 55.6s & 66.0s \\
    & & mean nodes & 3960 & 4502 \\
  $[100,3600)$ & 484 & fastest & 292 & 262 & 70 \\
    & & mean time & 185.3s & 231.9s \\
    & & mean nodes & 12620 & 17150 \\
  $[1000,3600)$ & 141 & fastest & 72 & 81 & 12 \\
    & & mean time & 803.5s & 623.5s \\
    & & mean nodes & 43345 & 39014 \\
  \bottomrule
 \end{tabular}
 \caption{Comparison of performance of SCIP with classic versus new handling of nonlinear constraints on MINLPLib.}
 \label{tab:minlplib_status}
\end{table}

Table~\ref{tab:minlplib_status} summarizes the results.
A run is considered as failed if the reported primal or dual bound conflicts with best known bounds for the instance, the solver aborted prematurely due to a fatal error (for example, failure in solving the LP relaxation of a node), or the solver did not terminate at the time limit.
For this comparison, runs where the final solution is not feasible are accounted separately.
One can observe that with the new version, for much fewer instances the final incumbent is not feasible for the original problem, that is, the issue discussed in Section~\ref{sect:nlmotivation} has been resolved for nonlinear constraints.
For the remaining 49 instances, typically small violations of linear constraints or variable bounds occur.
Further, the reduction in ``failed'' instances by half shows that the new version is also more robust regarding the computation of correct primal and dual bounds.
Finally, we see that the new version solves about 400 additional instances than the classic one, but also does no longer solve about 200 instances within the time limit.

Subset ``clean'' refers to all instances where both versions did not fail, i.e., either solved to optimality or stopped due to the time limit.
We count a version to be ``fastest'' on an instance if it is not more than 25\% slower than the other version.
Mean times were computed as explained in the beginning of Section~\ref{sect:performance}.
Due to the increase in the number of solved instances, a reduction in the mean time with the new version on subset ``clean'' can be observed, even though the new version is fastest on less instances than the classic one.

For the remaining subsets, $[t_1,t_2)$ refers to all instances where at least one version ran for $t_1$ or more seconds and both versions terminated in less than $t_2$ seconds.
That is, only instances that could be solved to optimality by both versions are considered.
For most of these subsets, the new version is still slower more often and on average than the classic version.
Further, for a third of the instances that can be solved, both versions perform similar.
Only on the (rather small) subset $[1000,3600)$ of difficult-but-solvable instances does the new version improve.

\begin{figure}
 \pgfplotsset{ legend style = { font = \small, rounded corners = 2pt } }
 \begin{tikzpicture}
  \begin{semilogxaxis}[
      name = solvedprof,
      width = 0.48\linewidth, height = 0.48\linewidth,
      title = {Time to solve},
      xmin = 1, xmax = 200,
      xlabel = {time factor to best ($\tau$)},
      axis y line = left,
      axis y line* = box,
      ylabel = {\#instances},
      ymin = 0, ymax = 5022,
      ytick = {500, 1000, 1500, 2000, 2500, 3000, 3500, 4000, 4500, 5000},
      ymajorgrids,
      legend pos = south east
      ]

      \addplot[thick, dashed] table[x=ratio, y=classic] {relpp_classicvsnew_Time.txt};
      \addlegendentry{classic};

      \addplot[thick] table[x=ratio, y=new] {relpp_classicvsnew_Time.txt};
      \addlegendentry{new};
  \end{semilogxaxis}
 \end{tikzpicture}
 \begin{tikzpicture}
  \begin{semilogxaxis}[
      width = 0.48\linewidth, height = 0.48\linewidth,
      title = {Gap at termination},
      xmin = 1, xmax = 200,
      xlabel = {gap factor to best ($\tau$)},
      axis y line = right,
      axis y line* = box,
      ymin = 0, ymax = 5022,
      ytick = {500, 1000, 1500, 2000, 2500, 3000, 3500, 4000, 4500, 5000},
      ymajorgrids,
      legend pos = south east
      ]

      \addplot[thick, dashed] table[x=ratio, y=classic] {relpp_classicvsnew_Gap.txt};
      \addlegendentry{classic};

      \addplot[thick] table[x=ratio, y=new] {relpp_classicvsnew_Gap.txt};
      \addlegendentry{new};
  \end{semilogxaxis}
 \end{tikzpicture}
 \caption{Performance profiles on time and gap at termination of SCIP with classic versus new handling of nonlinear constraints on MINLPLib. For a specific $\tau$, the ordinate shows the number of instance for which the corresponding version of the solver was at most this much worse (regarding time (left) or gap at termination (right)) as the best of both versions.
 For the time plot, all instances that were solved to optimality are considered for each version.
 For the gap plot, all instances that did not fail are considered for each version.}
 \label{fig:minlplib_perfprofile}
\end{figure}
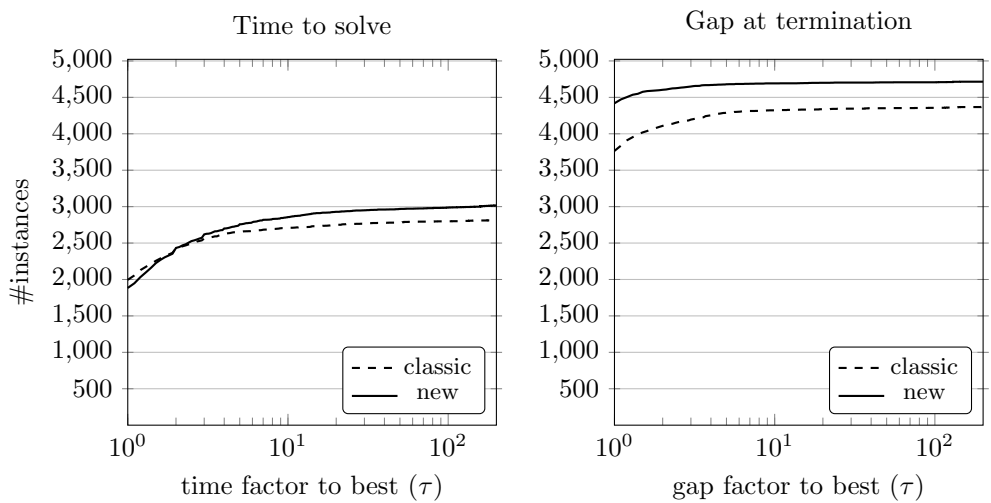

Figure~\ref{fig:minlplib_perfprofile} shows performance profiles that compare both versions w.r.t.\ the time to solve an instance and the gap at termination.
The time comparison visualizes what has already been observed in Table~\ref{tab:minlplib_status}: the new version solves more instances, but can be slower.
The gap comparison shows that on instances that are not solved, often the new version gives a smaller optimality gap than the classic version.

Appendix~\ref{appendixsect:consexpr} provides detailed results on the performance of both SCIP versions on the considered MINLPLib instances.
Further, information on the usage of the nonlinear handlers and separators that were described in this section is given.

\section{SoPlex}
\label{sect:soplex}

\subsection{Integration of \papilo in \soplex}
\label{subsect:spx_integration_papilo}

As described in Section~\ref{subsect:papilo_dual_postsolve}, version 2.0.0 of \papilo supports postsolving of dual LP solutions and basis information. This makes it possible to integrate \papilo fully as a presolving library into \soplex.  In version~6.0 of \soplex, \papilo is available as an additional option for presolving.  The previous presolving implementation continues to be the default.

The \papilo integration is handled similarly as in \scip \cite{SCIP7}.
\soplex calls a newly added simplifier plugin that converts the current problem to \papilo's data structure and then calls the presolve routine. The changes from \papilo's reduced problem are communicated by deleting the current matrix in \soplex and subsequently recreating it from \papilo's reduced problem, if the number of columns or rows decreased.

\subsection{Technical Improvements}

Several other smaller changes and improvements have been made in \soplexv. First, \soplex was extended by a C interface, as explained in Section \ref{sect:interfaces}.
Second, a rework of the internal data structures was necessary to fix warnings that were issued by current compiler versions.
Third, the dependency on the Boost program options library has been removed and the command line interface has been restored to its classic version.

Furthermore, it is now possible to use SoPlex rational solving mode without linking a GMP library, using Boosts internal implementation of rational numbers.
Finally, an ongoing LP solve of \soplex can now be interrupted from a different thread, by calling the \texttt{setInterrupt} function of SoPlex. On the \scip side, this is handled by calling \texttt{SCIPinterruptLP}.



\section{\papilo}
\label{sect:papilo}

\papilo, a C++ library, provides presolving routines for MILP and LP problems and was introduced with \scipopt~7.0 \cite{SCIP7}.
\papilo's transaction-based design generally allows presolvers to run in parallel without requiring expensive copies of the problem and without special synchronization in the presolvers themselves.
Instead of applying the results immediately, presolvers return their reductions to the core, where they are applied in a deterministic, sequential order. Modifications in the data structure are tracked to avoid applying conflicting reductions. These conflicting reductions are discarded.

The main new feature in \papilov is support for postsolving dual and basis information, which is described in Section~\ref{subsect:papilo_dual_postsolve}. This feature allows to use \papilo as an integrated presolving library in \soplex, see Section~\ref{subsect:spx_integration_papilo}.  Furthermore, \papilov comes with several improvements to the existing code base and presolving routines, described in  Section~\ref{subsect:papilo_minor_improvements}.  These changes result in a five percent improvement in the runtime when compared to the previous release.

\subsection{Postsolving Dual LP Solutions and Basis Information}
\label{subsect:papilo_dual_postsolve}

After removing, substituting, and aggregating variables from the original problem during presolving, the reduced problem (and solution) does not contain any information on missing variables.
To restore the solution values of these variables and obtain a feasible original solution, corresponding data needs to be stored during the presolving process.
The process of recalculating the original solution from the reduced one is called postsolving or post-processing~\cite{AchterbergBixbyGuetal.2019}.

Until version 1.0.2, \papilo supported only postsolving primal solutions.
In the latest version, \papilo supports postsolving also for the dual solutions, reduced costs, the slack variables of the constraints, and the basic status of the variables and constraints for the presolvers:
\texttt{DominatedColumns}, \texttt{Dualfix}, \texttt{ParallelCols}, \texttt{ParallelRows},\linebreak \texttt{Propagation}, \texttt{FixContinuous}, \texttt{ColSingleton}, and \texttt{SingletonStuffing}.
These form the majority of the LP presolvers.
The remaining presolvers are either only active in the presence of integer variables\footnote{Presolvers only active for MILP: \texttt{CoefficientStrenghtening}, \texttt{ImpliedInt}, \texttt{Probing},\texttt{SimpleProbing}, \texttt{SimplifyInequalities}}
or need to be disabled by the user\footnote{LP-Presolvers not supporting dual postsolving: \texttt{DualInfer}, \texttt{SimpleSubstitution}, \texttt{Substitution}, \texttt{Sparsify}, \texttt{ComponentDetection}, \texttt{LinearDependency}, see also settings file \texttt{lp\_presolvers\_with\_basis.set} in the \papilo repository}.

Furthermore, in dual postsolve mode \papilo only applies variable bound tightenings when they fix a variable. 
Otherwise, the solution to the reduced problem may correspond to a non-vertex solution in the original space and simple postsolving without an expensive crossover may not be possible.
If the basic information is irrelevant for the user, the variable tightening without fixing can be turned on by setting the parameter \param{calculate\_basis\_for\_dual} to false. An exception here is if a variable is unbounded. In this case, the bound of this variable is set to a finite value, which is slightly worse than the best possible bound so that the bound can not be tight in the reduced problem. 
This applies only to instances with no integer variables.
Variable tightening is still performed for mixed-integer programs.

For primal postsolving only information about removed, substituted and aggregated variables needs to be tracked. By contrast, dual postsolving needs to be informed about every modification found during presolving. \papilov keeps tracks of these changes and saves them in the postsolve stack analogously to primal postsolving. For example, a row-bound change can lead to changes in the dual solution due to complementary slackness.

After postsolving, \papilo checks if the original solution passes the primal and dual feasibility checks and fulfills the Karush-Kuhn-Tucker conditions \cite{kuhn1951nonlinear} for LP. The result of the checks is logged to the console. Since also infeasible solutions can be postsolved, \papilo does not abort if the checks fail and instead returns the result to the calling method.

For debugging purposes, this check can be performed after every step in the postsolve process. To activate this debugging feature, \papilo needs to be built in debug mode and the parameter \param{validation\-\_after\-\_every\-\_postsolving\_step} has to be turned on. This may be expensive because the problem at the current stage needs to be calculated from the original problem by applying all reductions until this point.

The introduction of dual postsolving allows using \papilo as presolving library in \soplex. Section~\ref{subsect:spx_integration_papilo} contains a brief description of the integration.

\subsection{Further Improvements}
	\label{subsect:papilo_minor_improvements}
	In this section we describe several smaller improvements in \papilov. These changes affect mostly only the performance of \papilo and rarely change the resulting reduced problem. All in all, these changes improve performance of \papilo since the last release by about five percent performance (using 16 threads) in terms of runtime and number of presolving rounds, see Table~\ref{tbl:papilo_rubberband_table} for details.
\begin{itemize}
	\item 
		When \papilo 1.0 is run with only one thread, the presolvers are executed in sequential order, but the reductions of every presolver are only applied at the end of the presolving round.  This is part of the parallel design of \papilo and helps to guarantee deterministic results independently of the number of threads used.  However, in sequential mode this does not guarantee best performance.
		
		Instead, when \papilov is run with only one thread, the reductions are applied before the next presolver, so that the next presolver can work on the modified problem.
		This feature can be turned off by setting the parameter \texttt{presolve.apply\-\_results\_\-immediately\_\-if\_\-run\_\-sequentially} to false.
	\item 
		 \texttt{DualFix} handles an additional case with two conditions: first, the objective value of the variable is zero; second, if the variable has only up-/down-locks, the lower-/upperbound is (negative) infinity. Then, the variable can be set to infinity and deleted from the model. \papilo removes the variable and marks all constraints containing the variable as redundant.
		 In postsolving, the variable is set to the maximum/minimum value such that the variable bounds and the constraints in which it appeared in the original problem are not violated and hence, the solution stays feasible.
    \item
    	 \papilo uses a transaction-based design to allow parallelization within the presolvers. This may generate conflicts when applying the reductions of the presolvers to the core. Conflicting reductions need to be discarded since it can not be ensured that the reduction is still valid.
         Conflicts make additional runs necessary to check if the discarded or a reformulated reduction can still be applied. 
         Therefore, we performed a detailed analysis of the most prominent conflict relationships, introduced a new reduction type in \papilov, and rearranged the order of presolving reductions reductions. In more detail, the improvements are as follows:
		\begin{itemize}
			\item 
				\texttt{Parallel\-Row\-De\-tection} could generate unnecessary conflicts mainly for \texttt{Parallel\-Col\-De\-tection}. To avoid these reductions and additional runs, two new reduction types \texttt{RHS\_\-LESS\_RESTRICTIVE} and \texttt{LHS\_LESS\_RESTRICTIVE} were introduced. In contrast to \texttt{RHS} and \texttt{LHS}, the columns of a row are not marked as modified, if the initial bound was (negative) infinity. 
			\item 
				\texttt{ParallelRowDetection}, \texttt{ParallelColDetection} and, \texttt{DominatedCol} could generate internal conflicts, if multiple rows/columns were parallel or dominating each other. To avoid these conflicts, bunches of parallel and dominating columns/rows are handled separately.
			\item 
				The order in which the reductions are applied to the core impacts the number of conflicts between the presolvers. We analyzed the conflicts between the presolvers and implement a new default order that minimizes the conflicts between the presolvers.
         \end{itemize}                  
      The positive impact of these changes can be observed in the reduced number of rounds reported in Table~\ref{tbl:papilo_rubberband_table}.
         	\item 
		\texttt{SimpleSubstitution} handles an additional case to detect infeasibility faster.
	\item 
		The loops in which the presolvers \texttt{Con\-straintPropagation}, \texttt{DualFix}, \texttt{Simplify\-Inequality}, \texttt{Coefficient\-Strengthening}, \texttt{SimpleSubstitution}, \texttt{SimpleProbing}, and \texttt{Implied\-Integer} scan the rows or columns of the problems were parallelized. Hence, these presolvers can distribute their workload on different threads and exploit multiple threads internally.
\end{itemize}

\begin{table}
	\caption{Performance comparison for \papilo on \miplib2017 benchmark}
	\label{tbl:papilo_rubberband_table}
	\scriptsize
	
	\begin{tabular*}{\textwidth}{@{}l@{\;\;\extracolsep{\fill}}ccccccc@{}}
		\toprule
				&           & \multicolumn{2}{c}{\papilo 1.0.3} & \multicolumn{4}{c}{\papilo 2.0.0} \\
		\cmidrule{3-4}\cmidrule{5-8}
		subset  & instances &  time [s] 	& rounds	 &  time [s] &  relative & rounds & relative \\
		\midrule
		\bracket{0}{tilim}    & 240 & 0.294557 & 19.57 & 0.281663 & 0.956 & 18.50 & 0.945 \\
		\bracket{0.01}{tilim} & 206 & 0.349799 & 22.15 & 0.334085 & 0.955 & 20.91 & 0.944 \\
		\bracket{0.1}{tilim}  & 111 & 0.693586 & 33.02 & 0.660184 & 0.952 & 30.88 & 0.935 \\
		\bracket{1}{tilim}    &  26 & 2.633462 & 43.69 & 2.477353 & 0.941 & 40.73 & 0.932 \\
		\bottomrule
	\end{tabular*}
	on a Intel Xeon CPU E5-2690 v4 @ 2.60GHz, 128GB - using 16 threads

\end{table}

Finally, two further features were introduced, improving transparency as well as the ability to debug:
\begin{itemize}
\item
	For analysis and debug purposes, \papilo can now log every transaction in the order they were applied to the problem  if verbosity level \texttt{kDetailed} is specified.
\item 
	\papilo provides an additional way to validate its correctness. A feasible debug solution can be passed via the command-line parameter \texttt{-b}. After presolving the corresponding instance, \papilo checks if the debug solution is still contained in the reduced problem and if the reduced solution can be postsolved to the same solution passed via command-line. It is recommended to turn off presolvers that use duality reasoning to (correctly) cut off optimal solutions.
\end{itemize}

\section{Interfaces}
\label{sect:interfaces}

SCIP is available via interfaces to several programming languages.
These interfaces allow users to programmatically call SCIP with an API close to the C one or leverage a higher-level syntax.
The following interfaces are available:
\begin{itemize}
\item The Python interface PySCIPOpt, which can now also be installed as a Conda package;
\item The AMPL interface that comes as part of the main SCIP library and executable;
\item The Julia package SCIP.jl;
\item C wrapper for SoPlex;
\item A Matlab interface.
\end{itemize}
We highlight below the main changes and development on interfaces to SCIP.


\subsection{AMPL}
\label{sect:ampl}


The AMPL interface of SCIP has been rewritten and moved from being a separate project (\texttt{interfaces/ampl}) to being a part of the main SCIP library and executable (\texttt{src/scip}).
The interface consists of a reader for \texttt{.nl} files as they are generated by AMPL and a specific AMPL-mode for the SCIP executable.

The \texttt{.nl} reader now relies on ampl/mp\footnote{\url{https://github.com/ampl/mp}} instead of the AMPL solver library (ASL) to read \texttt{.nl} files.
Required source files of ampl/mp are redistributed with SCIP.
Therefore, building the \texttt{.nl} reader and the AMPL interface is enabled by default.
The \texttt{.nl} reader supports linear and nonlinear objective functions and constraints, continuous, binary, and integer variables, and special-ordered sets.
More than one objective function is not supported by the interface.
A nonlinear objective function is reformulated into a constraint.
In nonlinear functions, next to addition, subtraction, multiplication, and division, operators for power, logarithm, exponentiation, sine, cosine, and absolute value are supported.
Variable and constraint flags (initial, separate, propagate, and others) can be set via AMPL suffixes.

If the SCIP executable is called with \texttt{-AMPL} as second argument, it expects the name of a \texttt{.nl} file (with \texttt{.nl} extension excluded) as first argument.
In this mode, a SCIP instance is created, a settings file \texttt{scip.set} is read, if present, the \texttt{.nl} file is read, the problem is solved, an AMPL solution file (\texttt{.sol}) is written, and SCIP exists.
Two additional parameters are available in the AMPL mode: boolean parameter \texttt{display/statistics} allows to enable printing the SCIP statistics after the solve; string parameter \texttt{display/logfile} allows to specify the name of file to write the SCIP log to.
If the problem is an LP, SCIP presolve has not run, and the LP was solved, then a dual solution is written to the solution file, too.

\subsection{Julia}\label{subsect:julia}

The Julia package \texttt{SCIP.jl} has been in development since SCIP 3
with several improvements since SCIP 7.
It contains a lower-level interface matching the SCIP public C API and a higher-level
interface based on \texttt{MathOptInterface.jl} (MOI)\cite{MathOptInterface-2021}.
The lower-level interface is automatically generated with \texttt{Clang.jl}
to match the public SCIP C API, allowing for the direct conversion of C programs using SCIP into Julia ones.
\texttt{MathOptInterface.jl} is a uniform interface for constrained structured optimization in Julia.
Solvers specify the types of constraints they support and implement those only.
Users can use the common interface for multiple solvers across different classes of problems including LP, MILP, (mixed-integer) conic optimization problems.
Higher-level modeling languages such as \texttt{JuMP.jl} are implemented on top of MOI, allowing practitioners to define their optimization model in a syntax close to the mathematical
specification and solve it through SCIP or swap solvers in a single line.

The \texttt{SCIP.jl} package can also automatically download the appropriate compiled binaries for SCIP and some of its dependencies on some platforms.
This removes the need for users to download and compile SCIP separately. Custom SCIP binaries can still be passed to the Julia package when building it.
This integration was made possible by cross-compiling SCIP through the \texttt{BinaryBuilder.jl} infrastructure, creating binaries for multiple combinations
of OS, architecture, C runtime, and compiler. The binaries are available through GitHub and versioned for other platforms to use outside of Julia.


\subsection{C Wrapper for SoPlex}\label{subsect:cwrapper}

With SCIP 8, there also comes a C wrapper for SoPlex. Since in some environments it is
much easier to interface with C code than it is with C++ (which is the language
SoPlex is written in), this wrapper paves the way for other projects to use
SoPlex as a standalone LP solver and not only through SCIP. By building a pure
C, simple shared library and header file, it is now possible to easily call
SoPlex through the foreign function interface from many other languages.


\subsection{Matlab}\label{subsect:matlab} 

In the past, two interfaces from Matlab to \scip existed. \scip came with a
rudimentary Matlab interface and there was the \texttt{OPTI Toolbox} by
Jonathan Currie, available at
\url{https://github.com/jonathancurrie/OPTI}. However, the development of
the OPTI Toolbox stopped. In order to retain the advantages of this
interface, a new interface was based on it. This new interface is available
through the Git repository
\begin{quote}
  \url{https://github.com/scipopt/MatlabSCIPInterface}.
\end{quote}
The following changes have been implemented:
\begin{itemize}
\item The new interface also runs under Linux and MacOS.
\item It works for Octave (but note that at the time of writing there is a
  bug in Octave that blocks the usage of the nonlinear part).
\item The support of other solvers has been stripped away, but the
  interface to \scip has been revised and adapted to \scipv.
\item The interface now also fully works for \scipsdp.
\end{itemize}
For the installation details, we refer to the corresponding web
page. However, the installation should only require calling
\begin{center}
\texttt{matlabSCIPInterface\_install.m}\quad or\quad 
\texttt{matlabSCIPSDPInterface\_install.m}
\end{center} from Matlab or Octave. Then the
interface will be build and you are possibly asked where to find the \scip or
\scipsdp installation (you can also supply this information through the
environment variables \texttt{SCIPDIR}/\texttt{SCIPOPTDIR} or
\texttt{SCIPSDPDIR}).

To highlight the advantages of this interface, we briefly show an
example. To solve the NLP
\[
\min_{x \in \R^2}\; \{(x_1 - 1)^2 + (x_2 - 1)^2 : 0 \leq x_1, x_2 \leq 2\},
\]
one can use the symbolic definition of the objective function:
\begin{verbatim}
  obj = @(x) (x(1) - 1)^2 + (x(2) - 1)^2;
\end{verbatim}
The remaining code is
\begin{verbatim}
lb = [0.0; 2.0];
ub = [0.0; 2.0];
x0 = [0.0, 0.0];
Opt = opti('obj',obj,'lb',lb,'ub',ub)
[x,fval,exitflag,info] = solve(Opt,x0)
\end{verbatim}
More examples are available in the repository.

\section{ZIMPL}
\label{sect:zimpl}

\zimpl 3.5.0 is released with the release of SCIP Optimization Suite 8.0.
\zimpl now allows also nonlinear objective functions.
There has been quite some work to increase the code quality further by
augmenting the code with compiler attributes.
Also the work has been
started to completely switch to C99 declarations, i.e, move the variable
declarations from the start of the functions further inside, use local
loop variables, const as much as possible, and in general try to initialize
variables once they are created.
Getting the maximum out of gcc, clang,
clang-analyzer, and pclint is an interesting experiment, which however
clearly shows that C was never mend to be validated.

A major additional feature in \zimpl is the ability to write out suitable
instances as a Quadratic Unconstrained Binary Optimization (QUBO) problem.
Unfortunately, there is no standard format for QUBO files yet, so we support
several small varieties of a sparse format for the moment.
Furthermore, we have started to implement the ability to automatically
convert constraints into quadratic binary objective
functions~\cite{glover2019quantum}.
With this release it is just a first experimental and limited ability,
but we plan to extend this continuously with the next releases.

\section{The UG Framework}
\label{sect:ug}

UG is a generic framework for parallelizing branch-and-bound based solvers in a distributed or shared memory computing environment. 
It was designed to parallelize powerful state-of-the-art branch-and-bound based solvers (we call these ``\emph{base solvers}''. 
Originally, the base solver is a branch-and-bound based solver, but in this release, it is redefined as any solver that is being parallelized by UG) externally in order to exploit their powerful performance. 
UG has been developed over 10 years as beta versions to have general interfaces for the base solvers. 
Internally, we have developed parallel solvers for SCIP~\cite{Shinano2012, Shinano18fiber, Shinano-ParaSCIP}, CPLEX (not developed anymore), FICO Xpress~\cite{Shinano-ParaXpress}, 
PIPS-SBB~\cite{Munguia2016,Munguia2019}, Concorde\footnote{\url{https://www.math.uwaterloo.ca/tsp/concorde.html}}, and QapNB~\cite{FujiiItoKimetal.2021}. 
In addition to the parallelization of these branch-and-bound base solvers, UG was used to develop MAP-SVP~\cite{MAP-SVP}, which is a solver for the Shortest Vector Problem (SVP), and whose algorithm does not rely on branch-and-bound. 
Developers of several solvers parallelized by UG needed to internally modify the UG framework itself, since UG could not handle the base solvers directly.
Especially, a success of MAP-SVP, which updated several records of the SVP challenge\footnote{\url{http://latticechallenge.org/svp-challenge}}, motivated us to develop \emph{generalized UG}, in which all solvers developed so far can be handled by a single unified framework. 
The generalized UG is included in this version of SCIP Optimization Suite as UG version 1.0.

\begin{figure}[tb]
	\begin{center}
		\includegraphics[width=120mm]{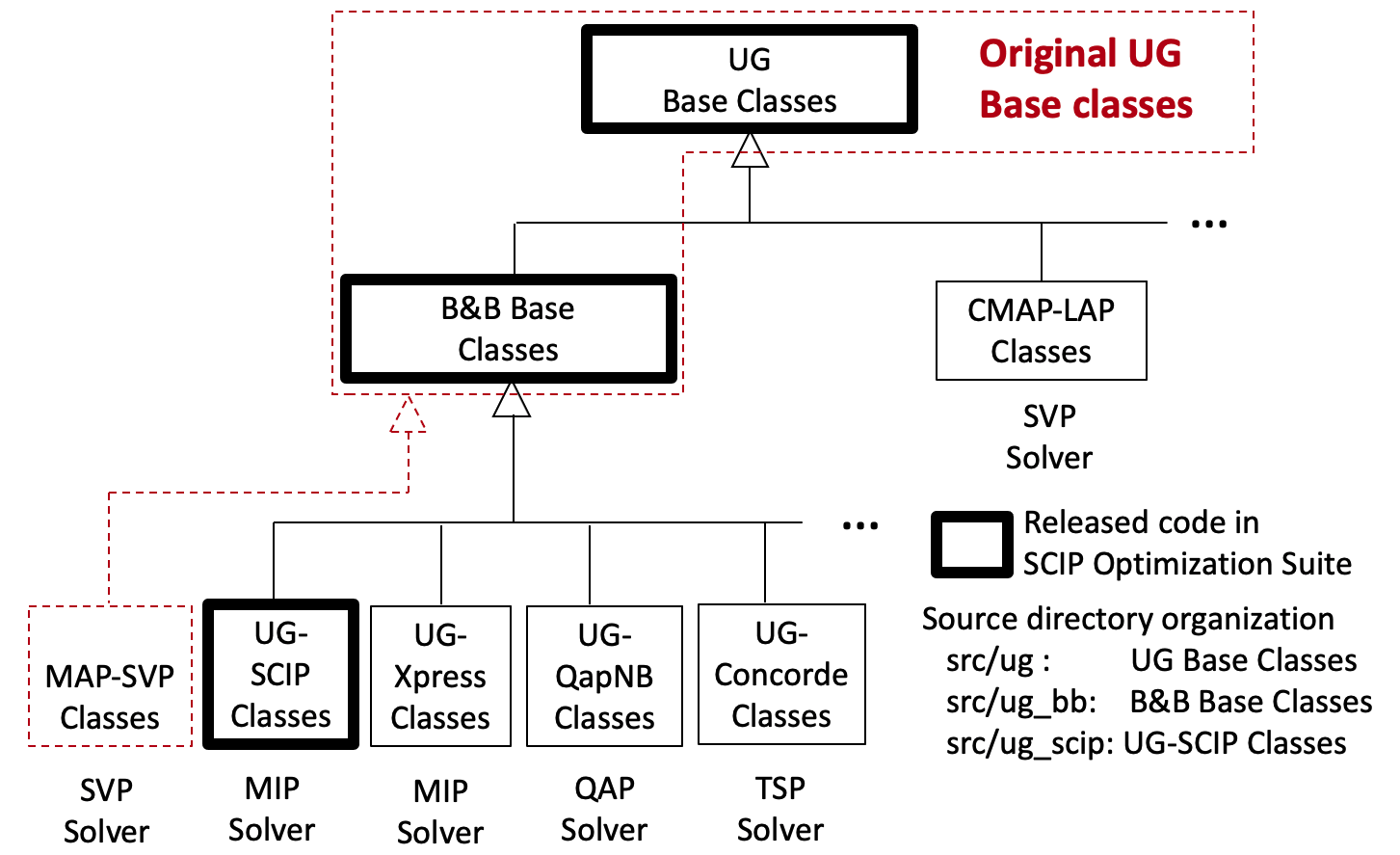}
		\vspace{-3mm}
		\caption{Class hierarchy and source code directory organization of the UG version 1.0}
		\label{fig:ug-classes-hierarchy}
	\end{center}
\end{figure}

UG version 1.0 is completely different from the previous versions internally, though its interfaces for branch-and-bound base solvers remain the same as far as possible.
Figure~\ref{fig:ug-classes-hierarchy} shows the class hierarchy of UG version 1.0. 
The original UG base classes are separated into branch-and-bound related codes and the others, so that non-branch-and-bound solvers can be parallelized naturally. 
In the original UG, the \texttt{ParaSolver} class, which wraps the ``base solver'', and the \texttt{ParaComm} class, which warps communication codes or parallelization libraries, are abstracted. 
On top of these abstractions, in UG version 1.0, the \texttt{ParaLoadCoordinator} class, which is a controller of the parallel solver, and the \texttt{ParaParamSet} class, which defines the parameter set,
are also abstracted so that a ``base solver'' specific parallel algorithm can be implemented flexibly with the ``base solver'' specific parameters. 
The flexibility of UG version 1.0 can be observed in the paper of CMAP-LAP (Configurable Massively Parallel solver framework for LAttice Problems)~\cite{TateiwaShinanoYamamuraetal.2021}, 
which is another parallel solver framework for lattice problems. On top of CMAP-LAP, CMAP-DeepBKZ~\cite{TateiwaShinanoYasudaetal.2021} has been developed, 
which is the successor of MAP-SVP and the first application of the generalized UG.

For the UG version 1.0,  proper documentation of software is started, and Doxygen style documentation is introduced. Moreover, a CMake build system is included.
In this opportunity, we made the following several modifications in FiberSCIP architecture and added a selfsplit ramp-up feature to FiberSCIP and ParaSCIP.

\subsection{Join \texttt{ParaSolver} Threads of FiberSCIP}
The {\em ramp-up} is a process which runs until all CPU cores become busy. For a general discussion of the ramp-up process for parallel branch-and-bound, see Ralphs et al.~\cite{ralphs2018parallel}
and for the ramp-up process of FiberSCIP see Shinano et al.~\cite{Shinano18fiber}.
One of the distinguishing features of FiberSCIP is racing ramp-up. 
FiberSCIP is composed of a \texttt{ParaLoadCoordinator} thread and several \texttt{ParaSolver} threads.
During the ramp-up phase, all ParaSolver threads solve the same root node with different parameter settings until certain termination criteria is met, 
that is, FiberSCIP generates multiple search trees in parallel and selects the winner within the \texttt{ParaSolver} threads, and afterward 
the winner search tree is solved in parallel.
FiberSCIP is composed of an LC thread and several ParaSolver threads.

In previous versions, the \texttt{ParaSolver} threads were detached.
The reason why the \texttt{ParaSolver} threads were detached is to enable terminating FiberSICP as soon as possible 
when one of the racing \texttt{ParaSolver} threads has solved the instance. 
When we developed FiberSCIP for the first time, it was very hard to interrupt SCIP when an LP solve is being executed. 
Therefore, all \texttt{ParaSolver} threads were detached, and the main thread exits when one of the \texttt{ParaSolver} threads has solved the instance. 
However, this mechanism leads to some instability in FiberSCIP.
The latest version of SCIP can interrupt solving appropriately while LP is running. 
In this version of FiberSCIP, all \texttt{ParaSolver} threads are joined and FiberSCIP is terminated cleanly.

\subsection{Time Limit Feature Implementation of FiberSCIP} 
In the previous versions of FiberSCIP, when a time limit is specified in the parameter, 
FiberSCIP created a \texttt{ParaTimeLimitMonitor} thread to create the time limit notification message to the \texttt{ParaLoadCoordinator}. 
The thread sleeps until the time limit, wakes up when time limit  is reached, and sends the notification message to the \texttt{ParaLoadCoordinator}. 
The \texttt{ParaLoadCoordinator} tries to interrupt all \texttt{ParaSolver} threads. 
However, these interruptions within a reasonable time could be failed when LP is running within SCIP, 
since it did not have a chance to receive the message.
From the performance point of view, creating the \texttt{ParaTimeLimitMonitor}  is not good, but \texttt{ParaLoadCoordinator} works as a kind of event-driven 
controller, and then an event to notify the time limit was needed. 

With UG version 1.0, this mechanism was changed to set a time limit when each SCIP solves a sub-MIP and hopes to detect the time limit on the \texttt{ParaSolver} side. 
How it works well depends on how accurately SCIP can terminate for the time limit setting. 
Unfortunately, it has several irregular timings and FiberSCIP needs to handle such cases currently.
However, from the performance point of view, it has a benefit of running FiberSCIP without the \texttt{ParaTimeLimitMonitor} thread.

\subsection{SelfSplit Ramp-up}
In the general terminology of parallel branch-and-bound, ramp-up is a part of the initialization phase of the computation. 
The initialization phase is summarized by Henrich~\cite{henrich93init}. The {\em enumerative initialization} broadcasts the root node to all processes,
which then perform an initial tree search according to the sequential algorithm. 
When the number of leaf nodes on each processor is at least the number of processes, processes
can stop expanding. The $i^{th}$ process then keeps the $i^{th}$ node and deletes the rest. 
In this method, all processes are working from the very beginning and no communication is required.
The {\em SelfSplit}~\cite{FISCHETTI2018101} is refined in the enumerative initialization phase so that the open nodes are ordered more accurately 
to each solver that has a similar amount of work and so that a parallel MILP solver can perform static load balancing more efficiently.
For QapNB parallelization, the enumeration initialization has a benefit, see Fujii et al.~\cite{FujiiItoKimetal.2021}.
In UG version 1.0, the initialization method is implemented in B\&B base classes as shown in Figure~\ref{fig:ug-classes-hierarchy}
so that any ``base solver'' can use it and name it SelfSplit ramp-up, that is, in the case of FiberSCIP, it is not limited to MILP, but also works for MINLP.
Unlike in the original SelfSplit, 
it is used as a ramp-up and can perform dynamic load balancing after ramp-up, 
also, in order to cooperate with the other features of UG, it does layered presolving for sub-MIPs generated by SelfSplit and the sub-MIPs are checkpointing. 
FiberSCIP and ParaSCIP included in this release have this SelfSplit ramp-up feature.

\subsection{Memory usage estimation}
To make FiberSCIP stable, one of the features needed is handling the memory limit, since memory overuse makes FiberSCIP abort. However, memory usage estimation is very hard for FiberSCIP. By using SCIP functions for memory usage estimation (plus Linux system feature of memory usage, in case of FiberSCIP runs on Linux), the FiberSCIP memory usage estimation feature is implemented. When the estimation is more than the system memory, the latest version of FiberSCIP terminates with "memory limit reached".

\section{The GCG Decomposition Solver}
\label{sect:gcg}


\pgfplotscreateplotcyclelist{timelist}{%
orange!60!black,fill=orange!40!white\\%
red!60!black,fill=red!40!white\\%
}

\pgfplotscreateplotcyclelist{bestlist}{%
olive!60!black,fill=olive!40!white\\%
olive!60!black,fill=olive!40!white,postaction={pattern=crosshatch dots}\\%
brown!60!black,fill=brown!40!white\\%
}

\pgfplotscreateplotcyclelist{nodelist}{%
orange!60!black,fill=orange!40!white,postaction={pattern=crosshatch dots}\\%
red!60!black,fill=red!40!white,postaction={pattern=crosshatch dots}\\%
}

\pgfplotstableread[row sep=\\,col sep=&]{
    setting        & geo. mean & arith. mean \\
    pseudo         & 0.44      & 0.39        \\
    random         & 0.71      & 0.53        \\
    mostfrac       & 0.49      & 0.4         \\
    SBw/oCG        & 0.73      & 0.47        \\
    SBw/CG         & 1         & 1           \\
    hierarchical   & 0.54      & 0.47        \\
    hybrid         & 0.72      & 0.67        \\
    reliable       & 0.65      & 0.55        \\
    reliable hier. & 0.47      & 0.48        \\
    hybrid hier.   & 0.32      & 0.27        \\
    }\origtime
    
\pgfplotstableread[row sep=\\,col sep=&]{
    setting        & time & nnodes & timeouts\\
    pseudo         & 14   & 5      & 5       \\
    random         & 8    & 4      & 5       \\
    mostfrac       & 8    & 5      & 5       \\
    SBw/oCG        & 4    & 7      & 2       \\
    SBw/CG         & 1    & 0      & 17      \\
    hierarchical   & 8    & 25     & 1       \\
    hybrid         & 1    & 18     & 7       \\
    reliable       & 4    & 11     & 4       \\
    reliable hier. & 6    & 14     & 3       \\
    hybrid hier.   & 22   & 11     & 1       \\
    }\orignbest
    
\pgfplotstableread[row sep=\\,col sep=&]{
    setting        & geo. mean & arith. mean \\
    pseudo         & 0.67      & 0.84        \\
    random         & 1         & 0.81        \\
    mostfrac       & 0.75      & 1           \\
    SBw/oCG        & 0.27      & 0.37        \\
    SBw/CG         & 0.10      & 0.32        \\
    hierarchical   & 0.20      & 0.44        \\
    hybrid         & 0.17      & 0.60        \\
    reliable       & 0.45      & 0.46        \\
    reliable hier. & 0.35      & 0.70        \\
    hybrid hier.   & 0.37      & 0.46        \\
    }\orignodes

\pgfplotstableread[row sep=\\,col sep=&]{
  setting       & geo. mean & arith. mean \\
  pseudo        & 0.61      & 0.58        \\
  hierarchical  & 0.75      & 0.70          \\
  hybrid        & 1         & 1               \\
  reliable      & 0.91      & 0.81            \\
  reliable hier.& 0.65      & 0.70         \\
  hybrid hier.  & 0.44      & 0.40        \\
  }\origadvtime
  
\pgfplotstableread[row sep=\\,col sep=&]{
  setting       & time & nnodes & timeouts\\
  pseudo        & 22   & 8      & 5       \\
  hierarchical  & 8    & 22     & 1       \\
  hybrid        & 1    & 20     & 7       \\
  reliable      & 7    & 13     & 4       \\
  reliable hier.& 11   & 14     & 3       \\
  hybrid hier.  & 24   & 10     & 1       \\
  }\origadvnbest
  
\pgfplotstableread[row sep=\\,col sep=&]{
  setting       & geo. mean & arith. mean \\
  pseudo        & 1         & 1           \\
  hierarchical  & 0.30      & 0.53        \\
  hybrid        & 0.25      & 0.71        \\
  reliable      & 0.67      & 0.55        \\
  reliable hier.& 0.53      & 0.83        \\
  hybrid hier.  & 0.56      & 0.55        \\
  }\origadvnodes

\gcg is an extension that turns \scip into a decomposition-based solver
for mixed-integer linear programs. While the focus of \gcg is
on Dantzig-Wolfe reformulation (DWR) and Lagrangian decomposition,
Benders decomposition (BD) is also supported. The philosophy
behind \gcg is that decomposition-based algorithms like branch-price-and-cut
(BP\&C) can be routinely applied to MILPs without the user's
interaction or even knowledge, just like branch-and-cut. To this end,
\gcg automatically detects a model structure that admits a
decomposition, and then performs the corresponding reformulation.
This results in a master problem and one or several subproblems, which
are usually formulated as MILP problems. The latter are solved as sub-\scip{s} or
using specialized solvers. Based on the reformulation, the linear
relaxation in every node is solved by column generation (in the DWR
case), respectively, Benders cut generation (in the BD case).  \gcg
features primal heuristics and separation of cutting planes, several
of which are adapted from \scip, but some are tailored to the
decomposition situation in which both, an original and a reformulated
model are available.

Since the last major release 3.0 in 2018, most development efforts
went into improving the usability for experts and beginners alike.
Besides few algorithmic features and code refactoring, the new release
comes with more interfaces, a much improved documentation, and a
collection of tools to support computational experiments and even more
their evaluation/visualization.
Even though not visible to the user, we improved the development
process, in particular by establishing continuous integration
pipelines etc.
We like to think of \gcg~\gcgversion\ as (the first half of) an \emph{ecosystem
  release}.

\subsection{Detection Loop Refactoring}
%


Decomposition-based algorithms rely on model structures, such as a
block-angular constraint matrix.
For an automatic identification of such structures, \gcg features a modular detection loop,
which was introduced with version 3.0. So-called \emph{detectors}
iteratively assign roles like ``master'' or ``block'' to variables
and/or constraints, potentially only to subsets, and possibly in
several rounds. This way, usually many different potential decompositions are
found.  We refer to the \scip Optimization Suite 6.0 release report~\cite{SCIP6} for
a more detailed overview.  Detectors are implemented as plugins such that new ones can
be added conveniently. In every round, each detector works on existing
(but possibly empty) partial or finished decompositions. An empirically
very successful detection concept builds on the classification of
constraints and variables, which is performed prior to the actual
detection process, using so-called \emph{classifiers}. Classifiers
group variables or constraints according to arbitrary
decomposition-relevant information, like type or name. Detectors
then may or may not use the resulting classes to create (partial) decompositions.  

\gcg~\gcgversion\ comes with a refactored detection loop. Large
parts of the code base were rewritten, moved, or renamed, while the
general idea and procedure of the modular detection loop were
maintained. Consequently, the interface changed significantly. First,
(partial) decompositions, which were propagated, finished, or
post-processed by the detectors, are represented by objects of the
class \texttt{PARTIALDECOMP} (formerly \texttt{Seeed}). Objects of the
class \texttt{DETPROBDATA} (formerly \texttt{Seeedpool}) manage all
decompositions that were created during a run of the detection loop.
Since the class \texttt{DETPROBDATA} should only manage the partial
decompositions, many functions of \texttt{Seeedpool} were moved.
Formerly, \texttt{Seeedpool} provided the functionality to classify
constraints and variables. With version~\gcgversion, variable and
constraint classifiers are implemented as plugins as well. All
classifiers are called at the beginning of the detection process. Each
classifier can produce partitions of constraints or variables, which
are represented by objects of the classes \texttt{ConsPartition}
(formerly \texttt{ConsClassifier}) and \texttt{VarPartition} (formerly
\texttt{VarClassifier}), respectively. Both classes implement the
abstract class \texttt{IndexPartition} (formerly
\texttt{IndexClassifier}). This modular design allows users to easily add
classifiers.

Furthermore, many parameters were changed. For a complete overview we
refer to the \texttt{CHANGELOG}. Users can enable or disable the entire detection
process using the parameter \texttt{detection/enabled}. Moreover, the
parameter \texttt{detection/postprocess} enables or disables the
post-processing of decompositions. Parameters related to
classification were moved to \texttt{detection/classification/}. 
By default, \gcg~\gcgversion\ preprocesses an instance, runs the detection, then solves.
If detection should run already on the un-preprocessed model, it must
be initiated manually before presolve starts; by default, a second round of
detection is then still performed on the preprocessed model.
With this more intuitive behavior,
the parameters \texttt{origenabled} of detectors and classifiers were
removed.

The legacy detection mode (for detectors prior to version~3.0) is no longer available. All corresponding
parameters and legacy detectors were removed. Users have to implement
detectors using the new callbacks and API.

\subsection{Strong Branching in Branch-and-Price}
\label{sect:gcgstrongbranch}


In \gcg, two general branching rules are implemented (branching on
original variables~\cite{VilleneuveDesrosiersLuebbeckeSoumis:05} and
Vanderbeck's generic branching~\cite{vanderbeck11}) as well as one
rule that applies only to set partitioning master problems
(Ryan and Foster branching~\cite{RyanFoster81}).
While these rules differ quite significantly (creating two child nodes
vs.\ several child nodes; branching on variables vs.\ on constraints),
the general procedure at a node comes in
two common stages: First, one determines the set of candidates we could
possibly branch on (called the branching rule here).
Second, the \emph{branching candidate selection heuristic\/} 
then actually selects one of the available candidates. In \scip the latter
is done by ranking the candidates according to a score.
%
Both the branching rule and the selection
heuristic can have a significant impact on the size of the
branch-and-bound tree, and hence on the runtime of the entire
algorithm. \gcg previously contained only pseudo cost, most fractional,
and random branching as selection heuristics for original variable
branching, and first-index branching for Ryan-Foster and Vanderbeck's
generic branching. In \gcg~\gcgversion, several new selection heuristics are added, all of
which are based on strong branching. For an overview of which selection
heuristics are available for which branching rules see
Table~\ref{tab:selection_heuristics}. In the following, we briefly
describe the new selection heuristics. For more detailed descriptions,
we refer to \cite{gaul2021branching}.

\begin{table}
    \caption{Branching candidate selection heuristics 
      available in \gcg~\gcgversion.} 
    \label{tab:selection_heuristics}
    \begingroup
    \centering
    \renewcommand\arraystretch{1.1}
    \begin{tabular}{|c|c|c|c|}
         \hline
         \diagbox[innerleftsep=.1cm,innerrightsep=.1cm,
         height=2.7\line, width=4.5cm]{selection heuristic}{branching rule}                               & \multirowcell{2}[0cm][c]{original} & \multirowcell{2}[0cm][c]{Ryan-Foster}     & \multirowcell{2}[0cm][c]{Vanderbeck} \\
         \hline
         \multirowcell{2}[0cm][c]{random/index-based\\ branching} & \multirowcell{2}[0cm][c]{\checkmark} & \multirowcell{2}[0cm][c]{\checkmark}       &    \multirowcell{2}[0cm][c]{\checkmark}                  \\
         &&&\\
         \hline
         \multirowcell{2}[0cm][c]{most fractional/infeasible \\ branching}           & \multirowcell{2}[0cm][c]{\checkmark} & \multirowcell{2}[0cm][c]{\checkmark$^2$} &                    \\
         &&&\\
         \hline
         pseudocost branching                                                        & \checkmark                           & \checkmark$^2$                         &                    \\
         \hline
         \multirowcell{2}[0cm][c]{strong branching \\ with column generation$^1$}    & \multirowcell{2}[0cm][c]{\checkmark} & \multirowcell{2}[0cm][c]{\checkmark}       &                    \\
         &&&\\
         \hline
         \multirowcell{2}[0cm][c]{strong branching \\ without column generation$^1$} & \multirowcell{2}[0cm][c]{\checkmark} & \multirowcell{2}[0cm][c]{\checkmark}       &                    \\
         &&&\\
         \hline
         hybrid branching$^1$ 
         & \checkmark                           & \checkmark$^3$                           &                    \\
         \hline
         reliability branching$^1$ 
         & \checkmark                           & \checkmark$^3$                           &                    \\
         \hline
         hierarchical branching$^1$ 
         & \checkmark                           & \checkmark$^3$                           &                    \\
         \hline
    \end{tabular}

    \endgroup

    \begin{footnotesize}
      $^1$The strong branching based heuristics
      can be combined. $^2$\gcg can only aggregate the respective
      scores of the (two) individual variables. $^3$These heuristics
      originally use both strong and pseudocost branching; however,
      pseudocost branching can also be substituted by 
      any other heuristic, with varying performance.

    \end{footnotesize}

\end{table}

\subsubsection{Branching Candidate Selection Heuristics Background}

Given a set of branching candidates, the selection heuristic
usually creates a ranking and selects a winner.
%
One ranking criterion is the expected gain, that is,
the improvement in dual bound in the child nodes when compared to the
current node. However, computing the exact gains amounts to performing
\emph{all full strong branching}. In a
branch-and-price context this means evaluating all branching candidates by
solving all child node LP relaxations with column
generation to optimality. With often hard pricing problems, this
variant is an even (much) larger computational burden than it is in
the standard branch-and-cut context.
Yet, strong branching has demonstrated potential in
branch-and-price for hard instances~\cite{pecin2014branching,roepke2012branching}. 
In particular, strong branching generally creates small trees
(compared to other branching rules).

To alleviate the computational effort, relaxations of strong branching
are considered, such as not evaluating all candidates, or not solving the
LP relaxations to optimality. In branch-and-price, one has even more
degrees of freedom. In particular, we can choose to not
perform column generation when evaluating the candidates. Thus, we
differentiate between \textit{strong branching without column
  generation} (SBw/oCG) and \textit{strong branching with column
  generation} (SBw/CG).
In principle, one has the entire spectrum in between (partial or
heuristic pricing, etc.) but this is beyond the scope here.

Another  alternative to exact computations are \emph{gain predictions}. A classical
approach is \textit{pseudocost branching}.
Pseudocosts measure the average
gain per eliminated fractionality 
for down and up branch separately, with the average being calculated
from candidates which we branched on in the past that correspond to
the same variable. Pseudocosts can be calculated very fast,
however, using them also creates larger trees than strong branching. In
particular, pseudocost scores are unreliable at the top of the
tree, as there is only little/no historic data from which the
averages can be calculated.

To combine the strengths and reduce the weaknesses of strong branching
and pseudocost branching, they can also be used together. The three
options for this that were added to \gcg~\gcgversion\ are \textit{hybrid
  strong/pseudocost branching}~\cite{LINDO2003manual}, 
\textit{reliability branching}~\cite{achterberg2005branching}, and
\textit{hierarchical strong branching}~\cite{gaul2021branching,pecin2014branching,roepke2012branching}. In hybrid
strong/pseudocost branching, strong branching is applied only for
nodes up to a given depth, and pseudocost branching to the rest.

For reliability branching, each candidate is assigned a reliability
score, which is simply the minimum of the number of down and up
branches in the tree of candidates corresponding to the same variable.
The reliability score of a candidate reflects how close the pseudocost
score for the candidate is likely to be to the actual gain, that is, how
\emph{reliable} the prediction is. Candidates whose reliability score
is below a certain threshold are evaluated with strong branching,
while the remaining candidates are again evaluated using pseudocosts.

For hierarchical strong branching, the selection process is divided
into three phases (also called a hierarchy in the literature), where
\begin{itemize}
\item in phase 0, the candidates are filtered based on some heuristic
  that is quick to compute (such as pseudocost, most fractional, or
  random branching), then
\item in phase 1, the remaining candidates are filtered based on their
  SBw/oCG scores, and finally
\item in phase 2, a candidate is selected out of the remaining
  candidates based on the score the candidates received from SBw/CG.
\end{itemize}

The effort 
at a given node depends on the assumed importance of evaluating the node
precisely (a larger estimated size of the subtree
gives more importance to a candidate), and on the difference in
computational effort vs.\ quality of predictions between the phases.
The intuition here is that only the most promising candidates (based
on the scores from the earlier heuristics) should receive the largest
evaluation effort (and best evaluation quality). 

In the same way that strong branching can be combined with pseudocost
branching to obtain hybrid strong/pseudocost branching and reliability
branching, we can also combine hierarchical strong branching with
hybrid strong/pseudocost branching and reliability branching to obtain
\textit{hybrid hierarchical strong/pseudocost branching} and
\textit{hierarchical reliability branching} \cite{gaul2021branching}:
for hierarchical reliability branching, we only perform strong
branching in phase 1 and 2 on candidates that are not yet reliable,
with different thresholds for phase 1 and phase 2. For hybrid
hierarchical strong/pseudocost branching, we only perform phase 0
starting from a given depth, and phases 1 and 2 up to a given depth
(again separate thresholds for each phase).

Strong branching is implemented using \scip's probing mode. The
columns generated when evaluating a node using SBw/CG are kept in
\gcg's column pool. The (potentially positive side) effect of this needs
still to be evaluated.

\subsubsection{Parameters for Strong Branching}

By default, like in \scip, strong branching is disabled.
It can be enabled for
original variable branching or Ryan-Foster branching by setting
\texttt{branching/orig/usestrong} or
\texttt{branching/\-ryanfoster/\-usestrong}, respectively, to
\texttt{TRUE}. By default, this performs a hybrid hierarchical
branching.
Furthermore, there are several parameters
that allow to completely change the behavior of the
heuristics. In fact, all of the strong branching heuristics use
the same implementation, just with different parameter settings.
Preset settings files for each of the previously described selection
heuristics as well as a template file can be found in \gcg's
\texttt{settings} folder. Table~\ref{tab:sbpars} lists
most of the available parameters. Further information can be found in
the paper by Gaul~\cite{gaul2021branching} and in \gcg's documentation.

\begin{table}
  \caption{Parameters for strong branching.}\label{tab:sbpars}
  \renewcommand{\arraystretch}{1.2}%
  \footnotesize
  \vspace{-.3cm}
  \begin{tabularx}{\textwidth}{@{}>{\raggedright\arraybackslash}p{0.3\textwidth}@{}X@{}}
    \toprule
    \textbf{Parameter}                            & \textbf{Effect} \\
    \midrule
    \multicolumn{2}{@{}l}{\texttt{branching/[orig,ryanfoster]/\dots}} \\
    \texttt{minphase[0,1]outcands}                & minimum number of output candidates from phase [0,1] \\
    \texttt{maxphase[0,1]outcands}                & maximum number of output candidates from phase [0,1] \\
    \texttt{maxphase[0,1]outcandsfrac}            & maximum number of output candidates from phase 0 as fraction of total candidates, takes precedence over \texttt{minphase[0,1]outcands} \\
    \texttt{phase[1,2]gapweight}                  & how much influence the nodegap has on the number of output candidates from phase [1,2]$-1$ \\
    \midrule
    \multicolumn{2}{@{}l}{\texttt{branching/bpstrong/\dots}} \\
    \texttt{histweight}         & fraction of candidates in phase 0 that are chosen based on historical strong branching performance \\
    \texttt{mincolgencands}    & minimum number of candidates for phase 2 to be performed, otherwise the best previous candidate will be chosen \\
    \texttt{maxsblpiters}, \texttt{maxsbpricerounds} & upper bound on number of simplex iterations/pricing rounds, sets upper bound to twice the average if set to 0 \\
    \texttt{immediateinf}                         & if set to \texttt{TRUE}, candidates with infeasible children are selected immediately \\
    $\texttt{reevalage}$                          & reevaluation age \\
    \texttt{maxlookahead}                         & upper bound for the look ahead \\
    \texttt{lookaheadscales}    & by how much the look ahead scales with the overall evaluation effort (currently \texttt{lookaheadscales} * \texttt{maxlookahead} is the minimum look ahead) \\
    \texttt{closepercentage}    & fraction of the chosen candidate's phase 2 score the phase 0 heuristic's choice needs to have in order to be considered close \\
    \texttt{maxconsecheurclose} & number of times in a row the phase 0 heuristic needs to be close for strong branching to be stopped entirely \\
    \texttt{minphase0depth},
    \texttt{maxphase[1,2]depth},
    \texttt{depthlogweight},
    \texttt{depthlogbase},
    \texttt{depthlogphase[0,2]frac}
  & $\kappa_1^+ = \lambda_{\text{depth}}^1 + \rho_{\text{depth}} \cdot \text{log}_{\lambda_{\text{base}}}(n_{\text{cands}})$, where $\kappa_i^+$ is the depth until which phase $i$ is performed, $\lambda_{\text{depth}}^1 = \texttt{maxphase1depth}$, $\rho_{\text{depth}} = \texttt{depthlogweight}$, $\lambda_{\text{base}} = \texttt{depthlogbase}$ and $n_{\text{cands}}$ is number of variables that we could branch on (usually all integer and binary variables). $\kappa_2^+ = \kappa_1^+ \cdot \texttt{depthlogphase2frac}$, but at most \texttt{maxphase2depth}. The minimum depth from which on phase 0 is performed is equal to $\kappa_1^+ \cdot \texttt{depthlogphase0frac}$, but at least \texttt{minphase0depth} \\
    \texttt{phase[1,2]reliable}                   & min count of pseudocost scores for a variable to be considered reliable in phase [1,2] \\
    \bottomrule
  \end{tabularx}
\end{table}

\subsubsection{Performance Evaluation}

A preliminary performance comparison (for
original variable branching only) regarding the number of nodes
and computation times can be
found in the tables 
in Appendix~\ref{appendsect:gcg} and
in Figure~\ref{fig:origres}. 
The computations
were performed using an Intel Xeon L5630 Quad Core with 2.13 GHz and
16 GB of RAM. The time limit was set to 3600 seconds. 
All instances in our testset are taken from the unreleased
strIPlib (\url{striplib.or.rwth-aachen.de}). Instances contained in strIPlib are all
known to contain a model structure to which a DWR is applicable.
We selected instances that need at least 1000 nodes to solve (with a
previous \gcg version). Diversity of the testset was controlled in a
manner similar to the creation of the MIPLIB 2017 benchmark set~\cite{miplib2017}.
We refer to Gaul~\cite{gaul2021branching} for further information about the
experiments.

The preliminary results suggest that it is conceivable to have strong branching components enabled
in branch-and-price by default. A potential explanation is that
solving the subproblems is very costly, and processing a node is
more expensive than in standard branch-and-bound; 
thus there is an ever stronger incentive to have a smaller tree in branch-and-price.
Certainly, this needs more investigation.


\begin{figure}[p]
  \centering
  \vspace{-.6cm}
  \scalebox{0.8}{
  \begin{tikzpicture}
      \begin{axis}[
              ybar,
              symbolic x coords={pseudo,random,mostfrac,SBw/oCG,SBw/CG,hierarchical,hybrid,reliable,reliable hier.,hybrid hier.},
              xticklabel style={rotate=-30},
              width=1.25\textwidth,
              height=.5\textwidth,
              legend style={at={(0.25,1)},
                  anchor=north,legend columns=-1},
              ylabel={time (relative to highest entry)},
              ymin=-0.1,
              cycle list name=timelist
          ]
          \addplot table[x=setting,y=geo. mean]{\origtime};
          \addplot table[x=setting,y=arith. mean]{\origtime};
          \legend{geo. mean, arith. mean}
      \end{axis}
  \end{tikzpicture}
  }
  \scalebox{0.8}{
  \begin{tikzpicture}
      \begin{axis}[
              ybar,
              symbolic x coords={pseudo,random,mostfrac,SBw/oCG,SBw/CG,hierarchical,hybrid,reliable,reliable hier.,hybrid hier.},
              xticklabel style={rotate=-30},
              width=1.25\textwidth,
              height=.5\textwidth,
              legend style={at={(0.45,1.05)},
                  anchor=north,legend columns=3},
              ylabel={\# instances},
              cycle list name=bestlist
          ]
          \addplot table[x=setting,y=time]{\orignbest};
          \addplot table[x=setting,y=nnodes]{\orignbest};
          \addplot table[x=setting,y=timeouts]{\orignbest};
          \legend{best time, best \# nodes w/o SBw/CG,\# timeouts}
      \end{axis}
  \end{tikzpicture}
  }
  \scalebox{0.8}{
  \begin{tikzpicture}
      \begin{axis}[
              ybar,
              symbolic x coords={pseudo,random,mostfrac,SBw/oCG,SBw/CG,hierarchical,hybrid,reliable,reliable hier.,hybrid hier.},
              xticklabel style={rotate=-30},
              width=1.25\textwidth,
              height=.5\textwidth,
              legend style={at={(0.75,1)},
                  anchor=north,legend columns=-1},
              ylabel={\# nodes (relative to highest entry)},
              ymin=-0.1,
              cycle list name=nodelist
          ]
          \addplot table[x=setting,y=geo. mean,]{\orignodes};
          \addplot table[x=setting,y=arith. mean]{\orignodes};
          \legend{geo. mean, arith. mean}
      \end{axis}
  \end{tikzpicture}
  }
  \caption{Visualizations for original variable branching based on
    tables from Appendix~\ref{appendsect:gcg}. 
    The upper plot
    shows the geometric and arithmetic mean for the amount of time
    needed relative to the highest in each category. The lower plot
    shows the same for the number of nodes needed. The middle figure
    shows the number of timeouts for each heuristic, and the number of
    times a heuristic was the best for a given instance regarding
    nodes (excluding full strong branching) and
    time.}\label{fig:origres} 
\end{figure}
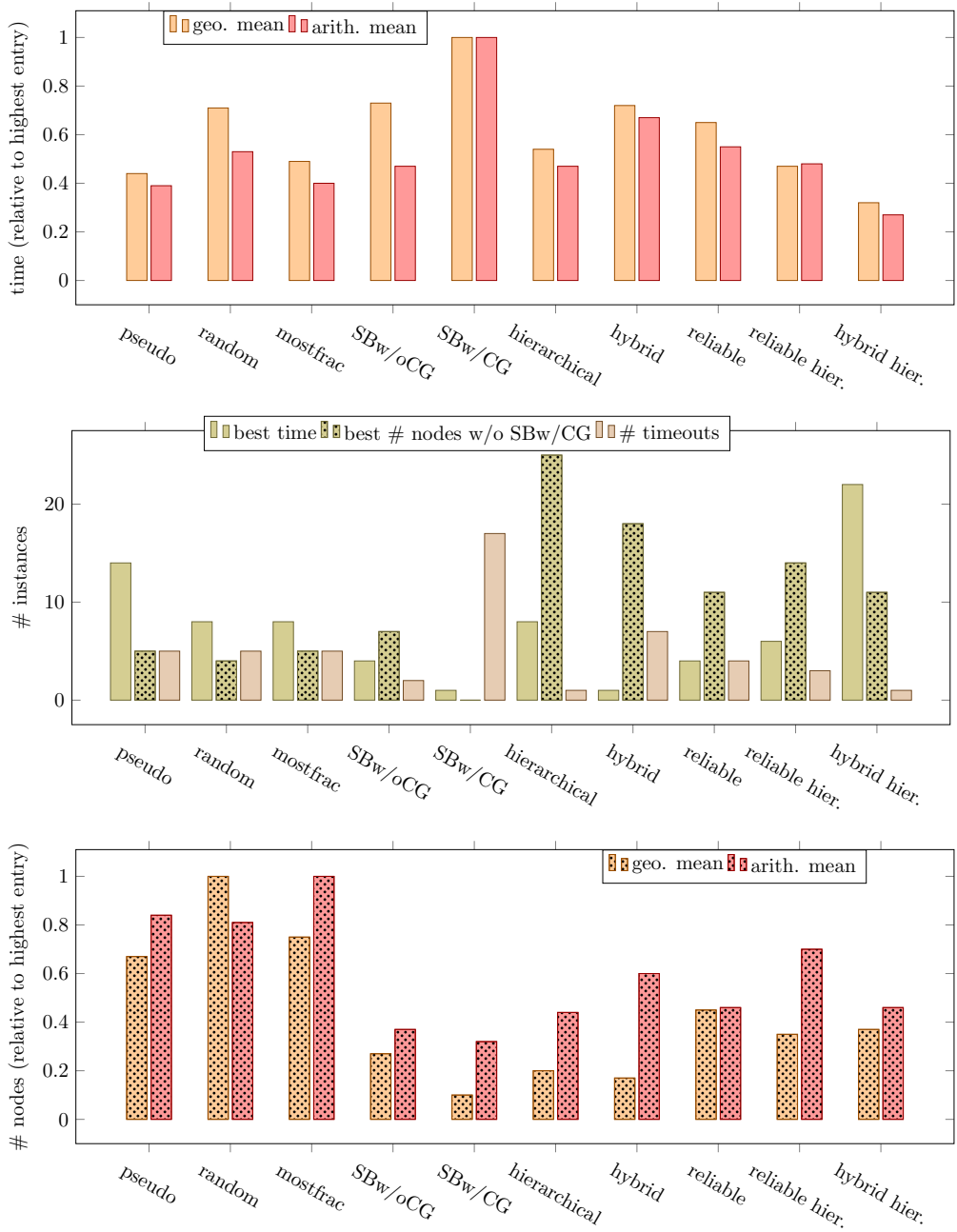

\subsection{Python Interface}
With \gcg~\gcgversion\
we introduce \textsc{PyGCGOpt} which extends \scip's existing Python
interface~\cite{MaherMiltenbergerPedrosoRehfeldtSchwarzSerrano2016}
for \gcg. It is implemented in Cython (\url{cython.org})
and is distributed  as a package independent from
the optimization suite
under  \url{github.com/scipopt/PyGCGOpt}. All the existing functionality for the
modeling 
of MILPs is inherited from \textsc{PySCIPOpt}. As a result, any MILP
modeled in Python can also be solved with \gcg without additional
effort. This lowers the technical hurdle to try out a branch-and-price
approach for any existing problem. In its first incarnation, the
interface supports specifying custom decompositions and exploration of
automatically detected decompositions. They can be visualized directly
within Jupyter notebooks. 
In addition, \gcg plugins for \texttt{detectors} and \texttt{pricing solvers} can be implemented in Python.

In the following code listing, the capacitated $p$-median problem
(CPMP) is modeled with \textsc{PySCIPOpt}'s expression syntax. The
specified textbook decomposition~\cite{Ceselli2005AProblem} is solved
by \gcg with Dantzig-Wolfe reformulation upon the call to
\texttt{m.optimize()}.
Note that the automatic structure detection functionality of \gcg
remains intact, so that the user does not need to (but can) specify a decomposition.

\begin{lstlisting}[backgroundcolor=\color{white}, basicstyle=\footnotesize\ttfamily, breakatwhitespace=false, breaklines=true, commentstyle=\color{green}, extendedchars=true, firstnumber=1, frame=single, keepspaces=true, keywordstyle=\color{blue}, numbers=left, numbersep=5pt, numberstyle=\tiny\color{gray}, rulecolor=\color{black}, showspaces=false, showstringspaces=false, showtabs=false, stepnumber=1, language=Python]
from pygcgopt import gcgModel, quicksum as qs

n_locs = 5
n_clusters = 2
distances = {0: {0: 0, 1: 6, 2: 54, 3: 52, 4: 19}, 1: {0: 6, 1: 0, 2: 28, 3: 75, 4: 61}, 2: {0: 54, 1: 28, 2: 0, 3: 91, 4: 40}, 3: {0: 52, 1: 75, 2: 91, 3: 0, 4: 28}, 4: {0: 19, 1: 61, 2: 40, 3: 28, 4: 0}}
demands = {0: 14, 1: 13, 2: 9, 3: 15, 4: 6}
capacities = {0: 39, 1: 39, 2: 39, 3: 39, 4: 39}

m = gcgModel()
x = {(i, j): m.addVar(f"x_{i}_{j}", vtype="B", obj=distances[i][j]) for i in range(n_locs) for j in range(n_locs)}
y = {j: m.addVar(f"y_{j}", vtype="B") for j in range(n_locs)}

conss_assignment = m.addConss(
  [qs(x[i, j] for j in range(n_locs)) == 1 for i in range(n_locs)])
conss_capacity = m.addConss(
  [qs(demands[i] * x[i, j] for i in range(n_locs)) <= capacities[j] * y[j] for j in range(n_locs)])
cons_pmedian = m.addCons(qs(y[j] for j in range(n_locs)) == n_clusters)

master_conss = conss_assignment + [cons_pmedian]
block_conss = [[cons] for cons in conss_capacity]
m.addDecompositionFromConss(master_conss, *block_conss)

m.optimize()

\end{lstlisting}

The Python interface required refactoring within the codebase of \gcg.
Before, a lot of core functionality of the solver was implemented
within dialog handlers. This made it hard to use \gcg as a library in
external programs. The functions \texttt{gcg\-transformProb()},
\texttt{gcg\-presolve()}, \texttt{gcg\-detect()},
\texttt{gcg\-solve()}, \texttt{gcg\-getDualbound()},
\texttt{gcg\-getPrimalbound()}, and \texttt{gcg\-getGap()} were added
to the public interface and are called from the dialog handlers as
well as the Python interface. As a side effect, \gcg can now be used
better as a C/C++ shared library.

\subsection{Visualization Suite}
\label{sec:gcg_visu_suite}



Visualizations of algorithmic behavior can yield understanding and
intuition for interesting parts of a solving process.
With \gcg~\gcgversion, we include a
\emph{visualization suite} that offers different
visualization scripts to show processes and results related to
detection, branching, or pricing, among others. These scripts are written in Python
3 and included in the folder \texttt{stats} and use the \texttt{.out},
\texttt{.res} and \texttt{.vbc} files generated when executing
\texttt{make test STATISTICS=true} (possible additional requirements
are given in the documentation). Furthermore, 
the suite also allows for two
additional ways of accessing the visualization scripts:

\begin{enumerate}
\item Reporting functionality: With two different scripts, callable
  via \texttt{make visu}, users can easily generate reports similar to
  the \emph{decomposition report} that was already available in \gcg~3.0, which
  offers an overview over all decompositions that \gcg found during its
  detection process. The generated documents include all
  visualizations offered by the suite along with descriptions of them
  in the captions. While the \emph{testset report} shows information
  about a single run of one selected testset, the \emph{comparison
    report} also compares two or more runs. Examples of both reports
  can now be found in the \gcg website documentation, see
  Section~\ref{sec:gcg_web_docu}.
\item Jupyter notebook: Since the scripts themselves already require
  a working installation of Python 3, we now added a
  \emph{visualization notebook} with which one can read data (sample
  data provided in the \gcg website documentation), clean, and filter it
  interactively, and visualize the results afterwards. The scripts of
  the visualization suite are imported and returned plots can be
  shown, exported, and even further edited.
\end{enumerate}

Just like \gcg should facilitate experimenting with a decomposition
approach without having to implement it, the visualization suite
should facilitate producing and presenting computational results and
algorithmic behavior. Also this is an ongoing long term effort.


\begin{figure}
    \centering
    \includegraphics[width=\textwidth]{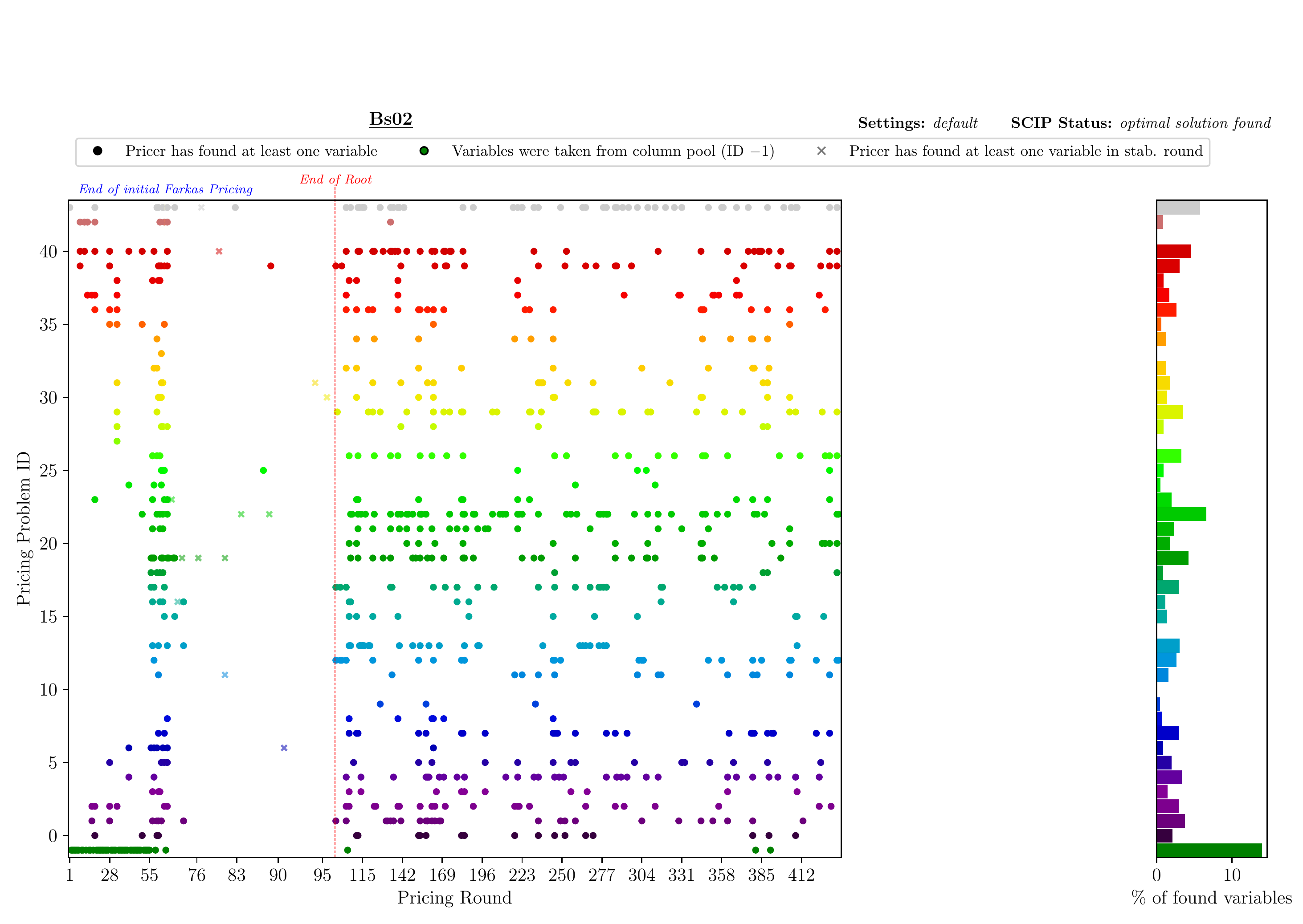}
    \caption{Bubble plot visualizing how the pricing problems
      performed during \gcg's Branch-and-Price process. This
      visualization was automatically generated using the new
      comparison report functionality.}
    \label{fig:pricing_bubbleplot}
\end{figure}
\begin{figure}
    \centering
    \includegraphics[width=\textwidth]{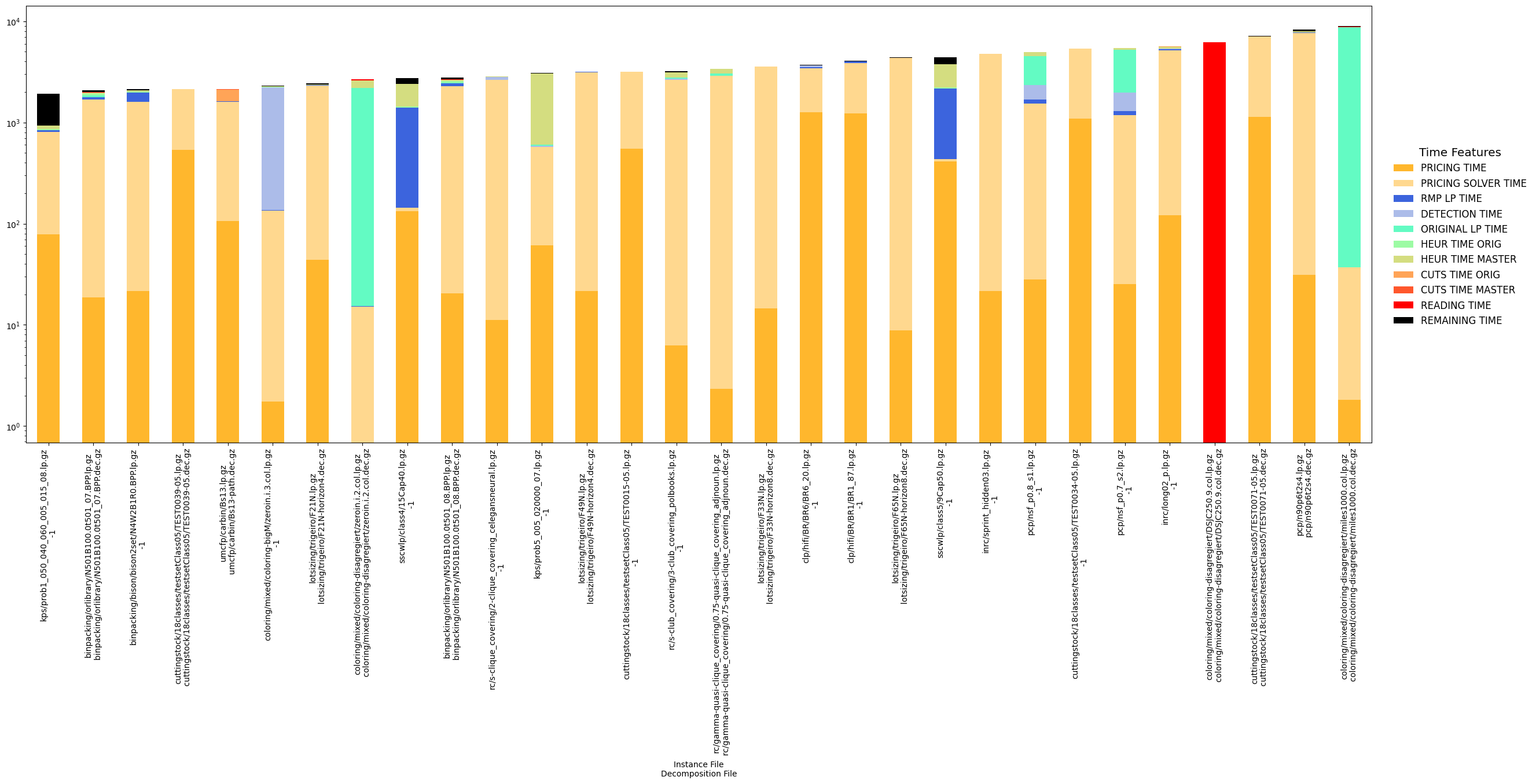}
    \caption{Time Distribution of a set of randomly drawn samples from
      the strIPlib. This visualization was generated using the
      included Jupyter notebook.}
    \label{fig:time_distribution}
\end{figure}

\subsection{Website Documentation}
\label{sec:gcg_web_docu}

The online documentation of \gcg was lagging behind the progress made
with the code itself. As part of this release, we offer a user-group
targeted website documentation. It enables users to make themselves
familiar with \gcg by means of very accessible feature descriptions for
functionality such as the explore menu or the visualization suite and
by a set of use cases to follow and reproduce.  For
developers, we now include a guide explaining the peculiarities of the
interplay between \gcg and \scip (``Getting Started: Developer's
Edition''). Within the ``Developer's Guide'', descriptions of existing
code and algorithmics such as detection, branching, and pricing are
provided to allow developers to familiarize themselves with them, if
required. Updates to the ``How to use'' (for instance, conducting experiments)
and ``How to add'' (for instance, adding branching rules) sections
completes the documentation.

\section{SCIP-SDP}
\label{sect:SCIP-SDP}


\scipsdp is a framework for solving mixed-integer semidefinite programs of
the following form
\begin{equation}\label{MISDP}
  \begin{aligned}
    \inf \quad & b^\T y \\
    \text{s.t.} \;\;\, & \sum_{k=1}^m A^k\, y_k - A^0 \succeq 0, \\
    & \ell_i \leq y_i \leq u_i && \forall\, i \in [m], \\
    & y_i \in \Z && \forall\, i \in I,
  \end{aligned}
\end{equation}
with symmetric matrices $A^k \in \R^{n \times n}$ for
$i \in \{0, \dots, m\}$, $b \in \R^m$, $\ell_i \in \R \cup \{- \infty\}$,
$u_i \in \R \cup \{\infty\}$ for all $i \in [m] \coloneqq \{1, \dots,
m\}$. The set of indices of integer variables is given by $I \subseteq [m]$
and $M \succeq 0$ denotes that a matrix $M$ is positive
semidefinite.

\scipsdp was initiated by Sonja Mars and Lars Schewe, see Mars~\cite{Mar13},
and then continued by Gally et al.~\cite{GallyPfetschUlbrich2018} and Gally~\cite{Gal19}. It
features interfaces to the SDP-solvers DSDP, Mosek, and SDPA.  In the
following, we briefly report on the changes since the last version 3.2.0.

\scipsdpv contains about \num{50000} lines of C-code. Since the last
version most of these lines have been touched. In particular, the interface
to the SDP-solvers has been completely revised. One benefit is that the
memory footprint of \scipsdp is now smaller for large instances. Moreover,
many bugs have been fixed.

Two important parameter changes that impact performance are:
\begin{itemize}
\item By default the number of used threads is 1 (it was previously set to
  ``automatic''). This change speeds-up the solution process by about 40\,\% for
  most smaller to medium sized SDPs.
\item The feasibility and optimality tolerances have been set to
  \num{e-5}. The exception is Mosek for which is is set to \num{e-6},
  because this leads to more reliable results.
\end{itemize}

\noindent
Further changes in a nutshell are the following:
\begin{itemize}
\item If the SDP-relaxation only has a single variable, it is solved using
  a semismooth Newton method. This slightly speeds up solution times and
  significantly decreases the times for heuristics. In particular, this
  holds for rounding heuristics on instances in which all integral
  variables are fixed except a single continuous variable. This continuous
  variable is often used for expressing the objective function, for example
  in cardinality least squares problems.
\item The LP-rows that are added to the LP-relaxation are strengthened
  using standard LP-preprocessing routines (coefficient tightening).
\item A new heuristic \texttt{heur\_fracround} has been added which
  iteratively rounds integer variables based on their fractional values in
  the last SDP-relaxation. In between, it performs propagation and solves a
  final SDP if unfixed continuous variables remain. This heuristic helps to
  significantly improve the overall running times. Propagation is now also
  used by the heuristic \texttt{heur\_sdprand} to improve its success rate
  for instances with additional linear constraints. Furthermore, both
  heuristics are correcting nearly integral values of integral variables in
  order to avoid small rounding errors, which might add up to significant
  amounts.
\item Several new presolving techniques have been introduced, which are
  discussed and evaluated in detail by Matter and Pfetsch~\cite{MatP21}. This includes two
  propagation methods to fix variables based on $2 \times 2$-minors and the
  upper bounds of other variables.
\item \scipsdp also allows to use LP-solving instead of SDP-relaxations
  using the parameter \texttt{misc/solvesdps}. It then generates so-called
  eigenvector cuts. The behavior of these cuts has been changed as
  follows. One can now add eigenvector cuts for all negative eigenvalues of
  the current infeasible relaxation. Moreover, SDP-relaxations can be
  solved in enforcing, that is, after all integer variables have integral
  values. Furthermore, the cuts can be sparsified.
\item The display of \scipsdp now changes depending on whether SDPs or LPs
  are solved for the relaxations. Moreover, the default settings are
  redefined for solving SDPs.
\item One can also generate a second-order cone relaxation, but so far this
  has not shown a run time improvement.
\item The readers for the SDPA and CBF formats have been completely revised
  (and rewritten for SDPA). They are now much faster and produce more
  warnings if errors occur.
\item \scipsdp can also handle rank-1 constraints, that is, the requirement
  that the resulting matrix has rank 1. This is achieved by adding
  quadratic constraints for $2 \times 2$-minors. Rank-1 constraints
  regularly appear in the literature, but are usually very hard to
  solve. The handling of these constraints has been revised.
\item The locking information (capturing whether the matrices $A^k$ are
  positive/negative semidefinite) is now copied to sub-SCIPs.
\item The statistics for solving SDP-relaxations has been extended and now
  reports more details.
\item There is a new file \texttt{scipsdpdef.h} that contains defines for
  the \scipsdp version. This enables code to depend on different \scipsdp
  versions.
\item It is now possible to add SDP-constraints within the solving process.
\item \scipsdp can now run concurrently, for example, by writing
  \texttt{concurrentopt} in the command line if \scip and \scipsdp are
  compiled using the \texttt{TPI}.
\item The updated Matlab Interface presented in
  Section~\ref{sect:interfaces} also allows to use SCIP-SDP.
\end{itemize}

Before we present some computational results, let us add some words of
caution. Although \scipsdp is numerically quite robust, accurately solving
SDPs is more demanding than solving LPs. This can lead to wrong results on
some instances\footnote{For instance, in seldom cases, the dual bound might
  exceed the value of a primal feasible solution.} and the results often
depend on the chosen tolerances. Technical reasons are that the SDPs are
solved using interior point solvers, which produce solutions with more
``numerical noise'' (since they do not have nonbasic variables). Moreover,
the solvers use relative tolerances, while \scipsdp uses absolute
tolerances. Finally, for Mosek, we use a slightly tighter tolerance than in
\scipsdp.

Table~\ref{tab:scipsdp40vs32} shows a comparison between \scipsdp 3.2 and
4.0 on the same testset as used by Gally et al.~\cite{GallyPfetschUlbrich2018}, which consists of 194 instances.
Reported are the number of optimally solved instances, as well as the
shifted geometric means of the number of processed nodes and the CPU time
in seconds. We use Mosek 9.2.40 for solving the continuous
SDP-relaxations. The tests were performed on a Linux cluster with 3.5 GHz
Intel Xeon E5-1620 Quad-Core CPUs, having 32~GB main memory and 10~MB
cache. All computations were run single-threaded and with a timelimit of
one hour.

As can be seen from the results, \scipsdp 4.0 is significantly faster than
\scipsdp 3.2, but we recall that we have relaxed the the tolerances (see
above). Nevertheless, the conclusion is that \scipsdp 4.0 has significantly
improved since the last version.

\begin{table}
  \caption{Performance comparison of \scipsdp 4.0 vs. \scipsdp 3.2}
  \label{tab:scipsdp40vs32}
  \small
  \begin{tabular*}{\textwidth}{@{}l@{\;\;\extracolsep{\fill}}ccc@{}}
    \toprule
    &  \# opt & \# nodes & time [s] \\
    \midrule
    \scipsdp 3.2 & 185 & 617.3 & 42.9 \\
    \scipsdp 4.0 & 187 & 497.3 & 26.6 \\
    \bottomrule
  \end{tabular*}
\end{table}

\section{\scipjack: Solving Steiner Tree and Related Problems}
\label{sect:SCIP-Jack}


Given an undirected, connected graph $G=(V,E)$, costs (or weights) $c: E \rightarrow
\R_{+}$ and a set $T \subseteq V$ of \textit{terminals}, the \textit{Steiner tree problem in graphs} (SPG)
asks for a tree $S = (V(S), E(S)) \subseteq G$ such that $T \subseteq V(S)$ holds and $\sum_{e \in E(S)} c(e)$ is minimized.
The SPG is a fundamental $\mathcal{NP}$-hard problem~\cite{Karp72}, and one of the most studied problems in combinatorial optimization. Moreover, many related problems have been extensively described in the literature and can be found in a wide range of practical applications~\cite{ljubic20}.

Since version 3.2, the \scipopt has contained \scipjack, an exact solver not only for the SPG but also for 11 related problems.  This release of the \scipopt contains the new \scipjack \scipjackversion\footnote{see also \url{http://scipjack.zib.de}}, which can handle two additional problem classes: the maximum-weight connected subgraph problem with budgets, and the partial-terminal Steiner tree problem.
Furthermore, \scipjack \scipjackversion~comes with major improvements on almost all problem classes it can handle.
Most importantly, the latest \scipjack outperforms the well-known SPG solver by Polzin and Vahdati~\cite{PolzinThesis, VahdatiThesis} on almost all nontrivial benchmark testsets from the literature. See the preprint~\cite{RehfeldtKoch20Conflicts} for more details. Notably, Polzin and Vahdati~\cite{PolzinThesis, VahdatiThesis} had remained out of reach for almost 20 years for any other SPG solver.

 The large number of newly implemented algorithms (and data structures) also results in an increase of the \scipjack code base by a factor of almost 3: To roughly 110\,000 lines of code. Additionally, the implementation of many existing methods has been improved.
We list several of the most important new features.

For SPG, a central new feature is a distance concept that provably dominates the well-known bottleneck Steiner distance from Duin and Volgenant~\cite{Duin1989}, see Rehfeldt and Koch~\cite{RehfeldtKoch21_ipco} for details. This distance concept is used in several (new) reduction methods implemented in \scipjack \scipjackversion. 
Also, the new \scipjack includes a full-fledged implementation of so-called \emph{extended reduction techniques}.
These methods are provably stronger than the state-of-the-art implementation~\cite{Polzin02}, and also yield strong practical results, see Rehfeldt and Koch~\cite{RehfeldtKoch20Conflicts} for details.
Furthermore, decomposition methods for the SPG have been implemented, for example to exploit the existence of biconnected components in the underlying graph.
Also, dynamic programming algorithms have been implemented to efficiently solve subproblems with special structures that sometimes arise after decomposition.

The improvements for SPG also have an immediate impact on problems that are transformed to SPG within \scipjack, such as the group Steiner tree problem. Even for the Euclidean Steiner tree problem, large improvements are possible (the SPG can be used for full Steiner tree concatenation after the discretization of the problem): \scipjack \scipjackversion~is able to solve 19 Euclidean Steiner tree problems with up to 100\,000 terminals for the first time to optimality, see Rehfeldt and Koch~\cite{RehfeldtKoch20Conflicts}. Notably, the state-of-the-art Euclidean Steiner tree solver GeoSteiner 5.1~\cite{juhl2018geosteiner} could not solve any of these instances even after one week of computation. In contrast, \scipjack \scipjackversion~solves all of them within 12 minutes, some even within two minutes.

Considerable problem-specific improvements have also been made for the prize-{\linebreak}collecting Steiner tree problem and (to a lesser extent) for the maximum-weight connect subgraph problem. 
For details on the improvements for the prize-collecting Steiner tree problem see Rehfeldt and Koch~\cite{RehfeldtKoch2021}, for the maximum-weight connect subgraph problem see Rehfeldt, Franz, and Koch~\cite{RehfeldtFranzKoch2020}.
The improvements encompass primal and dual heuristics as well as reduction techniques. As a result, \scipjack \scipjackversion~can solve many previously unsolved benchmark instances from both problem classes to optimality---the largest of these instances have up to 10 million edges. Additionally, for the prize-collecting Steiner tree problem \scipjack \scipjackversion~can solve most benchmark sets from the literature more than two times faster than its predecessor with respect to the shifted geometric mean (with a shift of 1 second).

\section{Final Remarks}
\label{sect:finalRemarks}

The SCIP Optimization Suite 8.0 release provides new functionality and improved performance and
reliability.
In SCIP, new symmetry handling features were added, including the handling of symmetries of general integer
and continuous variables, improving detection routines and adding a strategy for symmetry handling
routine selection.
Mixing cutting planes were implemented, which considerably improve the performance on chance constrained
programs.
A decomposition primal heuristic was updated to further improve the found solutions, and a new
decomposition primal heuristic was added.
A new cut strengthening procedure was added to the Benders decomposition framework, and
a new type of plugin for cut selection was introduced.

With this release also comes a thorough revision of how nonlinear constraints are handled in SCIP, in particular how extended formulations are created.
In the new version, the original formulation is preserved and the extended formulation is used for
relaxations only, which drastically improves the reliability of solutions.
High- and low-level nonlinear structures are now handled by plugins of different types in order
to avoid expression type ambiguity.
Simultaneously, a number of new MINLP features were introduced such as various new cutting planes,
symmetry detection for nonlinear constraints, support for sine and cosine functions, and others.


Regarding usability, the Julia package SCIP.jl was improved in several aspects and a new MATLAB interface to SCIP was
implemented.
UG was generalized to enable the parallelization of all solvers via a
unified framework, without the need to modify the framework for each solver; its internal
structure has been completely reworked.
The new version of GCG includes new algorithmic features and substantial ecosystem improvements,
such as extended interfaces, improved documentation, added utility for running and analyzing
computational experiments.
PaPILO features a new postsolving functionality for dual solutions when applied to pure LPs.
The handling of relaxations in SCIP-SDP was revised, and new heuristics and presolving methods were added.
The new version of SCIP-Jack can handle two additional problem classes and comes with major
performance improvements.

These developments yield a considerable performance improvement on nonconvex MINLP instances and
reduce the overall number of numerical failures on the MINLP testset, although a slowdown is
observed on convex instances due to a lack of recognition of a structure present in one instance
group. Nonetheless, we observe an overall runtime reduction of about $21\%$.
%
%
A substantial speed-up is observed on MILP instances with about a 50\% shorter runtime on the most challenging instances,
as well as special problem classes such as SDP and Steiner tree problems and their variants.
\subsection*{Acknowledgements}

The authors want to thank all previous developers and contributors to the \scip Optimization Suite and
all users that reported bugs and often also helped reproducing and fixing the bugs.
In particular, thanks go to Suresh Bolusani, Didier Chételat, Gregor Hendel, Andreas Schmitt,
Helena Völker, Robert Schwarz, Matthias Miltenberger, Matthias Walter, and Antoine Prouvoust and the Ecole team.
The Matlab-SCIP(-SDP) interface was set up with the big help of Nicolai Simon.

\subsection*{Contributions of the Authors}

The material presented in the article is highly related to code and software.
In the following we try to make the corresponding contributions of the authors and possible contact points more transparent.

JvD, CH, and MP are responsible for the changes of the symmetry handling routines (Section~\ref{sec:symmetry}).
The extension of symmetry handling to nonlinear constraints (Section~\ref{sect:symmetrynl}) is by FW.
WC and MP implemented the mixing cut separator (Section \ref{subsect:mixcuts}).
The update of \method{PADM} and the new \method{DPS} heuristic (Section~\ref{subsect:heuristics}) are due to KH and DW.
SJM is responsible for the updates to the Benders' decomposition framework (Section~\ref{subsect:benders}).
MT and FeS implemented the new cut selector plugin (Section~\ref{subsect:cutselectors}).
Various technical improvements (Section~\ref{subsect:further}) were added by MP and SV.
The new expressions framework (Section~\ref{sect:expr}) is by BM, FeS, FW, KB, and SV.
The rewritten handler for nonlinear constraints (Section~\ref{sect:consnl}) is by BM, FeS, KB, and SV.
The nonlinear handler for quadratic expressions (Section~\ref{sect:nlhdlrquad}) and the separator \texttt{sepa\_{\allowbreak}interminor} (Section~\ref{sect:interminor}) are by AC and FeS.
The nonlinear handler for second-order cones (Section~\ref{sect:nlhdlrsoc}) is by BM, FeS, and FW.
The nonlinear handler for bilinear expressions (Section~\ref{sect:nlhdlrbilin}) and the separator \texttt{sepa\_minor} (Section~\ref{sect:minor}) are by BM.
The nonlinear handler for convex and concave expressions (Section~\ref{sect:nlhdlrconvex}) are by BM, KB, and SV.
The nonlinear handler for quotients (Section~\ref{sect:nlhdlrquotient}) is by BM and FW.
The nonlinear handler for perspective reformulations (Section~\ref{sect:nlhdlrpersp}) is by KB.
The separator for RLT cuts (Section~\ref{sect:rlt}) is by FW and KB.
The separator for principal minors of $X\succeq xx^\T$ (Section~\ref{sect:minor}) is by BM and FW.
The separator for intersection cuts on the rank-1 constraint for a matrix (Section~\ref{sect:interminor}) is by AC and FeS.
The revised primal heuristic \texttt{subnlp} (Section~\ref{sect:subnlp}) and the updates to NLP, NLPI, and AD interfaces (Section~\ref{sect:nlp}) are by SV.

The changes to SoPlex (Section~\ref{sect:soplex}) are due to LE and AH.
AH and AG are responsible for the new dual postsolving functionality in PaPILO (Section~\ref{sect:papilo}).
The new AMPL interface of SCIP (Section~\ref{sect:ampl}) was implemented by SV.
The Julia interface \texttt{SCIP.jl} (Section~\ref{subsect:julia}) was extended and updated by MB, Erik Tadewaldt, Robert Schwarz, and Yupei Qi.
The \soplex C interface (Section~\ref{subsect:cwrapper}) was developed by AC, LE, and MB.
Nicolai Simon and MP updated the Matlab-interface (Section~\ref{subsect:matlab}).
The work on ZIMPL (Section~\ref{sect:zimpl}) was done by TK.
The updates to the UG framework (Section~\ref{sect:ug}) are by YS.
Concerning \gcg (Section~\ref{sect:gcg}), EM refactored the detector loop; the website documentation and
visualization suite is due to TD; OG created the strong branching
code; and SS implemented \textsc{PyGCGOpt}.
FM and MP implemented the changes in SCIP-SDP (Section~\ref{sect:SCIP-SDP}).
DR is responsible for \scipjack (Section~\ref{sect:SCIP-Jack}).

The work by FrS on the continuous integration system, regular test and benchmark runs, binary distributions, websites, and many more has been invaluable for all developments.

\renewcommand{\refname}{\normalsize References}
\setlength{\bibsep}{0.25ex plus 0.3ex}
\bibliographystyle{abbrvnat}

\begin{small}
\bibliography{scipopt}

\end{small}

\subsection*{Author Affiliations}

\hypersetup{urlcolor=black}
\newcommand{\myorcid}[1]{ORCID: \href{https://orcid.org/#1}{#1}}
\newcommand{\myemail}[1]{E-mail: \href{#1}{#1}}
\newcommand{\myaffil}[2]{{\noindent #1}\\{#2}\bigskip}

\small

\myaffil{Ksenia Bestuzheva}{
  Zuse Institute Berlin, Department AIS$^2$T, Takustr.~7, 14195~Berlin, Germany\\
  \myemail{bestuzheva@zib.de}\\
  \myorcid{0000-0002-7018-7099}}

\myaffil{Mathieu Besançon}{%
  Zuse Institute Berlin, Department AIS$^2$T, Takustr.~7, 14195~Berlin, Germany\\
  \myemail{besancon@zib.de}\\
  \protect\myorcid{0000-0002-6284-3033}}

\myaffil{Wei-Kun Chen}{%
  School of Mathematics and Statistics, Beijing Institute of Technology, Beijing 100081, China\\
  \myemail{chenweikun@bit.edu.cn}\\
  \myorcid{0000-0003-4147-1346}}

\myaffil{Antonia Chmiela}{%
  Zuse Institute Berlin, Department AIS$^2$T, Takustr.~7, 14195~Berlin, Germany\\
  \myemail{chmiela@zib.de}\\
  \myorcid{0000-0002-4809-2958}}

\myaffil{Tim Donkiewicz}{%
  RWTH Aachen University, Lehrstuhl f\"ur Operations Research, Kackertstr.~7, 52072~Aachen, Germany\\
  \myemail{tim.donkiewicz@rwth-aachen.de}\\
  \myorcid{0000-0002-5721-3563}}

\myaffil{Jasper van Doornmalen}{%
  Technische Universiteit Eindhoven, Department of Mathematics and Computer Science, P.O.\ Box 513, 5600 MB Eindhoven, The Netherlands\\
  \myemail{m.j.v.doornmalen@tue.nl}\\
  \myorcid{0000-0002-2494-0705}}

\myaffil{Leon Eifler}{%
  Zuse Institute Berlin, Department AIS$^2$T, Takustr.~7, 14195~Berlin, Germany\\
  \myemail{eifler@zib.de}\\
  \myorcid{0000-0003-0245-9344}}

\myaffil{Oliver Gaul}{%
  RWTH Aachen University, Lehrstuhl f\"ur Operations Research, Kackertstr.~7, 52072~Aachen, Germany\\
  \myemail{oliver.gaul@rwth-aachen.de}\\
  \myorcid{0000-0002-2131-1911}}

\myaffil{Gerald Gamrath}{%
  Zuse Institute Berlin, Department AIS$^2$T, Takustr.~7, 14195~Berlin, Germany\\
  and I$^2$DAMO GmbH, Englerallee 19, 14195~Berlin, Germany\\
  \myemail{gamrath@zib.de}\\
  \myorcid{0000-0001-6141-5937}}

\myaffil{Ambros Gleixner}{%
  Zuse Institute Berlin, Department AIS$^2$T, Takustr.~7, 14195~Berlin, Germany\\
  \myemail{gleixner@zib.de}\\
  \myorcid{0000-0003-0391-5903}}

\myaffil{Leona Gottwald}{%
  Zuse Institute Berlin, Department AIS$^2$T, Takustr.~7, 14195~Berlin, Germany\\
  \myemail{gottwald@zib.de}\\
  \myorcid{0000-0002-8894-5011}}

\myaffil{Christoph Graczyk}{%
  Zuse Institute Berlin, Department AIS$^2$T, Takustr.~7, 14195~Berlin, Germany\\
  \myemail{graczyk@zib.de}}

\myaffil{Katrin Halbig}{%
  Friedrich-Alexander Universität Erlangen-Nürnberg, Department of Data Science, Cauerstr.~11, 91058~Erlangen, Germany\\
  \myemail{katrin.halbig@fau.de}\\
  \myorcid{0000-0002-8730-3447}}

\myaffil{Alexander Hoen}{%
  Zuse Institute Berlin, Department AIS$^2$T, Takustr.~7, 14195~Berlin, Germany\\
  \myemail{hoen@zib.de}\\
  \myorcid{0000-0003-1065-1651}}

\myaffil{Christopher Hojny}{%
  Technische Universiteit Eindhoven, Department of Mathematics and Computer Science, P.O.\ Box 513, 5600 MB Eindhoven, The Netherlands\\
  \myemail{c.hojny@tue.nl}\\
  \myorcid{0000-0002-5324-8996}}

\myaffil{Rolf van der Hulst}{%
  University of Twente, Department of Discrete Mathematics and Mathematical Programming, P.O. Box 217, 7500 AE Enschede, The Netherlands\\
  \myemail{r.p.vanderhulst@utwente.nl}}

\myaffil{Thorsten Koch}{%
  Technische Universit\"at Berlin, Chair of Software and Algorithms for Discrete Optimization, Stra\ss{}e des 17. Juni 135, 10623 Berlin, Germany, and\\
  Zuse Institute Berlin, Department A$^2$IM, Takustr. 7, 14195~Berlin, Germany\\
  \myemail{koch@zib.de}\\
  \myorcid{0000-0002-1967-0077}}

\myaffil{Marco L\"ubbecke}{%
  RWTH Aachen University, Lehrstuhl f\"ur Operations Research, Kackertstr.~7, 52072~Aachen, Germany\\
  \myemail{marco.luebbecke@rwth-aachen.de}\\
  \myorcid{0000-0002-2635-0522}}

\myaffil{Stephen J.~Maher}{%
  University of Exeter, College of Engineering, Mathematics and Physical Sciences, Exeter, United Kingdom\\
  \myemail{s.j.maher@exeter.ac.uk}\\
  \myorcid{0000-0003-3773-6882}}

\myaffil{Frederic Matter}{%
  Technische Universität Darmstadt, Fachbereich Mathematik, Dolivostr.~15, 64293~Darmstadt, Germany\\
  \myemail{matter@mathematik.tu-darmstadt.de}\\
  \myorcid{0000-0002-0499-1820}}

\myaffil{Erik M\"uhmer}{%
  RWTH Aachen University, Lehrstuhl f\"ur Operations Research, Kackertstr.~7, 52072 Aachen\\
  \myemail{erik.muehmer@rwth-aachen.de}\\
  \myorcid{0000-0003-1114-3800}} 

\myaffil{Benjamin Müller}{%
  Zuse Institute Berlin, Department AIS$^2$T, Takustr.~7, 14195~Berlin, Germany\\
  \myemail{benjamin.mueller@zib.de}\\
  \myorcid{0000-0002-4463-2873}}

\myaffil{Marc E.~Pfetsch}{%
  Technische Universität Darmstadt, Fachbereich Mathematik, Dolivostr.~15, 64293~Darmstadt, Germany\\
  \myemail{pfetsch@mathematik.tu-darmstadt.de}\\
  \myorcid{0000-0002-0947-7193}}

\myaffil{Daniel Rehfeldt}{%
  Zuse Institute Berlin, Department A$^2$IM, Takustr.~7, 14195~Berlin, Germany\\
  \myemail{rehfeldt@zib.de}\\
  \myorcid{0000-0002-2877-074X}
  }

\myaffil{Steffan Schlein}{%
  RWTH Aachen University, Lehrstuhl f\"ur Operations Research, Kackertstr.~7, 52072~Aachen, Germany\\
  \myemail{steffan.schlein@rwth-aachen.de}}

\myaffil{Franziska Schlösser}{%
  Zuse Institute Berlin, Department AIS$^2$T, Takustr.~7, 14195~Berlin, Germany\\
  \myemail{schloesser@zib.de}}

\myaffil{Felipe Serrano}{%
  Zuse Institute Berlin, Department AIS$^2$T, Takustr.~7, 14195~Berlin, Germany\\
  \myemail{serrano@zib.de}\\
  \myorcid{0000-0002-7892-3951}}

\myaffil{Yuji Shinano}{%
  Zuse Institute Berlin, Department A$^2$IM, Takustr.~7, 14195~Berlin, Germany\\
  \myemail{shinano@zib.de}\\
  \myorcid{0000-0002-2902-882X}}

\myaffil{Boro Sofranac}{%
  Zuse Institute Berlin, Department AIS$^2$T, Takustr.~7, 14195~Berlin, Germany and\\
  Technische Universit\"at Berlin, Stra\ss{}e des 17. Juni 135, 10623~Berlin, Germany\\
  \myemail{sofranac@zib.de} \\
  \myorcid{0000-0003-2252-9469}}

\myaffil{Mark Turner}{%
  Zuse Institute Berlin, Department A$^2$IM, Takustr.~7, 14195~Berlin, Germany\\
  and Chair of Software and Algorithms for Discrete Optimization, Institute of Mathematics, Technische Universität Berlin, Straße des 17. Juni 135, 10623 Berlin, Germany \\
  \myemail{turner@zib.de}\\
  \myorcid{0000-0001-7270-1496}}

\myaffil{Stefan Vigerske}{%
  GAMS Software GmbH, c/o Zuse Institute Berlin, Department AIS$^2$T, Takustr.~7, 14195~Berlin, Germany\\
  \myemail{svigerske@gams.com}}

\myaffil{Fabian Wegscheider}{%
  Zuse Institute Berlin, Department AIS$^2$T, Takustr.~7, 14195~Berlin, Germany\\
  \myemail{wegscheider@zib.de}}

\myaffil{Philipp Wellner}{%
  \myemail{p.we@fu-berlin.de}
  }

\myaffil{Dieter Weninger}{%
  Friedrich-Alexander Universität Erlangen-Nürnberg, Department of Data Science, Cauerstr.~11, 91058~Erlangen, Germany\\
  \myemail{dieter.weninger@fau.de}\\
  \myorcid{0000-0002-1333-8591}}

\myaffil{Jakob Witzig}{%
  Zuse Institute Berlin, Department AIS$^2$T, Takustr.~7, 14195~Berlin, Germany\\
  \myemail{witzig@zib.de}\\
  \myorcid{0000-0003-2698-0767}}



\clearpage
\newgeometry{textwidth=.8\paperwidth}
\begin{appendices}




\section{Detailed Computational Results to Section~\ref{sect:perfconsexpr} (Performance Impact of Updates for Nonlinear Constraints)}
\label{appendixsect:consexpr}

The following table lists the results for running the classic code (SCIP 7) and the new code (SCIP 8) for each considered instance of MINLPLib.
Only results for the non-permuted instances are given.

Column ``time/gap'' gives the time it took to solve the instance to optimality (with respect to specified gap tolerances)
or the gap at termination if solving stopped at the time limit.
If the time or gap of a version is not worse than the other version, the time or gap is printed in a bold font.
If a version did not return a result or the reported bounds conflict with best known bounds, then ``fail'' is printed.

For the classic version, columns ``quad'', ``soc'', ``abspow'', and ``nonlin'' give the number of quadratic, second-order cone, abspower, and nonlinear constraints, respectively, after presolve.
Due to the reformulations that are applied in presolve, nonlinear constraints are sums of convex or concave functions, quadratic terms (including bilinear products) are parts of quadratic constraint, unless an upgrade to a soc constraint was possible.
Monomials of odd degree, signpowers, and monomials of even degree with fixed sign are handled by the abspower constraint handler.

For the new version, columns ``quad'', ``bilin'', ``soc'', ``convex'', ``concave'', ``quot'', ``persp'', ``def'' give the number of expressions for which the detection algorithm of the nonlinear handlers quadratic, bilinear, soc, convex, concave, quotient, perspective, and default, respectively (see Sections~\ref{sect:nlhdlrquad}--\ref{sect:nlhdlrpersp} and~\ref{sect:nlhdlrdefault}) reported success, that is, registered themselve for domain propagation or separation after presolve.
Recall that by default the quadratic nonlinear handler only gets active for propagable quadratic expressions (see Section~\ref{sect:nlhdlrquaddetect}) and the convex and concave nonlinear handlers only handle nontrivial expressions (Section~\ref{sect:nlhdlrconvexdetect}).
Further, the nonlinear handler for bilinear expressions (Section~\ref{sect:nlhdlrbilin}) currently registers itself for any product of two non-binary variables (original or auxiliary) and only checks when called later whether linear inequalities in the two variables are available, because the latter are computed after the extended formulations are initialized.
Columns ``minor'' and ``rlt'' report the number of cuts that were generated by the respective separators (Sections~\ref{sect:minor} and~\ref{sect:rlt}) and got added to the LP.
Here, if cuts were generated but not applied, a zero is printed.

The last row summarizes on how many instances each constraint handler, nonlinear handler, or separator was used, i.e., the number of nonzeros in each column.

{%
\tiny
\setlength{\extrarowheight}{0.5pt}
\begin{longtable}{@{}l|@{\,}r@{\,}r@{\,}r@{\,}r@{\,}r@{\,}|@{\,}r@{\,}r@{\,}r@{\,}r@{\,}r@{\,}r@{\,}r@{\,}r@{\,}r@{\,}r@{\,}r@{}}
\toprule
& \multicolumn{5}{c|@{\,}}{classic} & \multicolumn{11}{c}{new} \\
instance & time/gap & quad&soc&abspow&nonlin & time/gap & quad&bilin&soc&convex&concave&quot&def&persp & minor&rlt \\ \midrule
\endhead
4stufen & 61.4\% & 26 &  &  & 16 & \textbf{14.6\%} & 11 & 26 &  & 5 & 5 & 12 & 153 &  &  &  \\
abel & 0.8s & 1 &  &  &  & \textbf{0.0s} &  &  &  &  &  &  & 58 &  &  &  \\
alan & \textbf{0.1s} & 1 &  &  &  & \textbf{0.1s} & 1 &  &  & 1 &  &  & 4 &  &  &  \\
alkyl & \textbf{0.2s} & 9 &  &  &  & 0.3s & 3 & 8 &  &  &  &  & 32 &  &  &  \\
alkylation & 1.2\% & 7 &  & 1 &  & \textbf{0.13\%} & 5 & 7 &  & 1 &  & 3 & 28 &  &  &  \\
arki0001 & \textbf{0.1s} & 1 &  &  &  & 0.3s &  &  &  &  &  &  & 1026 &  &  &  \\
arki0002 & \textbf{$\infty$} & 1 &  &  & 1824 & \textbf{$\infty$} &  &  &  & 1824 & 912 &  & 6258 &  &  &  \\
arki0003 & 7.5s & 4669 &  & 360 &  & \textbf{2.1s} & 1080 & 77 &  &  &  &  & 6243 &  &  &  \\
\bottomrule
\end{longtable}

}

\clearpage
\section{Detailed Computational Results to Section~\ref{sect:gcgstrongbranch} on Strong
  Branching in GCG}
\label{appendsect:gcg}

\newcommand\NameEntry[1]{%
  \multirow{1}*{%
    \begin{varwidth}{42em}%
    \flushleft #1%
    \end{varwidth}}}

\newcommand{\grrrr}{\raisebox{0.25ex}{\tiny $>$}}
\newcommand{\spc}{\hspace{2em}}

{\makeatletter
\setlength{\@fptop}{0pt}
\setlength{\@fpsep}{\floatsep}
\setlength{\@fpbot}{0pt plus 1fil}
\makeatother
\setlength{\extrarowheight}{0.6pt}
\setlength{\tabcolsep}{2pt}

The following table gives the time in seconds needed to solve each problem instance with original variable branching. Entries that perform better than \texttt{pseudocost} are \textit{italic}, the best entry in each row is in \textbf{bold face}.

  \tiny
  \noindent
\begin{tabular*}{\textwidth}{@{\extracolsep{\fill}}lrrrrrrrrrr@{}}
      \toprule
                                          & \NameEntry{\texttt{pseudo-\\cost}}&\texttt{random}&\NameEntry{\texttt{most-\\frac}}&\texttt{SBw/oCG}&\texttt{SBw/CG}&\NameEntry{\texttt{hierar-\\chical}} & \texttt{hybrid}&\NameEntry{\texttt{reli-\\able}}&\NameEntry{\texttt{rel.\\hier.}}&\NameEntry{\texttt{hybrid\\hier.}}\\[2ex]
      \midrule
            12Cap10                             &              237.7  &             3433.9  &              841.7  &             2071.0  &             1245.0  & \textit{     223.7} & 635.1                & 3197.4               & \textbf{126.9}       & \textit{     174.8} \\
      14Cap10                             &              100.9  & \textit{      60.5} &              176.4  &              452.1  &              401.3  & \textbf{      59.3} & 213.9                & \textit{60.2}        & 135.7                & \textit{      93.2} \\
      20Cap10                             &              171.2  &              613.6  &              402.4  &             1569.5  &             1043.5  & \textit{     117.0} & 458.3                & 625.8                & \textit{129.5}       & \textbf{      87.7} \\
      NU\_1\_0010\_05\_3                  & \textbf{      40.9} &               41.2  &               41.7  &               44.5  &               48.5  &               47.2  & 47.5                 & 46.7                 & 41.4                 &               45.1  \\
      NU\_1\_0010\_05\_7                  & \textbf{      23.5} &               49.4  &               23.7  &               36.7  &               30.4  &              181.5  & 27.1                 & 27.6                 & 40.4                 &               32.4  \\
\bottomrule
    \end{tabular*}%
}

\setcounter{totalnumber}{1}


\clearpage
{\makeatletter
\setlength{\@fptop}{0pt}
\setlength{\@fpsep}{\floatsep}
\setlength{\@fpbot}{0pt plus 1fil}
\makeatother
\setlength{\extrarowheight}{0.6pt}
\setlength{\tabcolsep}{2pt}
  The following table gives the number of nodes needed to solve each problem instance with original variable branching. The best entry in each row except full strong entries is in \textbf{bold face}.
  
  \tiny
  \noindent
  \centering
  \begin{tabular*}{\textwidth}{@{\extracolsep{\fill}}lrrrrrrrrrr@{}}
      \toprule
      & \NameEntry{\texttt{pseudo-\\cost}}&\texttt{random}&\NameEntry{\texttt{most-\\frac}}&\texttt{SBw/oCG}&\texttt{SBw/CG}&\NameEntry{\texttt{hierar-\\chical}} & \texttt{hybrid}&\NameEntry{\texttt{reli-\\able}}&\NameEntry{\texttt{rel.\\hier.}}&\NameEntry{\texttt{hybrid\\hier.}}\\[2ex]
      \midrule
      12Cap10                             &             1\,009  &             7\,083  &             2\,234  &                351  &                137  &                186  &         143          & 7\,083                 &   \textbf{137}         &                280  \\
      14Cap10                             &                643  &                422  &             1\,056  &                102  &                 52  &                 68  & \textbf{47}          &           422          & 643                    &                128  \\
      20Cap10                             &                795  &             2\,575  &             1\,569  &                318  &                112  &                116  & \textbf{95}          & 2\,575                 & 795                    &                135  \\
      NU\_1\_0010\_05\_3                  &             4\,672  &             4\,439  &             4\,672  &             4\,352  &             4\,397  &             4\,397  &         4\,397       &           4\,397       &           4\,439       & \textbf{    4\,284} \\
      NU\_1\_0010\_05\_7                  & \textbf{    2\,495} &             4\,576  & \textbf{    2\,495} &             3\,713  &             2\,643  &            12\,338  & 2\,652               & 2\,643                 & 3\,730                 &             3\,278  \\
      NU\_3\_0010\_05\_1                  &             1\,368  &             4\,229  &             1\,490  &             1\,958  &                415  & \textbf{       566} &         800          & 4\,229                 &           920          &             2\,253  \\
\bottomrule
    \end{tabular*}
}

\end{appendices}

\end{document}